\documentclass[final,5p,times,sort&compress]{elsarticle}
\usepackage{amsmath, amsfonts, amssymb}
\usepackage{graphicx,epsfig}
\usepackage{color}
\usepackage{booktabs}
\usepackage[section]{placeins}
\numberwithin{equation}{section}

\usepackage{array}
\usepackage{enumitem}
\usepackage{dsfont}
\graphicspath{{images1/}}

\usepackage{amsthm}
\theoremstyle{plain}\newtheorem{remark}{Remark}[section]
\theoremstyle{definition}\newtheorem{definition}{Definition}[section]

\newcommand{\wt}{\widetilde}
\DeclareMathOperator{\argmax}{argmax}
\DeclareMathOperator{\re}{Re}

\DeclareMathOperator{\std}{std}

\definecolor{shadegray}{rgb}{0.8,0.8,0.8}

\journal{Digital Signal Processing}
\usepackage{hyperref}

\begin{document}

\begin{frontmatter}

\title{Linear and synchrosqueezed time-frequency representations revisited. \\Part I: Overview, standards of use, related issues and algorithms.}
\tnotetext[support]{Work supported by the Engineering and Physical Sciences Research Council, UK}

\author{Dmytro Iatsenko}
\ead{dmytro.iatsenko@gmail.com}

\author{Peter V. E. McClintock\corref{}}
\ead{p.v.e.mcclintock@lancaster.ac.uk}

\author{Aneta Stefanovska\corref{cor1}}
\cortext[cor1]{Corresponding author}
\ead{aneta@lancaster.ac.uk}

\address{Department of Physics, Lancaster University, Lancaster LA1 4YB, UK}


\begin{abstract}

Time-frequency representations (TFRs) of signals, such as the windowed Fourier transform (WFT), wavelet transform (WT) and their synchrosqueezed variants (SWFT, SWT), provide powerful analysis tools. However, there are many important issues related to the practical use of TFRs that need to be clarified. Here we present a thorough review of these TFRs, summarizing all theoretical, practical and numerical aspects of their use, reconsidering some conventions and introducing new concepts and procedures. The purposes of this work are: (i) to provide a consistent overview of the computation, properties, and use of the (S)WFT/(S)WT methods; (ii) to establish general standards related to their use, both theoretical and practical; and (iii) to provide clean and optimized algorithms and MatLab codes, appropriate for any window or wavelet.

\end{abstract}


\end{frontmatter}


\section{Introduction}\label{sec:introduction}

Time-frequency analysis is devoted to the problem of identifying and quantifying the oscillatory components present in a signal. This is done by projecting the signal onto the time-frequency plane, which allows one to study its properties simultaneously in time and frequency. Such projections are called time-frequency representations (TFRs) and there are a number of different kinds, depending on the exact way in which the projection is carried out. Time-frequency analysis is especially useful for signals containing many oscillatory components with time-varying amplitudes and/or frequencies, which is a very common scenario for real-life signals. The paradigmatic example is the heart beat: although the cardiac frequency is usually localized around $\sim$1\,Hz, it varies continuously about this average; these variations are hard to analyze in the time domain, or by using the Fourier transform (FT), but they can easily be traced in the time-frequency plane \cite{Delprat:92,Carmona:97,Carmona:99,Iatsenko:cardio,Thakur:13,Daubechies:11}. TFRs have thus become established as very powerful tools that are routinely applied in almost every area of science, from image processing and finance to geophysics and the life sciences \cite{Mallat:08,Addison:10,Iatsenko:cardio,Shiogai:10,Wacker:13,Stefanovska:99b,Torrence:98,Kumar:97,Stefanovska:07,Labat:05a,Labat:05b,Aguiar:12}.

There exist two main types of TFR: linear and quadratic. Although each of these has its own advantages, in what follows we will consider only linear TFRs, such as the windowed Fourier transform (WFT) and wavelet transform (WT). They are to be preferred because, first, they are more readily interpretable on account of being additive (the TFR of a sum of signals = the sum of the TFRs for each signal) and, secondly, linear TFRs offer the possibility of extracting and reconstructing individual components, which can be problematic for quadratic representations. Additionally, we also consider the synchrosqueezed WFT and WT (abbreviated as SWFT and SWT) \cite{Daubechies:11,Thakur:11,Thakur:13,Auger:13,Daubechies:96}, which represent a particular nonlinear transformation of the original WFT and WT to increase their concentration. Although not additive, the SWFT/SWT still allows for straightforward extraction and reconstruction of the signal's components.

General aspects of time-frequency analysis, and the properties of the different existing types of TFR, have been thoroughly discussed in a number of excellent books \cite{Mallat:08,Stankovic:13,Hlawatsch:10,Boashash:03,Cohen:95,Addison:10,Daubechies:92,Grochenig:01,Kaiser:94} and reviews \cite{Hlawatsch:92,Cohen:89,Torrence:98,Qian:99}. In a sense, however, there is too much information available and some aspects, e.g.\ those related to the synchrosqueezed transforms, have still not been considered in sufficient detail. The purposes of the present work are therefore: (i) to bring together and summarize the minimal but sufficient knowledge needed to understand and apply the WFT and WT effectively; (ii) to cover some aspects of their use that are not discussed consistently in the literature; (iii) to consider in detail the synchrosqueezed transforms and their properties; (iv) to advance the theory and extend the applicability of the (S)WFT and (S)WT by introducing certain improvements, procedures and concepts; (v) to provide clean and user-friendly algorithms and MatLab codes, appropriate for any window or wavelet; and (vi) to investigating some important related issues.

The review falls into two parts. In the present Part I we provide a thorough discussion of the TFRs mentioned and their properties, concentrating on the most practically relevant aspects and not going too deep into the mathematical theory. In Part II we will deal with more sophisticated issues, such as the optimal choice of window/wavelet parameters, the performance of different reconstruction methods, and the advantages/drawbacks of synchrosqueezing. In both parts, particular emphasis is placed on the TFR-based estimation of the properties of the components present in the signal, and on the underlying caveats.

In Sec.\ \ref{sec:analyticsignal} below, we discuss the analytic signal approach and its applicability. Different TFRs are introduced and discussed in Sec.\ \ref{sec:tfr}, questions of time and frequency resolution are considered in Sec.\ \ref{sec:tfres}, and the extraction of individual components from the TFRs is discussed in Sec.\ \ref{sec:tfsrec}. Practical issues associated with signal preprocessing, frequency discretization, boundary effects and frequency limits are dealt with in Sec.\ \ref{sec:pract}. The step-by-step algorithms for each TFR and the corresponding (freely downloadable) MatLab codes are summarized in Sec.\ \ref{sec:stepalg}. \ref{app:nomenclature} summarizes the notation, abbreviations, terminology and conventions used throughout the paper to supplement the local introduction of the individual items as the text develops. \ref{app:asigerr}--\ref{app:winwav} fill in some technical details that are omitted from the main text to avoid interrupting the smooth flow of ideas. Finally, common window/wavelet types, their resolution properties and expressions for the related quantities are considered in \ref{app:winwav}.

\begin{remark}
It should be noted, that many advanced methods based on different TFRs have been developed over the last decades. Examples include wavelet bispectral analysis \cite{Jamsek:07,Jamsek:10}, harmonic identification \cite{Sheppard:11}, wavelet coherence \cite{Grinsted:04,Lachaux:02,Labat:05a} and phase coherence \cite{Bandrivskyy:04,Bloomfield:04,Sheppard:12}. However, it seems in principle impossible to review appropriately all such techniques in one work, given their large number, and we do not consider them all here. Rather, we concentrate on the in depth study of the basic time-frequency representations, which provide the foundation for all these advanced methods.
\end{remark}

\section{Analytic signal}\label{sec:analyticsignal}

One of the basic notions of time-frequency analysis is the \emph{AM/FM component} (or simply \emph{component}), which is defined as a function of time $t$ of form
\begin{equation}\label{amfm}
x(t)=A(t)\cos\phi(t)\;\;\;\big(\forall t:\;A(t)>0,\nu(t)\equiv\phi'(t)>0\big).
\end{equation}
The time-dependent values $A(t)$, $\phi(t)$ and $\nu(t)\equiv\phi'(t)$ are then called the instantaneous amplitude, phase and frequency of the component (\ref{amfm}) (for a more detailed discussion of their definitions and related issues see \cite{Picinbono:97,Boashash:92a,Boashash:92b}).

Given that the signal is known to be of the form (\ref{amfm}), the natural question is how to find its associated $A(t)$, $\phi(t)$ and $\nu(t)$. The most convenient way of doing this is the analytic signal approach. However, before considering it, a few additional notions should be introduced. Thus, for an arbitrary function $f(t)$, its Fourier transform (FT), positive and negative frequency parts, time-average and standard deviation will be denoted as $\hat{f}(\xi)$, $f^{+}(t)$, $f^{-}(t)$, $\langle f(t)\rangle$ and $\std[f(t)]$, respectively:
\begin{equation}\label{nt}
\begin{gathered}
\hat{f}(\xi)\equiv\int_{-\infty}^{\infty} f(t)e^{-i\xi t}dt\\
f(t)=\frac{1}{2\pi}\int_{-\infty}^{\infty}\hat{f}(\xi)e^{i\xi t}d\xi=\langle f(t)\rangle+f^+(t)+f^{-}(t),\\
f^{+}(t)\equiv \int_{0^+}^\infty \hat{f}(\xi)e^{i\xi t}d\xi,\quad
f^{-}(t)\equiv \int_{-\infty}^{0^-} \hat{f}(\xi)e^{i\xi t}d\xi,\\
\langle f(t)\rangle=\frac{\int f(t)dt}{\int dt},\quad
\std[f(t)]\equiv\sqrt{\langle [f(t)]^2\rangle-[\langle f(t)\rangle]^2},\\
\end{gathered}
\end{equation}
where, here and in what follows, the integrals are taken over $(-\infty,\infty)$ if unspecified (or, in practice, over the full time duration of $f(t)$). A simple example is $f(t)=a+b\cos\nu t$, for which one has $\hat{f}(\xi)=a\delta(\xi)+b\delta(\xi-\nu)e^{i\nu t}/2+b\delta(\xi+\nu)e^{-i\nu t}/2$, $\langle f(t)\rangle=a$, $f^{\pm}(t)=be^{\pm i\nu t}/2$ and $\std[f(t)]=b/\sqrt{2}$. Note, that if $f(t)$ is real, then $f(t)-\langle f(t)\rangle=2\re f^{\pm}(t)$ and $\hat{f}(\xi)=[\hat{f}(-\xi)]^*\Rightarrow f^{+}(t)=[f^{-}(t)]^*$, where the star denotes complex conjugation.

For a given signal $s(t)$ (which is always assumed to be real in this work), its doubled positive frequency part is called its \emph{analytic signal} and will be denoted as $s^{a}(t)$:
\begin{equation}\label{asig}
s^{a}(t)\equiv 2s^{+}(t)\quad\big(s(t)=\langle s(t)\rangle+{\rm Re}[s^{a}(t)]\big).
\end{equation}
The analytic signal is complex, so its dynamics can easily be separated into amplitude and phase parts. For signals represented by a single component (\ref{amfm}), the analytic amplitude and phase $A^{a}(t),\phi^{a}(t)$ match closely the true amplitude and phase $A(t),\phi(t)$, thus providing an easy way to estimate them:
\begin{equation}\label{ap}
A(t)\approx A^{a}(t)\equiv|s^{a}(t)|,\;\;\phi(t)\approx\phi^{a}(t)\equiv{\rm arg}[s^{a}(t)].\\
\end{equation}
The approximate equality (\ref{ap}) will be called the analytic approximation, and it can alternatively be formulated as $[A(t)\cos\phi(t)]^+\approx A(t)e^{i\phi(t)}/2$. As stipulated by the Bedrosian theorem \cite{Bedrosian:63}, this approximation is exact when the spectrum of $A(t)$ lies lower than the spectrum of $e^{i\phi(t)}$, and there are no intersections between the two. For example, in the case of amplitude modulation only, $s(t)=A(t)\cos (\nu t+\varphi)$, (\ref{ap}) gives the exact amplitude and phase if all the spectral content of $A(t)$ lies lower than $\nu$, i.e.\ $\hat{A}(\xi\geq\nu)=0$. Usually, however, there is a small discrepancy between the true and analytic amplitude/phase (considered in detail in \ref{app:asigerr}), but it is often very small. Thus, there is arguably still nothing better than the analytic signal approach for amplitude and phase estimation in the case of a single AM/FM component (\ref{amfm}) \cite{Vakman:96} (except, maybe, the recently proposed direct quadrature method \cite{Huang:09}).

However, real-life signals rarely consist of only one component, and they usually also contain noise. In this case, the analytic signal will represent a mix of the amplitude and phase dynamics of all components contained in the signal (additionally corrupted by noise), so their individual parameters cannot be recovered from it. One should therefore employ more sophisticated techniques, able to distinguish the different components within a single time-series. This can be done by using the TFR-based approaches that will be described below. The main assumption behind them (and time-frequency analysis more generally) is that the signal is represented by a sum of AM/FM components, each of which satisfies the analytic approximation (\ref{ap}), plus some noise $\eta(t)$:
\begin{equation}\label{sigc}
\begin{aligned}
&s(t)=\sum_i x_i(t)+\eta(t)=\sum_i A_i(t)\cos\phi_i(t)+\eta(t),\\
&\forall t,i:\;A_i(t)>0,\,\phi_i'(t)>0,\,[A_i(t)\cos\phi_i(t)]^+\approx A_i(t)e^{i\phi_i(t)}/2.
\end{aligned}
\end{equation}
Although signal representation (\ref{sigc}) is not unique, in practice one aims at the sparsest among such representations, i.e.\ the one characterized by the smallest number of components $x_i(t)$. It is also important to note, that the most accurate estimates of the components' parameters obtainable using any TFR-based method are the corresponding analytic estimates (\ref{ap}), so the (best achievable) goal of such methods is to extract the separate analytic signals $x^{a}_i(t)$ for each of the chosen $x_i(t)$ in (\ref{sigc}).

\begin{figure*}[t!]
\begin{center}
\includegraphics[width=0.9\linewidth]{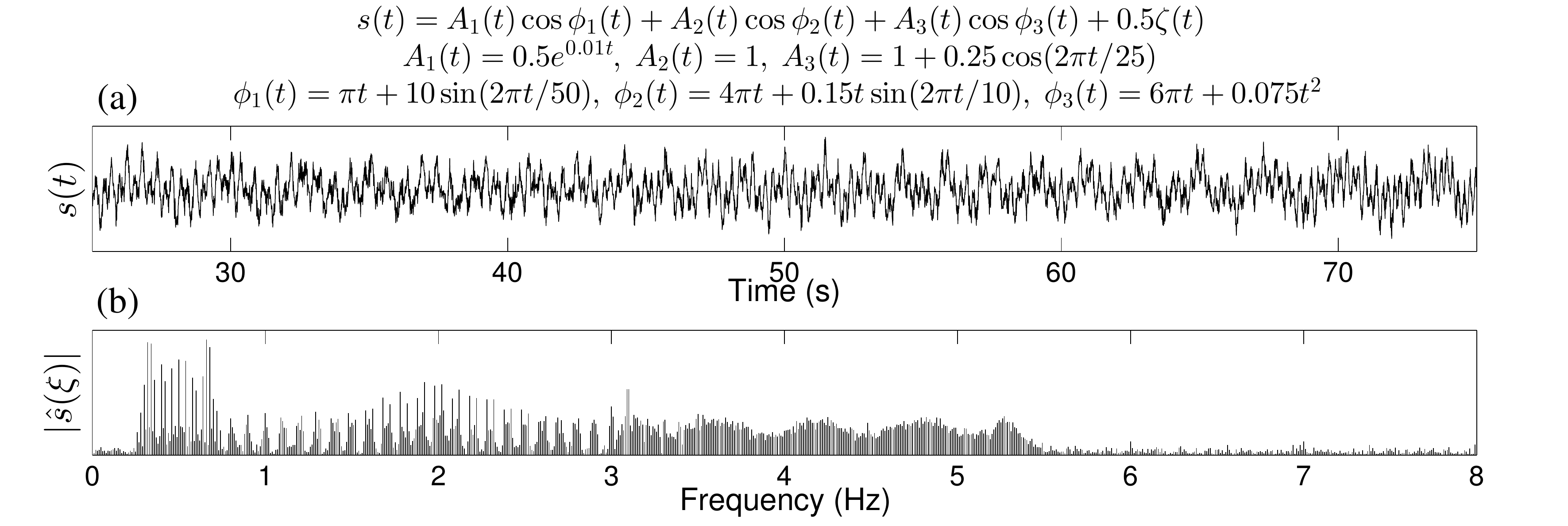}
\includegraphics[width=0.9\linewidth]{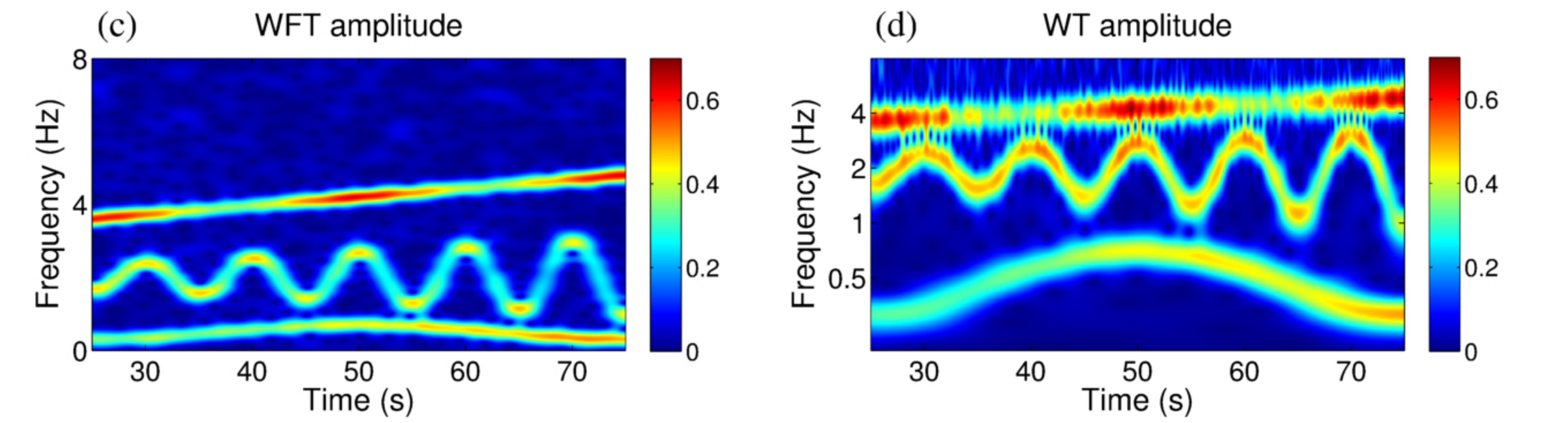}
\end{center}
\caption{Different representations of a signal composed of three AM/FM components and corrupted by noise $0.5\zeta(t)$, where $\zeta(t)$ denotes Gaussian white noise of unit deviation; the parameters of the components are described above the top panel. (a) Signal in the time domain; (b) signal in the frequency domain, given by its Fourier transform; (c,d) signal in the time-frequency domain, given by its WFT and WT (see Secs.\ \ref{sec:tfrWFT} and \ref{sec:tfrWT} below), respectively.}
\label{fig:trfrtfr}
\end{figure*}

The above discussion was related to problems that arise while considering multi-component signals in the time domain. To make the story complete, the problems that arise when trying to treat such signals in the frequency domain should also be addressed. In practice, any signal can be represented by its Fourier expansion (discrete FT), i.e.\ as the sum of tones -- monochromatic signals $A_i\cos(\nu_i t+\varphi_i)$ -- with different constant amplitudes $A_i$, phase shifts $\varphi_i$ and frequencies $\nu_i$. The same is obviously true for a single AM/FM component, and its Fourier expansion is directly related to the spectrum of the corresponding amplitude and frequency modulations:
\begin{equation}\label{amfmft}
\begin{aligned}
x&(t)=A\big(1+\sum_ar_a\cos(\nu_a t+\varphi_a)\big)\cos\big(\nu t+\sum_br_b\sin(\nu_b t+\varphi_b)\big)\\
=&\big/\mbox{expressing }e^{ia\sin \phi}=\sum_{n=-\infty}^\infty J_n(a)e^{in\phi}\big/\\
=&Ae^{i\nu t}\Big[1+\sum_a\frac{r_a\cos(\nu_a t+\varphi_a)}{2}\Big]
\prod_b\sum_{n_b}J_{n_b}(r_b)e^{i(n_b\nu_b t+n_b\varphi_b)}+c.c.\\
=&\sum_{\{n_b\}=\{n_1,n_2,...\}\in\mathds{Z}}
\Big[
\widetilde{A}_{\{n_b\}}e^{i\widetilde{\varphi}_{\{n_b\}}}e^{i\widetilde{\nu}_{\{n_b\}}t}+\\
&+\sum_a \frac{r_a}{2}\widetilde{A}_{\{n_b\}}
\big(e^{i(\widetilde{\varphi}_{\{n_b\}}+\varphi_a)}e^{i(\widetilde{\nu}_{\{n_b\}}+\nu_a)t}
+e^{i(\widetilde{\varphi}_{\{n_b\}}-\varphi_a)}e^{i(\widetilde{\nu}_{\{n_b\}}-\nu_a)t}\big)
\Big]+c.c.,
\end{aligned}
\end{equation}
where
\begin{equation}\label{amfmftapf}
\begin{gathered}
\widetilde{A}_{\{n_b\}}\equiv A\prod_b J_{n_b}(r_b)=AJ_{n_1}(r_1)J_{n_2}(r_2)...,\\
\widetilde{\nu}_{\{n_b\}}\equiv\nu+\sum_b n_b\nu_b,\quad \widetilde{\varphi}_{\{n_b\}}\equiv\sum_b n_b\varphi_b,
\end{gathered}
\end{equation}
and $J_n(r_b)=(-1)^nJ_{-n}(r_b)$ denote Bessel functions of the first kind, while $c.c.$ stands for the complex conjugate of the preceding expression. Thus, any component (\ref{amfmft}) can be represented as a sum of tones with frequencies $|\nu+\sum_b n_b\nu_b|$ and $|\nu\pm\nu_a+\sum_b n_b\nu_b|$ (for all possible combinations of $a\in\mathds{N}$ and $\{n_b\}=\{n_1,n_2,...\}\in\mathds{Z}$). But, obviously, not any sum of tones can be represented as the AM/FM component. What make this possible in (\ref{amfmft}) are the specific relationships between the tones' amplitudes, phases and frequencies, which encode the corresponding amplitude and frequency modulations of the component.

Clearly, the representation of $x(t)$ (\ref{amfmft}) as a single entity is more compact and meaningful than its representation as the sum of tones, which corresponds to its FT. Thus, for a multi-component signal $s(t)$ (\ref{sigc}), each component is encoded as a set of entries in the associated frequency representation $\hat{s}(\xi)$, but it is usually unclear which entries correspond to what component, and which ones are attributable to noise. This is illustrated in Fig.\ \ref{fig:trfrtfr}, which shows an example of signal representation in the time, frequency, and time-frequency domains, with the latter being given by its WFT and WT, to be discussed below. As can be seen, although all representations by definition contain the same amount of information about the signal, in the case of Fig.\ \ref{fig:trfrtfr} the most readily interpretable view of this information is provided in the time-frequency domain.

\section{Time-frequency representations (TFRs)}\label{sec:tfr}

Instead of studying a signal in either of the (one-dimensional) time ($s(t)$) or frequency ($\hat{s}(\xi)$) domains, it is often more useful to study its properties in time and frequency simultaneously. This can be done by considering specifically constructed projections of the signal onto a (two-dimensional) time-frequency plane, and these are called time-frequency representations (TFRs). Such an approach gives the possibility of tracking the evolution of the frequency content in time, which is especially suitable for the analysis of multicomponent and nonstationary time-series, characterized by time-varying spectral properties. An essential feature of a TFR is the possibility of extracting the oscillatory components (possibly with time-varying properties) present in the signal. In what follows, we will concentrate exclusively on those TFRs that are suitable for such a procedure.

In what follows, the notions of time, frequency and time-frequency resolution of the TFR will sometimes be used. The time (frequency) resolution can be defined quite generally as the reciprocal of the minimum time (frequency) difference between two time events -- delta-peaks (two frequency events -- tones), for which they still can be resolved reliably in the TFR. Thus, the higher the time resolution of the TFR, the faster the changes in time (including amplitude/frequency modulation) that can accurately be reflected in it; in the same manner, the higher the frequency resolution, the closer in frequency are the components that can be resolved. The time and frequency resolutions are not independent, being inversely proportional to each other: increasing time resolution of the TFR while preserving its general form is accompanied by simultaneous decrease of a frequency resolution, and {\it vice versa}. The (joint) time-frequency resolution, which determines the quality of this trade-off, is (for now) intuitively defined as the reciprocal of the minimum time-frequency area where two time and two frequency events can simultaneously be resolved in the TFR. A more detailed consideration of these issues is postponed to Sec.\ \ref{sec:tfres}.

\subsection{Windowed Fourier Transform (WFT)}\label{sec:tfrWFT}

The windowed Fourier transform (WFT), also called the short-time Fourier transform or (in a particular form) the Gabor transform, is one of the oldest and thus most-investigated linear TFRs. The WFT $G_s(\omega,t)$ of the signal $s(t)$ is usually constructed as:
\begin{equation}\label{wft0}
G_{s}(\omega,t)=\int s(u)g(u-t)e^{-i\omega u}du=e^{-i\omega t}\frac{1}{2\pi}\int e^{i\xi t}\hat{s}(\xi)\hat{g}(\omega-\xi)d\xi,
\end{equation}
where $g(u)$ is the specified \emph{window function} and $\hat{g}(\xi)$ is its FT. Without loss of generality, we consider $|\hat{g}(\xi)|$ to have a maximum at $\xi=0$:
\begin{equation}\label{omg}
\omega_g\equiv\operatorname{argmax}|\hat{g}(\xi)|=0,
\end{equation}
which ensures that, for a single frequency component $s(t)=\cos\nu t$, the amplitude of the WFT has a maximum at the frequency $\omega=\nu$ (see (\ref{Nwft}) below). If this is not the case, one should always set $\omega_g$ to zero manually by considering the demodulated window function $\{g(t),\hat{g}(\xi)\}\rightarrow\{g(t)e^{-\omega_g t},\hat{g}(\xi+\omega_g)\}$. Note that one is usually much more interested in the form of $\hat{g}(\xi)$ than in $g(t)$: the former is used in many formulas and, furthermore, the numerical computation of a WFT is performed in the frequency domain (see below), whereas the time-domain form $g(t)$ is used mainly for estimation of the time resolution (see Sec.\ \ref{sec:tfres}) and boundary effects (see Sec.\ \ref{sec:pract3}).

The negative frequency part of the $\hat{s}(\xi)$ can influence the quality of the resultant WFT and can impose restrictions on the lowest frequency and on the parameters of the window function (see \ref{app:posfreq}). Additionally, to make the WFT phase equal to the phase of the signal, even in the simplest case of a pure tone $s(t)=A\cos(\nu t+\varphi)$, one should multiply (\ref{wft0}) by $e^{i\omega t}$. The WFT should therefore be considered in the slightly modified form
\begin{equation}\label{wft}
G_{s}(\omega,t)=\int s^{+}(u)g(u-t)e^{-i\omega (u-t)}du=\frac{1}{2\pi}\int_0^\infty e^{i\xi t}\hat{s}(\xi)\hat{g}(\omega-\xi)d\xi.
\end{equation}
The difference between the (\ref{wft0}) and (\ref{wft}) is usually negligible, apart from a phase shift, but it becomes important for large spreads of $\hat{g}(\xi)$ (see \ref{app:posfreq}).

It is easy to see that, for a signal represented as a sum of pure frequencies (which always can be effected through the FT), the corresponding WFT (\ref{wft}) will take the form
\begin{equation}\label{Nwft}
s(t)=\sum_nA_n\cos(\nu_n t+\varphi_n)\Rightarrow G_s(\omega,t)=\frac{1}{2}\sum_n A_n\hat{g}(\omega-\nu_n)e^{i(\nu_n t+\varphi_n)},
\end{equation}
from which one can gain a first impression of how it works. It should be noted, that although peaks in $|W_s(\omega,t)|$ are all at positive frequencies, as follows from (\ref{Nwft}) and the assumptions made, the WFT (\ref{wft}) is defined on the full interval $\omega\in(-\infty,\infty)$ and one can calculate it for $\omega<0$ as well tracing the peak's tails. Nevertheless, it is often enough to calculate the WFT at positive frequencies only.

The WFT is an invertible transform, so that the original signal in both time and frequency domains can be recovered from it as (see \ref{app:derdirrec}):
\begin{equation}\label{iwft}
\begin{gathered}
s^{a}(t)=C_g^{-1}\int G_s(\omega,t)d\omega,\;s(t)=\langle s(t)\rangle+\re[s^{a}(t)],\\
\hat{s}(\omega>0)=\widetilde{C}_g^{-1}\int G_s(\omega,t)e^{-i\omega t}dt,\;\hat{s}(-\omega)=\hat{s}^*(\omega),\\
C_g\equiv\frac{1}{2}\int \hat{g}(\xi)d\xi=\pi g(0),\;\widetilde{C}_g\equiv\int g(t)dt=\hat{g}(0).
\end{gathered}
\end{equation}
Furthermore, $\int\int |G_s(\omega,t)|^2d\omega dt=\frac{1}{2}\int |s(t)|^2 dt\int|\hat{g}(\xi)|^2d\xi$ \cite{Mallat:08,Cohen:95,Grochenig:01}, so that the squared WFT modulus $|G_s(\omega,t)|^2$ -- the spectrogram -- has the energy conservation property (up to normalization). It does not, however, satisfy the marginal conditions, i.e.\ the time (frequency) average of the spectrogram is not proportional to $|\hat{s}(\xi)|^2$ ($|s(t)|^2$), being the smoothed version of the latter \cite{Cohen:95,Hlawatsch:10}. Nonetheless, both $|s(t)|^2$ and $|\hat{s}(\xi)|^2$ can easily be recovered from the signal's WFT using (\ref{iwft}).

In numerical applications, given the signal $s(t_n),n=1,..,N$ sampled at frequency $f_s$, its WFT is calculated using the FFT algorithm, utilizing the frequency domain form of (\ref{wft}). One first divides the frequency axis into equispaced bins $\omega_k=(k-k_0)\Delta\omega$, $k\in \mathds{Z}$, where $\Delta\omega$ is the frequency bin width (its optimal choice and related issues are discussed in Sec.\ \ref{sec:pract2}). One then computes the signal's discrete FT $\hat{s}(\xi)=\{\hat{s}(\xi_n)\}$, where $\xi_n=n f_s/N,n=0,\pm1,..,\pm N/2$ are the frequencies of the discrete Fourier transform, and sets $\hat{s}(\xi_i\leq0)=0$. Finally, taking the inverse FT of $\{s(\xi_n)\hat{g}(\omega_k-\xi_n)\}$ at each frequency $\omega_k$ gives the full WFT time evolution $G_{s}(\omega_k,t_{n=1,..,N})$ for that frequency bin. The complete numerical implementation, including all issues, is given in Sec.\ \ref{sec:stepalg}. As to the units of a WFT, the window function $g(u)$ is dimensionless (but $\hat{g}(\xi)$ is not, being in ${\rm Hz}^{-1}$), so the units of $G_s(\omega,t)$ and $\langle G_s(\omega,t)\rangle$ are [signal units]$\times {\rm Hz}^{-1}$, as can be seen from (\ref{wft}); the $C_g$ (\ref{iwft}) used in reconstruction is dimensionless.

\subsubsection*{Gaussian window}

Unless otherwise specified, all considerations and formulas in this work apply for an arbitrary window $g(t),\hat{g}(\xi)$, while some particular window forms are given in \ref{app:winwav} together with the expressions for the related quantities and their resolution characteristics. However, in what follows we pay particular attention to the Gaussian window function
\begin{equation}\label{gw}
\begin{aligned}
&g(u)=\frac{B_{f_0}}{\sqrt{2\pi}f_0}e^{-u^2/2f_0^2}\Leftrightarrow \hat{g}(\xi)=B_{f_0}e^{-f_0^2\xi^2/2},\\
\end{aligned}
\end{equation}
which we use for simulations. It is commonly used on account of its unique property of maximizing the ``classic'' joint time-frequency resolution of the transform, see e.g. pp.\ 43-45 in \cite{Mallat:08} (different resolution properties of the TFR will be discussed in Sec.\ \ref{sec:tfres} below). The \emph{resolution parameter} $f_0$ controls the spread of a Gaussian window in time and frequency, thus determining the trade-off between the time and frequency resolution of the WFT: the smaller it is, the faster changes in time are reflected in the transform, but the harder it becomes to distinguish between two components that are close in frequency. In the literature, the choice of $f_0$ is usually based on visual inspection of the resultant WFT. The choice of window parameters for a general $g(t),\hat{g}(\xi)$ will be discussed in detail in Part II.

$B_{f_0}$ is a normalization factor that generally depends on $f_0$ and can be chosen in a number of ways based on convenience. For example, given a monochromatic signal $s(t)=A\cos\nu t$, the normalization $B_{f_0}=2$ makes the peak in the WFT amplitude equivalent to the amplitude at $\nu$ ($\max_\omega |G_s(\omega,t)|=A,\forall t$), whereas the use of $B_{f_0}=2f_0/\sqrt{2\pi}\Rightarrow C_g=1$ makes total area under the peak equal to the amplitude ($|\int W_s(\omega,t)d\omega|=A,\forall t$); on the other hand, the value $B_{f_0}=\pi^{1/4}\sqrt{2f_0}$ normalizes the window function energy $\int|g(u)|^2du=1$ ($=\frac{1}{2\pi}\int |\hat{g}(\xi)|^2d\xi$ due to Parseval's identity). For simplicity, we use $B_{f_0}=1$.

\subsection{Wavelet transform (WT)}\label{sec:tfrWT}

The (continuous) wavelet transform (WT) is the other well-known linear TFR which, in contrast to the WFT, has logarithmic frequency resolution; in other respects the two TFRs are quite similar. The WT $W_s(\omega,t)$ of a signal $s(t)$ for a chosen wavelet function $\psi(u)$ is constructed as
\begin{eqnarray}\label{wt0}
W_s(\omega(a),t)&=&\int s(u)\psi^*{\Big(}\frac{u-t}{a}{\Big)}\frac{du}{a} \nonumber \\
&=&\int s(au+t)\psi^*(u)du \\
&=&\frac{1}{2\pi}\int e^{i\xi t}\hat{s}(\xi)\hat{\psi}^*(a\xi)d\xi, \nonumber
\end{eqnarray}
where $a=a(\omega)$ denote scales, inversely proportional to the frequencies $\omega$ with a law defined by the chosen wavelet. They can be expressed through the position of the peak $\omega_{\psi}$ in the wavelet's FT as
\begin{equation}\label{wtsc}
a(\omega)=\frac{\omega_{\psi}}{\omega}\mbox{,   }\omega_\psi\equiv\operatorname{argmax}|\hat{\psi}(\xi)|,
\end{equation}
which establishes that, in the simplest case $s(t)=A\cos(\nu t+\varphi)$, the modulus of the WT $|W_s(\omega,t)|$ will have the peak exactly at the tone frequency $\omega=\nu$. Like the WFT, the wavelet is mainly characterized by its form in the frequency domain $\hat{\psi}(\xi)$, with the time-domain form $\psi(t)$ being less important.

\begin{remark}\label{rem:wtrem}
There are some important conventions for normalizing the WT. Generally one has $du/a^{p}$ in the first integral in (\ref{wt0}): we choose $p=1$, but quite often in the literature $p=1/2$ is used, corresponding to $W_s'=\sqrt{a}W_s$. The latter choice is motivated by the invariance of the \emph{integrated} wavelet power $\sim\int |W_s'(\omega,t)|^2d\omega$, which will be the same e.g.\ for $\sin 2\pi t$ and $\sin4\pi t$, while $\int |W_s(\omega,t)|^2d\omega$ will differ. However, this comes at the expense of bias in the wavelet amplitude and power at each frequency towards the lower frequencies (higher scales). For example, if one has a two-tone signal $s(t)=\cos\nu_1 t+\cos \nu_2 t$ with $\nu_1\ll\nu_2$, then there will be two equal peaks in $|W_s(\omega,t)|$ at $\omega=\nu_{1,2}$, but peaks in $|W_s'(\omega,t)|$ will have different heights: that corresponding to the lower frequency $\nu_1$ will be higher, even though the amplitudes of the tones are the same. Moreover, the peaks in the former will be located exactly at $\omega=\nu_{1,2}$, while those in latter will be shifted. All this makes the interpretation of the resultant WT for $p=1/2$ more complicated and less straightforward than for $p=1$, and it can be quite misleading for inexperienced users. Additionally, the WT formula (\ref{wt0}), corresponding to $p=1$, in the frequency domain is completely analogous to that for the WFT (\ref{wft0}), while for $p=1/2$ there is an additional multiplier $\sqrt{a}$; the same considerations apply to reconstruction of the components parameters from WFT/WT peaks (ridge reconstruction, see Sec.\ \ref{sec:tfsrecB}). The bias in $\int |W_s(\omega,t)|^2d\omega$ for $p=1$ occurs because the linear frequency scale is not natural for the WT and, as also with other measures, the power should be integrated on a logarithmic scale, i.e.\ redefined as $\int |W_s(\omega,t)|^2\frac{d\omega}{\omega}$ (with a constant normalization multiplier), which is unbiased. We suggest, therefore, that the choice $p=1$, corresponding to (\ref{wt0}), is much more natural and convenient than $p=1/2$, as well as being simpler in terms of the analytical derivation. Nevertheless, the two normalizations differ only in terms of convenience and ease of understanding: the correct usage of both will by definition give the same results.
\end{remark}

To remove the contributions of the negative frequency part, and the consequent restrictions on wavelet parameters (see \ref{app:posfreq}), instead of full signal $s(t)$ one should use only its positive frequency part $s^{+}(t)$ in (\ref{wt0}). With this in mind, and using (\ref{wtsc}), the ``final'' WT formula becomes
\begin{eqnarray}\label{wt}
W_s(\omega,t)&=&\int s^{+}(u)\psi^*{\Big(}\frac{\omega(u-t)}{\omega_\psi}{\Big)}\frac{\omega du}{\omega_\psi} \nonumber \\
&=&\int s(u){\Big[}\psi^{+}{\Big(}\frac{\omega(u-t)}{\omega_\psi}{\Big)}{\Big]}^*\frac{\omega du}{\omega_\psi} \nonumber\\
&=&\frac{1}{2\pi}\int_0^\infty e^{i\xi t}\hat{s}(\xi)\hat{\psi}^*(\omega_{\psi}\xi/\omega)d\xi\\
&&\big(\omega_{\psi}\equiv \operatorname{argmax}|\hat{\psi}(\xi)|\big),\nonumber
\end{eqnarray}
where $\psi^{+}(t)$ is the positive frequency part of the wavelet. Thus, (\ref{wt}) can be viewed as the original, unmodified WT (\ref{wt0}), but with a slightly different wavelet $\psi^{+}(t)$ in place of $\psi(t)$. The difference between (\ref{wt}) and (\ref{wt0}) becomes significant only for wavelets with considerable $\hat{\psi}(\xi<0)$, while many wavelets (called analytic) {\it a priori} have $\hat{\psi}(\xi<0)=0$, making (\ref{wt}) and (\ref{wt0}) equal.

The analogue of (\ref{Nwft}) for WT (\ref{wt}) is
\begin{equation}\label{Nwt}
s(t)=\sum_nA_n\cos(\nu_n t+\varphi_n)\Rightarrow W_s(\omega,t)=\frac{1}{2}\sum_nA_n \hat{\psi}^*{\big(}\omega_\psi\frac{\nu_n}{\omega}{\big)}e^{i(\nu_n t+\varphi_n)}.
\end{equation}
As can be seen, in contrast to the WFT, which considers the frequency difference $\nu-\omega$ (\ref{Nwft}), the WT considers the ratio $\nu/\omega$ (or the difference between logarithms $\log \nu/\omega=\log\nu-\log\omega$), which amounts to the definition of the logarithmic frequency resolution. The natural frequency scale is therefore also logarithmic for the WT, in contrast to the linear one for the WFT.

\begin{remark}\label{rem:reschwt}
It is often not recognized that, in contrast to the case of the WFT, rescaling the wavelet function $\{\psi(t),\hat{\psi}(\xi)\}\rightarrow\{\psi(rt),\hat{\psi}(\xi/r)\}$ has no effect on the resultant WT. What matters is $\hat{\psi}(\omega_\psi\nu/\omega)$, while direct rescaling changes both $\hat{\psi}(\xi)$ and the peak frequency $\omega_\psi\rightarrow r\omega_\psi$, thus having no overall effect when substituted in (\ref{wt}) and (\ref{Nwt}). As an illustration, the parameter $k$ in the wavelet $\hat{\psi}(\xi)=\xi^ae^{-k\xi^b}$ is completely redundant, as this wavelet can be presented in the form $\hat{\psi}(\xi)=k^{-a/b}(k^{1/b}\xi)^ae^{-(k^{1/b}\xi)^b}$, which up to the constant multiplier is equivalent to $\hat{\psi}(\xi)=\xi^ae^{-\xi^b}$. To really change the tradeoff between time and frequency resolution of the wavelet, one needs to change either its $\omega_\psi$ while preserving the spread of $|\hat{\psi}(\xi)|$, or {\it vice versa}. Furthermore, the two approaches are completely equivalent, as can be seen from (\ref{wt}), so it is redundant to define parameters controlling both the peak frequency of a wavelet and its bandwidth.
\end{remark}

The reconstruction formulas in the case of the WT become (see \ref{app:derdirrec})
\begin{equation}\label{iwt}
\begin{gathered}
s^{a}(t)=C_\psi^{-1}\int_0^\infty W_s(\omega,t)\frac{d\omega}{\omega},\;s(t)=\langle s(t)\rangle+\re[s^{a}(t)],\\
\hat{s}(\omega>0)=\widetilde{C}_\psi^{-1}\int W_s(\omega,t)e^{-i\omega t}dt,\;\hat{s}(-\omega)=\hat{s}^*(\omega),\\
C_\psi\equiv\frac{1}{2}\int_0^\infty \hat{\psi}^*(\xi)\frac{d\xi}{\xi},\;\widetilde{C}_\psi\equiv\Big[\int \psi(t)e^{-i\omega_\psi t}dt\Big]^*=\hat{\psi}^*(\omega_\psi),
\end{gathered}
\end{equation}
where, in contrast to the WFT (\ref{iwft}), the signal is reconstructed from its WT by integrating $W_s(\omega,t)$ over the frequency logarithm $d\omega/\omega=d\log\omega$, which is standard for WT-based measures. It can also be shown that $\int\int |W_s(\omega,t)|^2\frac{d\omega}{\omega} dt=\frac{1}{2} \int |s(t)|^2dt$ $\times\int|\hat{\psi}(\xi)|^2\frac{d\xi}{\xi}$ \cite{Mallat:08,Hlawatsch:10,Grochenig:01}, i.e.\ the squared WT modulus $|W_s(\omega,t)|^2$ -- the normalized scalogram -- has the energy conservation property on a logarithmic frequency scale (up to normalization). The normalized scalogram does not, however, satisfy the marginal conditions, as its time and frequency averages represent smoothed versions of $|\hat{s}(\xi)|^2$ and $|s(t)|^2$, respectively \cite{Hlawatsch:10}. Nonetheless, the latter can easily be recovered using (\ref{iwt}).

From (\ref{Nwt}) and (\ref{iwt}) it is clear, that for a meaningful WT the wavelet FT $\hat{\psi}(\xi)$ should vanish at zero frequency (or, which is the same, integral of the wavelet $\psi(t)$ over time should be zero):
\begin{equation}\label{adc}
\hat{\psi}(0)=\int \psi(t)dt=0,
\end{equation}
which is called the \emph{admissibility condition}. Evidently, for wavelets not satisfying (\ref{iwt}) the value of $C_\psi$ in (\ref{iwt}) is infinite, and signal reconstruction from the WT becomes impossible. Even more importantly, from (\ref{Nwt}) it follows that, if $\hat{\psi}(0)\neq0$, then each component will be spread over the whole frequency range of the WT, leading to a highly corrupted representation. Thus, for a single tone signal $s(t)=A\cos\nu t$, the WT amplitude approaches $|W_s(\omega,t)|=A|\hat{\psi}(\omega_\psi\nu/\omega)|\rightarrow A|\hat{\psi}(0)|$ as $\omega\rightarrow\infty$. Hence $|\hat{\psi}(0)|$ determines the minimum level to which the WT amplitude for each component decays: if it is nonzero then even components with infinitely distant frequencies will interfere with each other. Drawing an analogy with the WFT, use of the wavelet with $\hat{\psi}(0)\neq 0$ corresponds to use of a window with $\hat{g}(\xi\rightarrow\infty)\neq 0$ in (\ref{wft}), which is obviously inappropriate. Each wavelet should therefore be admissible, i.e.\ should satisfy (\ref{adc}).

In numerical applications, given the signal $s(t_n),n=1,..,N$ sampled at frequency $f_s$, its WT is calculated using the frequency domain form of (\ref{wt}) and taking advantage of the FFT. One first divides the frequency axis into equilogspaced bins, which can be taken in the form $\omega_k/2\pi=2^{(k-k_0)/n_v}$. The numerical parameter $n_v$, called the \emph{number-of-voices}, represents the number of frequency bins in each diadic interval, thus determining the fineness of frequency binning (its selection is discussed in Sec.\ \ref{sec:pract2}). One then computes the signal's discrete FT $\hat{s}(\xi)=\{\hat{s}(\xi_n)\}$, where $\xi_n=n f_s/N,n=0,\pm1,..,\pm N/2$ are the frequencies of the discrete Fourier transform, and sets $\hat{s}(\xi_n\leq0)=0$. Finally, taking the inverse Fourier transform of $\{\hat{s}(\xi_n)\hat{\psi}^*(\omega_\psi\xi_n/\omega_k)\}$ at each frequency $\omega_k$ gives the full time-evolution of the WT $W_s(\omega_k,t_{n=1,..,N})$ for this frequency bin. The complete numerical implementation, including all issues, is given in Sec.\ \ref{sec:stepalg}. As to the units of WT, because the wavelet function $\psi(u)$ is dimensionless (while $\hat{\psi}(\xi)$ is then in ${\rm Hz}^{-1}$), it follows from (\ref{wt}) that $W_s(\omega,t)$ and $\langle W_s(\omega,t)\rangle$ are measured in [signal units]$\times{\rm Hz}^{-1}$, similarly to the WFT; however, $C_\psi$ (\ref{iwt}) is not dimensionless and has units of $\hat{\psi}(\xi)$, i.e.\ ${\rm Hz}^{-1}$.

\subsubsection*{Morlet wavelet}

Except where otherwise specified, all considerations and formulas in this work apply for an arbitrary wavelet $\psi(t),\hat{\psi}(\xi)$. Some particular examples of wavelet forms are listed in \ref{app:winwav} together with the expressions for related quantities and their resolution characteristics. However, one of the most commonly used wavelets (at least in terms of the continuous WT considered here), the Morlet wavelet \cite{Morlet:83}, is worthy of special consideration. It is constructed in analogy with the Gaussian window (\ref{gw}) and takes the form
\begin{eqnarray}\label{mw}
\psi(u)&=&\frac{B_{f_0}}{\sqrt{2\pi}}\left(e^{i2\pi f_0u}-e^{-(2\pi f_0)^2/2}\right)e^{-u^2/2}, \nonumber \\
\hat{\psi}(\xi)&=&B_{f_0}e^{-(\xi-2\pi f_0)^2/2}\left(1-e^{-2\pi f_0\xi}\right),
\end{eqnarray}
where, in analogy with (\ref{gw}), $f_0$ is the \emph{resolution parameter} and $B_{f_0}$ is a normalization constant. Except where specified, we perform all WT simulations using the (\ref{mw}) with $B_{f_0}=1$. In other contexts we leave the normalization unspecified in order to preserve generality. The fact that the Morlet wavelet is used so commonly for the continuous wavelet transform is attributable to the belief that it maximizes the time-frequency resolution. However, this is in fact not exactly correct because the WT has logarithmic frequency resolution, whereas the ``classic'' maximization property of the Gaussian window is valid only for the case of linear resolution (see Sec.\ \ref{sec:tfres} and \ref{app:winwav}). Similarly to the WFT case, the choice of $f_0$ is usually based on visual inspection of the WT, or often set {\it a priori} to $f_0=1$. The effects of different wavelet parameters, and their choice for arbitrary $\psi(t),\hat{\psi}(\xi)$, will be discussed in detail in Part II.

\begin{remark}As already mentioned, there are two ways of parametrizing the wavelet: either through the parameter controlling its peak frequency $\omega_\psi$, or through that responsible for its spread in time and frequency. For the Morlet wavelet, the former approach corresponds to (\ref{mw}) and is utilized here since is so commonly used; in the latter case the Morlet wavelet takes an equivalent form $\psi(u)=\frac{B_{f_0}}{\sqrt{2\pi}f_0}\left(e^{i2\pi u}-e^{-(2\pi f_0)^2/2}\right)e^{-u^2/2f_0^2}\Leftrightarrow \hat{\psi}(\xi)=B_{f_0}e^{-f_0^2(\xi-2\pi)^2/2}\left(1-e^{-2\pi f_0^2\xi}\right)$. Nevertheless, as previously discussed, both definitions are completely equivalent and, as can be easily checked, both parametrizations with the same $f_0$ will give exactly the same WTs. It should be also noted, that the parameter $f_0$ in (\ref{mw}) is often called the central frequency of the Morlet wavelet. However, this appellation, although widespread, is not entirely appropriate. Thus, in (\ref{mw}) $f_0$ indeed controls the wavelet peak frequency $\omega_\psi$ but, due to the admissibility term (see below), it never \emph{exactly} coincides either with $\omega_\psi/2\pi$, or with the wavelet mean frequency (although for $f_0\gtrsim 0.5$ it is very close to both). Taken additionally the duality of possible wavelet parametrization and the similarities with Gaussian window (\ref{gw}), ``resolution parameter'' seems to be more suitable name for $f_0$.\end{remark}

The second term in $\psi(u)$ (\ref{mw}) is called the \emph{admissibility term}, and is needed to establish the wavelet admissibility condition (\ref{adc}). Its relative contribution $e^{-(2\pi f_0)^2/2}$ is usually very small, being $<0.01$ for $f_0>0.5$ (so it is often inappropriately neglected), but for small $f_0$ it becomes considerable. The admissibility term represents the ``unavoidable curse'' of the Morlet wavelet, greatly complicating its form and hence the analysis, leading e.g.\ to nonlinear proportionality between $\omega_\psi$ and $f_0$. The most significant drawback, however, is that the larger $e^{-(2\pi f_0)^2/2}$ is, the lower is the joint time-frequency resolution of the Morlet wavelet (see Sec.\ \ref{sec:tfres} and \ref{app:winwav}). This is in contrast to Gaussian window WFT (\ref{gw}), which preserves its properties however small $f_0$ becomes. The existence of considerable wavelet power at negative frequencies $\hat{\psi}(\xi<0)$ for a non-negligible admissibility term is commonly regarded as an additional serious drawback of the Morlet wavelet, but this problem does not arise when using the ``correct'' WT formula (\ref{wt}), implicitly setting $\hat{\psi}(\xi<0)=0$.

\subsubsection*{Lognormal wavelet}

Most wavelets, such as Morlet, are constructed in a similar way to windows for WFT, i.e.\ considering $\hat{\psi}(\xi)$ on a linear frequency scale. However, as mentioned earlier, the WT has logarithmic frequency resolution, so it seems more appropriate to construct a wavelet using $\log\xi$ as its argument. Therefore, a more correct WT analogy to the Gaussian window (\ref{gw}) would be probably not Morlet, but the lognormal wavelet
\begin{equation}\label{lw}
\hat{\psi}(\xi)\sim e^{-(2\pi f_0\log\xi)^2/2},\;\xi>0
\end{equation}
As investigated in \ref{app:winwav}, the resolution properties of the wavelet (\ref{lw}) are generally better than that of the Morlet, and it has a variety of other useful properties, which make lognormal wavelet a preferable choice. Thus, in contrast to the Morlet wavelet, this wavelet is ``infinitely admissible'', i.e.\ all its moments $\int \xi^{-n}\hat{\psi}(\xi)d\xi/\xi,n\geq0$ are finite, which gives the possibility to reconstruct any order time-derivative of the component's amplitude and phase from its WT (see Sec.\ \ref{sec:tfsrecB} and \ref{app:derdirrec}). Next, being symmetric on a logarithmic scale (i.e.\ $\hat{\psi}(\omega_\psi\omega)=\hat{\psi}(\omega_\psi/\omega)$), it also provides some advantages in terms of components reconstruction, as well as makes behavior of the WT better defined (see Part II). Finally, the lognormal wavelet is more analytically tractable, allowing to obtain $C_\psi$ and many other quantities in the explicit form. Nevertheless, in spite of all these advantages, we will still employ the Morlet wavelet (\ref{mw}) in our simulations just because, apart from being the most widespread choice, it has more in common with the other wavelet forms and thus better demonstrates what one typically gets.

\subsection{Synchrosqueezed WFT (SWFT) and WT (SWT)}\label{sec:tfrSS}

Synchrosqueezing \cite{Daubechies:11,Thakur:11,Thakur:13,Auger:13,Daubechies:96} provides a way to construct more concentrated time-frequency representation from the WFT or WT. The underlying idea is very simple, namely to join all WFT/WT coefficients having same instantaneous frequency (the first derivative of the unwrapped-over-time WFT/WT phase) into one SWFT/SWT coefficient. In mathematical terms, the definition of the WFT/WT instantaneous frequency $\nu_{G,W}(\omega,t)$ is
\begin{equation}\label{iftfr}
\begin{aligned}
\nu_G(\omega,t)=&\frac{\partial}{\partial t}\operatorname{arg}[G_s(\omega,t)]
=\operatorname{Im}\Big[G_s^{-1}(\omega,t)\frac{\partial G_s(\omega,t)}{\partial t}\Big],\\
\nu_W(\omega,t)=&\frac{\partial}{\partial t}\operatorname{arg}[W_s(\omega,t)]
=\operatorname{Im}\Big[W_s^{-1}(\omega,t)\frac{\partial W_s(\omega,t)}{\partial t}\Big].\\
\end{aligned}
\end{equation}

The SWFT $V_s(\omega,t)$ \cite{Thakur:11} and SWT $T_s(\omega,t)$ \cite{Daubechies:11} are then constructed as
\begin{equation}\label{stfr}
\begin{aligned}
&V_s(\omega,t)=C_g^{-1}\int \delta(\omega-\nu_G(\tilde{\omega},t))G_s(\tilde{\omega},t)d\tilde{\omega},\\
&T_s(\omega,t)=C_\psi^{-1}\int_0^\infty\delta(\omega-\nu_W(\tilde{\omega},t))W_s(\tilde{\omega},t)\frac{d\tilde{\omega}}{\tilde{\omega}},\\
\end{aligned}
\end{equation}
where $C_g,C_\psi$ are defined in (\ref{iwft}), (\ref{iwt}). Similarly to the underlying WFT/WT themselves, their synchrosqueezed versions also represent invertible transforms: integrating (\ref{stfr}) over $d\omega$ and using (\ref{iwft}), (\ref{iwt}), one can show that signal can be reconstructed from its SWFT/SWT as
\begin{equation}\label{istfr}
s^{a}(t)=\int_0^\infty V_s(\omega,t)d\omega=\int_0^\infty T_s(\omega,t)d\omega,\;(s(t)=\operatorname{Re}[s^{a}(t)]).
\end{equation}
However, in contrast to the WFT/WT, there is no possibility of reconstructing directly the signal's FT from its SWFT/SWT. Furthermore, the squared magnitudes of the synchrosqueezed transforms lack the energy conservation property (and therefore do not satisfy the marginal conditions either). This can clearly be seen by considering the SWFT/SWT of a finite-duration tone, for which $\int\int|V_s(\omega,t)|^2d\omega dt$ and $\int\int|T_s(\omega,t)|^2d\omega dt$ will be infinite (as $\int|\delta(x)|^2dx=\infty$) even though the signal's energy $\int |s(t)|^2dt$ is finite.

Taking account of (\ref{Nwft}) and (\ref{Nwt}), for a single tone signal $s(t)=A\cos\nu t$ it follows from (\ref{iftfr}) that $\nu_{G,W}(\omega,t)=\nu,\forall t,\omega$, so in this case (\ref{stfr}) gives $V_s(\omega,t)=T_s(\omega,t)=A\delta(\omega-\nu)$, i.e.\ the SWFT/SWT is the perfect representation. For more complicated signals, however, it does not take the simple form like (\ref{Nwft}), (\ref{Nwt}) for WFT/WT. This is mainly because, although constructed from linear additive WFT/WT, the synchrosqueezed TFRs are not linear, i.e.\ the SWFT/SWT of a sum of signals does not equal the sum of the SWFTs/SWTs for each signal separately. This greatly complicates analytical treatment of the synchrosqueezed representations.

\begin{remark}Note, that we introduced $C_g,C_\psi$ directly into the synchrosqueezing formulas (\ref{stfr}) instead of using them in reconstruction (as is done in \cite{Thakur:11,Daubechies:11}). This seems slightly more convenient (e.g.\ in the case of a single tone signal the \emph{numerical} SWFT/SWT amplitude will have a peak equal to the tone amplitude, see below), as well as makes the units of $V_s$ and $T_s$ the same.\end{remark}

Numerically, one computes not the true SWFT/SWT as defined in (\ref{stfr}), but the integrals of $V_s(\omega,t)$ and $T_s(\omega,t)$ over each frequency bin (denote as $\widetilde{V}_s(\omega_k,t)$ and $\widetilde{T}_s(\omega_k,t)$). This is because the SWFT and SWT are generally not analytic, e.g.\ in theory for a single tone signal they are $\delta$-functions, which, obviously, cannot be reliably represented in practice; but the integral of the $\delta$-function is finite and well-defined. Practically, one first calculates the WFT/WT at frequencies $\tilde{\omega}_k$, as previously discussed, and estimates the WFT/WT frequency $\nu_{G,W}(\tilde{\omega}_k,t)$ (\ref{iftfr}). The SWFT/SWT is calculated for the same frequency bins $\omega_k=\tilde{\omega}_k$ as those used for the underlying WFT/WT. At each time $t$ one finds the $\tilde{\omega}_j$ for which $\nu_{G,W}(\tilde{\omega}_j,t)$ lies in the frequency bin centered at $\omega_k$, and then joins all the corresponding WFT/WT coefficients to a single SWFT/SWT entry (\ref{stfr}):
\begin{equation}\label{numstfr}
\begin{aligned}
C_g\widetilde{V}_s(\omega_k,t)\;\equiv\;&\int_{(\omega_{k-1}+\omega_{k})/2}^{(\omega_{k}+\omega_{k+1})/2}V_s(\omega,t)d\omega
=\int_{\omega_k-\Delta\omega/2}^{\omega_{k}+\Delta\omega/2}V_s(\omega,t)d\omega\\
=\;&\sum_{j:\;\frac{\omega_{k-1}+\omega_k}{2}<\nu_G(\tilde{\omega}_j,t)\leq\frac{\omega_k+\omega_{k+1}}{2}}
G_s(\tilde{\omega}_j,t)\frac{\tilde{\omega}_{j+1}-\tilde{\omega}_{j-1}}{2}\\
=\;&\sum_{j:\;-\Delta\omega/2<\nu_G(\tilde{\omega}_j,t)-\omega_k\leq\Delta\omega/2}G_s(\tilde{\omega}_j,t)\Delta\omega,\\
C_\psi\widetilde{T}_s(\omega_k,t)\;\equiv\;&\int_{\sqrt{\omega_{k-1}\omega_{k}}}^{\sqrt{\omega_{k}\omega_{k+1}}}T_s(\omega,t)d\omega
=\int_{2^{-1/2n_v}\omega_k}^{2^{1/2n_v}\omega_{k}}V_s(\omega,t)d\omega\\
=\;&\sum_{j:\;\sqrt{\omega_{k-1}\omega_k}<\nu_W(\tilde{\omega}_j,t)\leq\sqrt{\omega_k\omega_{k+1}}}
W_s(\tilde{\omega}_j,t)\frac{1}{2}\log\frac{\tilde{\omega}_{j+1}}{\tilde{\omega}_{j-1}}\\
=\;&\sum_{j:\;2^{-1/2n_v}<\nu_W(\tilde{\omega}_j,t)/\omega_k\leq2^{1/2n_v}}W_s(\tilde{\omega}_j,t)\frac{\log 2}{n_v},
\end{aligned}
\end{equation}
where we have used the fact that $\omega_k=(k-k_0)\Delta\omega$ for the WFT and $\omega_k/2\pi=2^{(k-k_0)/n_v}$ for the WT. Note that, for the SWT, the summation is done on a logarithmic scale to establish good accuracy (these and related issues are discussed in Sec.\ \ref{sec:pract2}). The SWFT/SWT is usually very sparse, i.e.\ most coefficients are zero so, in numerical applications, one can make advantage of sparse matrices representations. The fast algorithm for synchrosqueezing, with computational time of the same order as for the usual WFT/WT (at least in MatLab), is given in detail in Sec.\ \ref{sec:stepalg}. From (\ref{stfr}) and (\ref{numstfr}), taking into account the units of the WFT/WT and $C_{g,\psi}$, the units of $V_s(\omega,t),T_s(\omega,t)$ are $\{\mbox{signal units}\}\times{\rm Hz}^{-1}$, being the same as for the $G_s(\omega,t),W_s(\omega,t)$ (note that the arguments of the $\delta$-function in (\ref{stfr}) are in Hz, so it is not dimensionless and has units of ${\rm Hz}^{-1}$); the units of $\widetilde{V}_s(\omega,t),\widetilde{W}_s(\omega,t)$ and their time-averages are then simply in the same units as the signal.

It is clear, that for a single tone signal $s(t)=A\cos \nu t$ numerical SWFT/SWT (\ref{numstfr}) will be $\widetilde{V}_s(\omega_k,t)=A\mbox{ if }\nu-\omega_k\in[-\Delta\omega/2,\Delta\omega/2]$, and $=0$ otherwise (similarly for SWT, except logarithmic scale). However, when amplitude and/or frequency modulation is present, the component will be spread over some frequency range in the $V_s(\omega,t)$ and $T_s(\omega,t)$. The amplitudes of $\widetilde{V}_s(\omega_k,t)$ and $\widetilde{T}_s(\omega_k,t)$ will then depend on the bin widths. Importantly, this is reflected in the reconstruction formula (\ref{istfr}), which in terms of a numerical SWFT/SWT (\ref{numstfr}) takes the form of a simple sum, without multiplying by frequency bin widths:
\begin{equation}\label{numistfr}
s^{a}(t)=\sum_k \widetilde{V}_s(\omega_k,t)=\sum_k \widetilde{T}_s(\omega_k,t).
\end{equation}

\begin{remark}To avoid dependence of the SWFT/SWT on frequency discretization, one can additionally divide (\ref{numstfr}) on the frequency bin widths $\Delta\omega_k$, which will also make the numerical and theoretical SWFT/SWT the same. However, the SWFT/SWT amplitude will still depend on the frequency binning, but now in the case when there are one or few distinct (possibly interfering) tones. Moreover, since bin widths depend on frequency for the SWT, such redefinition will lead to a bias in $|\widetilde{T}_s(\omega,t)|$ towards lower frequencies. Therefore, it seems more convenient to leave (\ref{numstfr}) as it is.\end{remark}

\begin{remark}Before performing synchrosqueezing, the calculated WFT/WT used for it is often noise-filtered by hard-thresholding \cite{Daubechies:11,Thakur:11,Thakur:13}, i.e.\ all WFT/WT coefficients having amplitude lower than some specified threshold are set to zero. However, such denoising is useful only in the idealistic case of white Gaussian noise, having same power at all frequencies, whereas in real life this situation is rather an exception than the rule. Thus, generally noise is colored (e.g.\ Brownian), having different powers at different frequencies, and it is often described by distributions other than normal. Therefore, in practise hard-threshold denoising usually does not provide significant advantages; but neither is it harmful. In addition, the choice of an appropriate threshold is usually not obvious. In \cite{Thakur:13} it was proposed to estimate its value from the median absolute deviation of the finest level wavelet coefficients, but even in the case of pure white noise it can highly overestimate the threshold if the signal, apart from the noise, contains AM/FM components with high enough frequencies. Filtering using an overestimated threshold, on the other hand, might be more deleterious than the noise itself. However, denoising, whether by hard-thresholding or by some other filtration scheme, is a useful addition rather than a necessary part of the synchrosqueezing procedure, and it provides similar benefits for the SWFT/SWT and WFT/WT separately. To avoid complications, therefore, we will not perform any denoising in what follows.\end{remark}

\subsection{Difference between the (S)WFT and the (S)WT}\label{sec:dwftwt}

\begin{figure*}[t!]
\begin{center}
\includegraphics[width=0.95\linewidth]{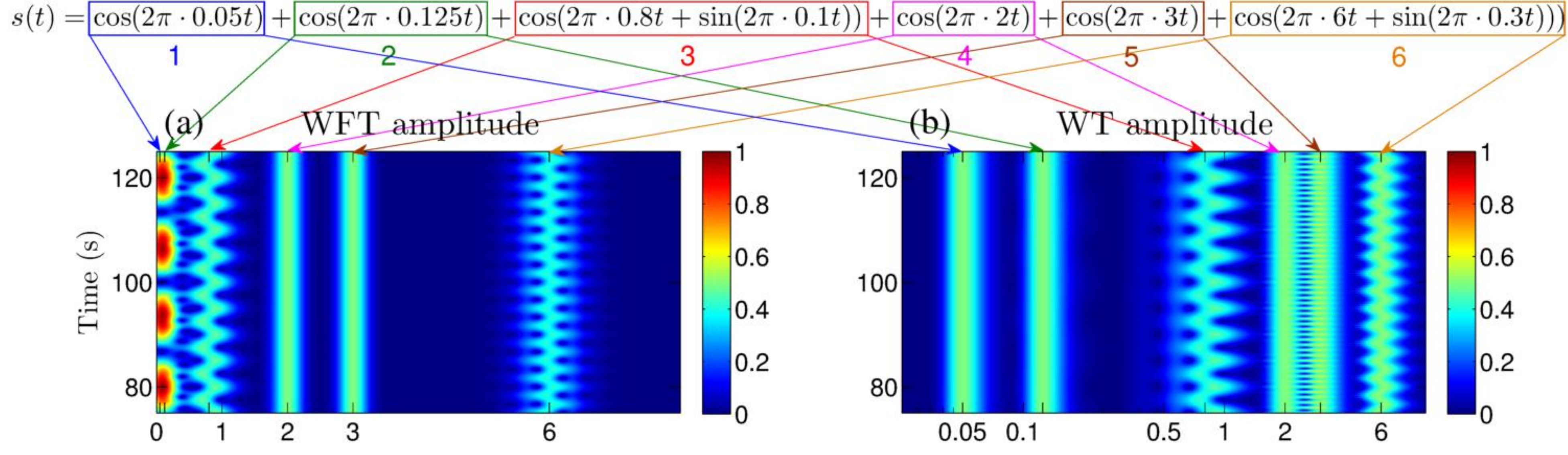}\\
\includegraphics[width=0.95\linewidth]{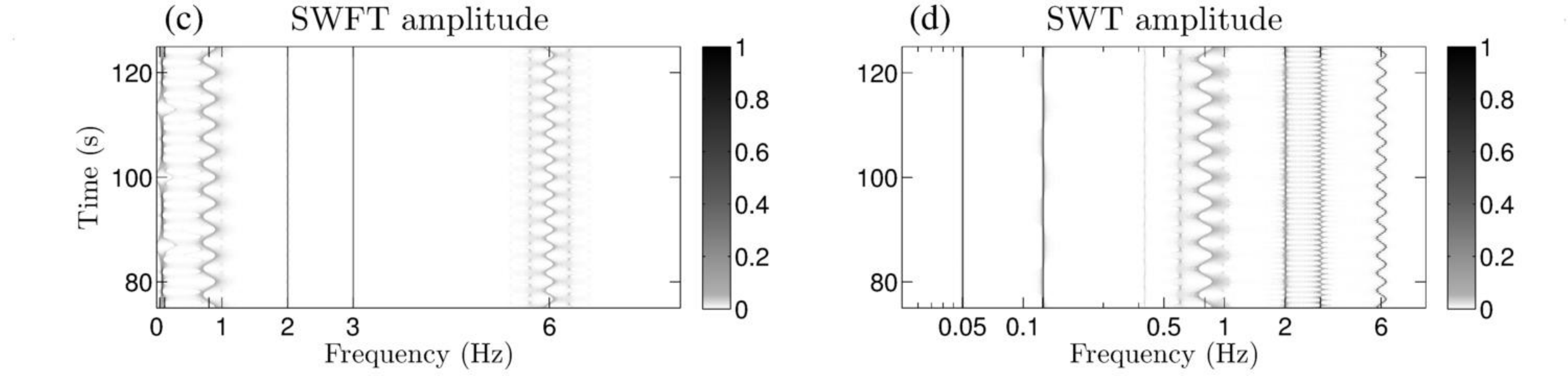}\\
\end{center}
\caption{Examples of different time-frequency representations for the same signal consisting of 6 components, as shown at the figure top: (a) WFT amplitude $|G_s(\omega,t)|$, calculated using Gaussian window (\ref{gw}) with $f_0=1$; (b) WT amplitude $|W_s(\omega,t)|$, calculated using Morlet wavelet (\ref{mw}) with $f_0=1$; (c) SWFT calculated from WFT in (a); (d) SWT calculated from WT in (b). The signal was sampled for 500 s at 100 Hz.}
\label{fig:tfrex}
\end{figure*}

Figure \ref{fig:tfrex} presents examples of the four discussed TFRs calculated for the same signal. Although it is not really the case, we will assume for simplicity of the following illustrative considerations that the WFT and WT in Figure \ref{fig:tfrex} have same time and frequency resolution properties at $\omega=2\pi$ (see Sec.\ \ref{sec:tfres}). First consider the difference between the WFT and WT (compare (a) and (b) in the figure). As mentioned before, the main distinction lies in the different kinds of frequency resolution: linear for WFT, and logarithmic for WT. Thus, if the signal contains two tones $\cos\nu_{1,2}t$, then WFT takes account of their frequency difference $\nu_2-\nu_1$, whereas the WT asks what is the ratio of their frequencies $\nu_2/\nu_1$ (or difference of their logarithms). For example, if there are two oscillations with frequencies $0.05$ and $0.125$ Hz, corresponding to periods $20$ and $8$ s (components 1,2 in Figure \ref{fig:tfrex}), then they will be much better resolved in the WT ($\nu_2/\nu_1-1=1.5$) than in the WFT ($(\nu_2-\nu_1)/2\pi=0.075$). At the same time, if there are two tones with frequencies $2$ and $3$ Hz, corresponding to periods $0.5$ and $\approx0.33$ s respectively (components 4 and 5 in Figure \ref{fig:tfrex}), then they will be much better resolved in the WFT ($(\nu_2-\nu_1)/2\pi=1$ Hz) than in the WT ($\nu_2/\nu_1-1=0.5$).

On the other hand, the WFT has same time resolution at each frequency, while that for WT is $\sim \omega$, i.e.\ increases for increasing frequencies. This means that, for fixed window/wavelet parameters, the higher the frequency, the more time-variability is allowed for components to be reliably represented in a WT, while the WFT does not discriminate in this respect. For example, if one has a frequency-modulated component $\cos(\nu t+\sin\nu_b t)$, the WFT accounts for the value of $\nu_b$ (the smaller the better), while the WT considers $\nu_b/\nu$. Thus, if one has a component with a mean frequency of $0.8$ Hz subject to sinusoidal modulation at $0.3$ Hz (component 3 in Figure \ref{fig:tfrex}), it will be better represented in the WFT ($\nu_b/2\pi=0.2$) than in the WT ($\nu_b/\nu=0.25$). On the other hand, if one has component with mean frequency $6$ Hz and frequency modulation at $0.3$ Hz (component 6 in Figure \ref{fig:tfrex}), it will be represented more reliably in the WT ($\nu_b/\nu=0.05$) than in the WFT ($\nu_b/2\pi=0.3$).

Obviously, one can always adjust $f_0$ (or other window/wavelet parameters, if present) to represent any one chosen component reliably. However, there is often no universal $f_0$ suitable for all components (see Part II), so the choice between the WFT and the WT depends on how many such components can in principle be represented reliably. Therefore, contrary to what is sometimes thought, the WT is not in general superior to the WFT: it just considers all on a logarithmic frequency scale, and whether or not this is more useful than the linear frequency resolution of the WFT depends on the signal structure. Thus, the WT is most suitable when the lower-frequency components are less time-varying and closer in frequency than the components at higher frequencies, i.e.\ what matters is the ratios of periods of the underlying oscillations and their relative modulation; when there are no reasons to think that this might be the case, the WFT seems to be the more appropriate.

\begin{figure*}[t!]
\begin{center}
\includegraphics[width=0.95\linewidth]{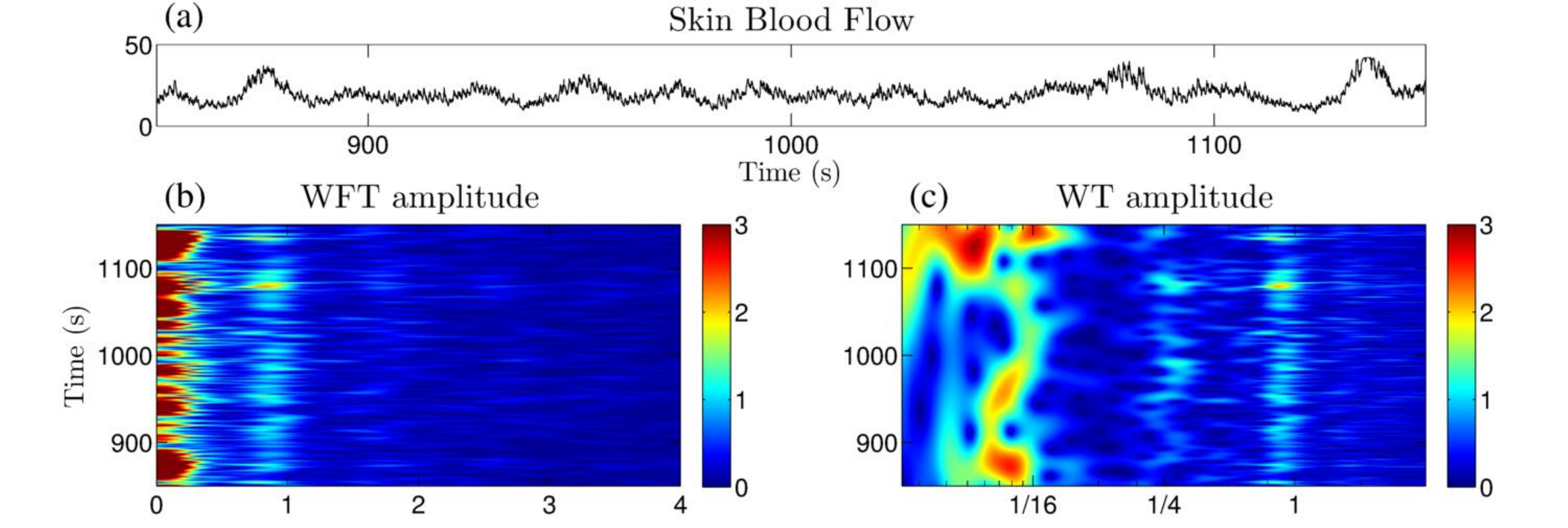}\\
\includegraphics[width=0.95\linewidth]{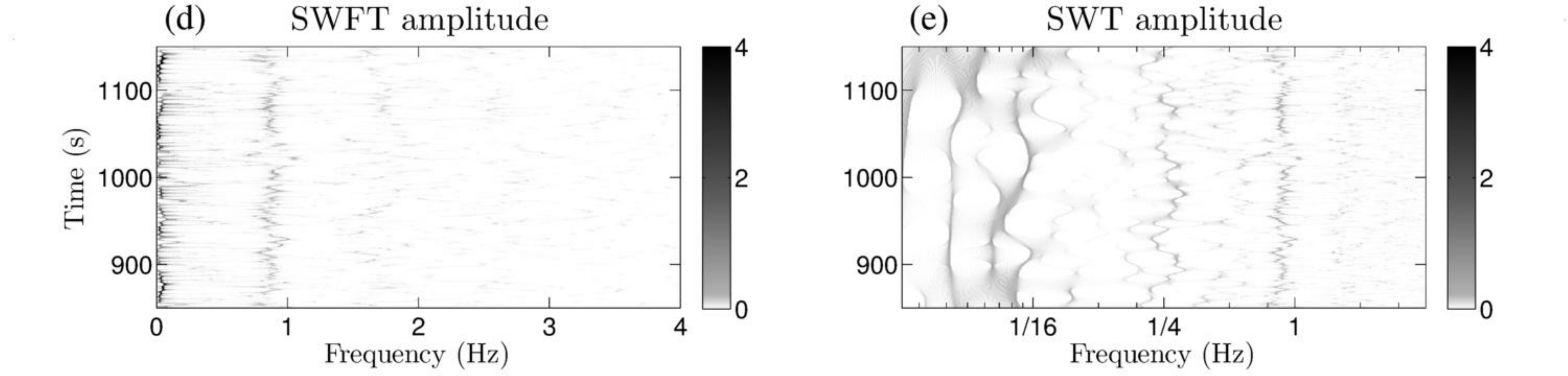}\\
\end{center}
\caption{(a) The human skin blood flow signal (measured by laser-Doppler flowmetry \cite{Nilsson:80} with the probe over the right wrist \emph{caput ulna} \cite{Stefanovska:99b,Shiogai:10}) and (b-e) its different time-frequency representations. The signal was recorded for 30 min. and sampled at 40 Hz.}
\label{fig:tfrSBF}
\end{figure*}

\begin{figure*}[t!]
\begin{center}
\includegraphics[width=0.95\linewidth]{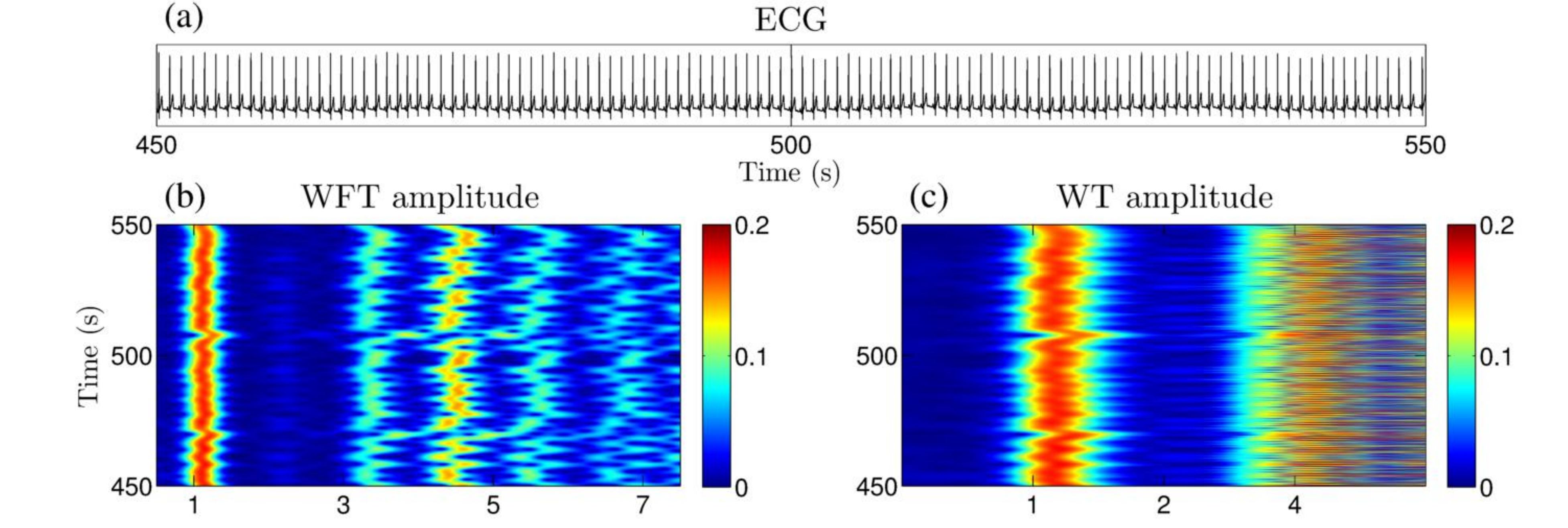}\\
\includegraphics[width=0.95\linewidth]{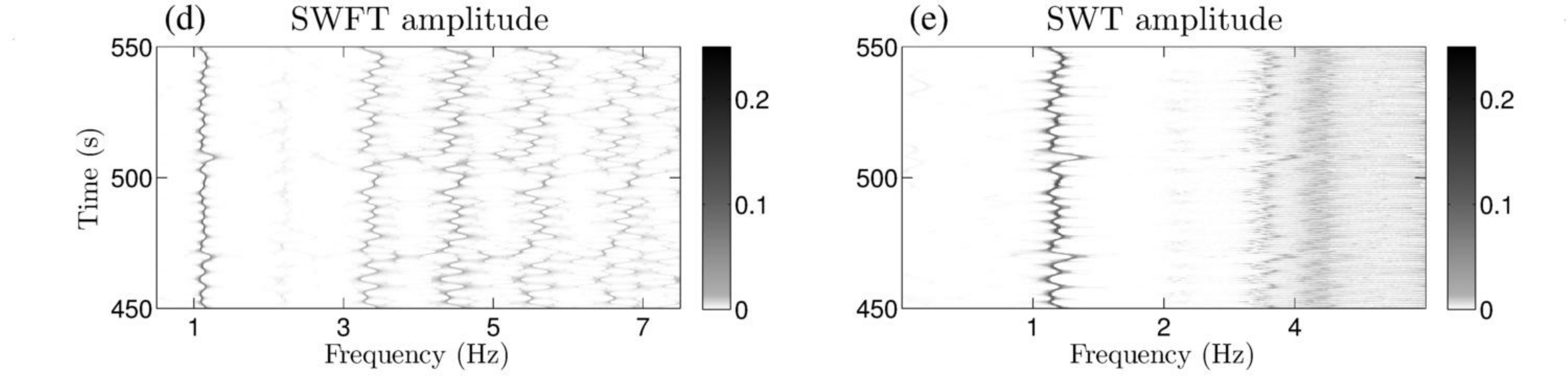}\\
\end{center}
\caption{(a) The human ECG signal (3-lead, with electrodes placed on shoulders and the lowest left rib \cite{Lotric:00}) and (b-e) its different time-frequency representations. The signal was recorded for 30 min. and sampled at 400 Hz.}
\label{fig:tfrECG}
\end{figure*}

This issue is further illustrated in Figs.\ \ref{fig:tfrSBF} and \ref{fig:tfrECG} in relation to two real signals taken as examples -- human microvascular blood flow and an electrocardiogram (ECG), respectively (see also pp.\ 126-134 in \cite{Mallat:08} for an additional discussion and examples). The blood flow signal, shown in Fig.\ \ref{fig:tfrECG}(a), contains a variety of physiologically meaningful oscillations on different scales \cite{Shiogai:10,Stefanovska:99a,Stefanovska:99b}, so that the ratios of their periods are what matters. Hence, the WT is the ``right'' TFR type for such signals and, using it, one can separate the cardiac ($\sim$\,1 Hz), respiratory ($\sim$\,0.25 Hz) and other low-frequency oscillations (Fig.\ \ref{fig:tfrSBF}(c)). In contrast, the WFT of the blood flow (Fig.\ \ref{fig:tfrSBF}(c)) represents well only the cardiac oscillations, while the other ones are merged together and cannot be resolved.

On the other hand, the ECG signal, shown in Fig.\ \ref{fig:tfrECG}(a), has the kind of structure for which a linear frequency scale is more suitable. Thus, the ECG is characterized by a complex waveform, so that its dynamics can be represented as a sum of harmonics -- the AM/FM components with frequencies that are multiples of the fundamental heart frequency \cite{Wu:13}. These harmonics are represented by a number of ``curves'' in the time-frequency plane and, because the frequency difference between them does not scale with frequency, the harmonic structure is better represented in the WFT (Fig.\ \ref{fig:tfrECG}(b)) than in the WT (Fig.\ \ref{fig:tfrECG}(c)), where all the higher harmonics interfere strongly with each other. It should be noted, however, that one is usually interested in the properties of the fundamental component around 1 Hz \cite{Iatsenko:cardio}, which is well represented in both the WFT and WT.

In general, it seems that many real life signals have structures that are more suited to studies based on the WT, which is the reason why this TFR type became so popular. Additionally, one can often analyze a variety of time-series of different kinds using the same wavelet parameters, while for the WFT the window parameters should be adjusted for each particular case; usually, the WT is also computationally cheaper due to its logarithmic frequency scale, requiring fewer bins to cover the same frequency range. On the other hand, the WFT has better resolution properties (see \ref{app:winwav}) and, if $\hat{g}(\xi)$ is symmetric, provides additional advantages for component reconstruction (see Part II).

The difference between the WFT/WT and the SWFT/SWT is nicely illustrated by comparing Fig. \ref{fig:tfrex}(a) with (c), and (b) with (d) (see also Figs.\ \ref{fig:tfrSBF} and \ref{fig:tfrECG}). Thus, synchrosqueezing effectively ``compresses'' the TFR, making the representation denser. This might at first glance seem to imply that the SWFT/SWT is able to resolve closer in frequency components than is possible with the WFT/WT, i.e.\ that synchrosqueezing increases the frequency resolution. This, however, appears to be an illusion: the components which interfere in the WFT/WT (e.g.\ 3,4 in Fig. \ref{fig:tfrex}(a) or 1,2 in (b)) are interfering in the SWFT/SWT as well, so the frequency resolution does not appear to be enhanced. The time resolution also seems not to be improved by synchrosqueezing (compare e.g.\ component 6 in Fig.\ \ref{fig:tfrex}(a) and (c)). To understand why greater concentration does not necessarily imply better resolution properties, consider the WFT/WT with only amplitude peaks left: while obviously not having better time or frequency resolution than the original WFT/WT, such a representation will be extremely dense (and in fact will look very much the same as the SWFT/SWT). Therefore, it is questionable whether synchrosqueezing provides any significant advantages apart from a more visually appealing picture. This question will be investigated thoroughly in Part II.

\subsection{Other time-frequency representations}

The (S)WFT and (S)WT represent only a few out of the many existing types of TFR. These include e.g.\ different quadratic representations, such as the Wigner-Ville \cite{Wigner:32,Ville:48}, Rihaczek \cite{Rihaczek:68} and Choi-Williams \cite{Choi:89} distributions (see \cite{Hlawatsch:10,Boashash:03,Hlawatsch:92} for a comprehensive lists of quadratic TFRs and their properties). The advantage of such representations is that they can provide almost perfect time-frequency localization for particular types of signals, such as chirps. At the same time, quadratic TFRs do not satisfy the additivity property, so that the TFR of a sum of signals does not equal the sum of the TFRs for each signal separately. Consequently, the analysis and interpretation of multicomponent signals are greatly complicated by the appearance of cross-terms bearing no physical meaning. Furthermore, quadratic representations do not allow straightforward signal reconstruction, as opposed to the (S)WFT and (S)WT which allow one to recover the original signal using (\ref{iwft}), (\ref{iwt}) or (\ref{istfr}).

Another TFR family is represented by various WFT/WT modifications, with the most popular among them being the windowed fractional Fourier transform \cite{Tao:10,Almeida:94,Capus:03} (or more generally the local polynomial FT \cite{Katkovnik:95,Stankovic:13}) and the chirplet transform \cite{Mann:95,Mann:91}. The basic idea behind these two modifications is to introduce a rotation in the time-frequency plane by considering $g(u)\rightarrow g(u)e^{i\alpha u^2}$ and $\psi(u)\rightarrow\psi(u)e^{i\alpha u^2}$ in the WFT (\ref{wft}) and WT (\ref{wt}), respectively. This leads to a better representation of chirps if the parameter $\alpha$ is well-matched to the signal's structure while, by allowing this parameter to depend on time $t$ in (\ref{wft}), (\ref{wt}) and choosing it appropriately, one can generally improve the representation for any signal; such a choice, however, is in general highly nontrivial, and is often also computationally expensive \cite{Stankovic:13}. Note that the aforementioned TFRs can be viewed as the WFTs/WTs with modified window/wavelet, so that the signal can still be reconstructed using (\ref{iwft}), (\ref{iwt}); but if the chirp rate $\alpha$ varies for each signal point, then $C_{g,\psi}$ in (\ref{iwft}), (\ref{iwt}) also become time-dependent, while direct reconstruction of $\hat{s}(\xi)$ ceases to be so straightforward.

There are also techniques that allow one to improve certain properties of the TFR by applying some sort of processing to it. For example, the ``sharpness'' of different TFRs -- both linear and quadratic -- can be increased by use of the so-called reassignment methods \cite{Auger:95,Auger:13}, which operate by joining the TFR coefficients at different times and frequencies into a fewer time-frequency entries according to the chosen law. This family includes synchrosqueezing as a particular case, where the reassignment is performed only over the frequency variable but not the time variable, thus retaining the possibility of recovering the time-domain form of the original signal (\ref{istfr}). The conventional reassignment, however, is performed over both time and frequency, which often provides a more concentrated TFR than synchrosqueezing but, at the same time, eliminates the possibility of straightforward signal reconstruction even if the original TFR is linear \cite{Auger:13}.

Finally, many techniques have been proposed for adaptive signal representation, e.g.\ the method of frames \cite{Daubechies:88}, matching pursuit \cite{Mallat:93} (subsequently improved to orthogonal matching pursuit \cite{Pati:93}), best orthogonal basis \cite{Coifman:92} and the basis pursuit (BP) \cite{Chen:01} (for overviews see \cite{Chen:01} and \cite{Mallat:08}). Such methods try to find the sparsest decomposition of a signal in terms of some chosen set of basis functions (the dictionary of time-frequency elements), with BP being the most accurate but, at the same time, most computationally expensive among them. The idea of BP has also inspired the emerging field of compressive sampling \cite{Donoho:06,Candes:06} (for an extensive list of publications in this area see \cite{CSweb}). However, the methods mentioned are more related to discrete TFRs (while here we study continuous ones), and are mainly devoted to fields such as optimal information encoding and compression rather than to estimating the properties of the components present in a signal. Nevertheless, the main problem is that the computational cost of BP is $O(N^3)$, greatly exceeding $O(N\log N)$ of the (S)WFT/(S)WT, while the other methods of this kind, apart from usually being relatively expensive computationally, do not always represent the AM/FM components reliably \cite{Chen:01}.

In addition to the particular TFRs already discussed, there exist many other types and related techniques (see e.g.\ \cite{Stankovic:13if,Stankovic:03,Kwok:00,Kawahara:99,Stankovic:94,Barkat:01}), including those devoted exclusively to some specific goal, such as instantaneous frequency estimation. Each representation, however, has its own issues and advantages/drawbacks, thus deserving a separate review. To avoid an excessive complication, in this work we concentrate solely on the (S)WFT and (S)WT, which we have chosen mainly because they provide an easy way to retrieve the components present in the signal and estimate their characteristics, see Sec.\ \ref{sec:tfsrec} below; but note that they are not the only TFR types to possess this useful property.

\section{Time-, frequency- and time-frequency resolution}\label{sec:tfres}

We have made extensive use of the notions of time, frequency and time-frequency resolutions, although having defined them only briefly. In this section different resolution characteristics are considered in detail, and the definitions are made more precise. Since synchrosqueezing does not seem to change the resolution properties of the transform (see Sec.\ \ref{sec:dwftwt}), in what follows we will concentrate mainly on the WFT and WT, assuming that the the same qualitative and quantitative considerations are applicable respectively to the SWFT and SWT. The validity of such an assumption is demonstrated in Part II, where this issue is studied in more detail.

\subsection{General formulation}\label{sec:tfresA}

Consider a signal $s_{\nu\nu}(t)$ ($s_{\tau\tau}(t)$) consisting of two frequency events -- tones (time events -- delta-peaks), so that its WFT (\ref{wft}) and WT (\ref{wt}) are
\begin{equation}\label{sfreq}
\begin{aligned}
&s_{\nu\nu}(t)=\cos(\nu t)+\cos((\nu+\Delta\nu)t+\Delta\varphi)\\
&\Rightarrow
\left\{\begin{array}{l}
G_s(\omega,t)= \frac{1}{2}\big[\hat{g}(\omega-\nu)+\hat{g}(\omega-\nu-\Delta\nu)e^{i\Delta\nu t}e^{i\Delta\varphi}\big]e^{i\nu t},\\
W_s(\omega,t)= \frac{1}{2}\Big[\hat{\psi}^*\Big(\frac{\omega_\psi\nu}{\omega}\Big)
+\hat{\psi}^*\Big(\frac{\omega_\psi(\nu+\Delta\nu)}{\omega}\Big)e^{i\Delta\nu t}e^{i\Delta\varphi}\Big]
e^{i\nu t},
\end{array}\right.
\end{aligned}
\end{equation}
\begin{equation}\label{stime}
\begin{aligned}
&s_{\tau\tau}(t)=\delta(t-\tau)+\delta(t-\tau-\Delta \tau)e^{i\Delta\varphi}\\
&\Rightarrow
\left\{\begin{array}{l}
G_s(\omega,t)\cong \big[g(\tau-t)+g(\tau+\Delta\tau-t)e^{-i\omega \Delta\tau}e^{i\Delta\varphi}\big]e^{i\omega (t-\tau)},\\
W_s(\omega,t)\cong \frac{\omega}{\omega_\psi}\Big[\psi^*\Big(\frac{\omega(\tau-t)}{\omega_\psi}\Big)
+\psi^*\Big(\frac{\omega(\tau+\Delta\tau-t)}{\omega_\psi}\Big)e^{i\Delta\varphi}\Big].
\end{array}\right.
\end{aligned}
\end{equation}
Everywhere in this section, e.g.\ in (\ref{stime}), the symbol ``$\cong$'' denotes equality up to an error associated with the difference between the WFT/WT (\ref{wft}), (\ref{wt}) and their forms as calculated using the full signal $s(t)$ instead of its positive frequency part $s^{+}(t)$; note that, if redefining WFT/WT to use the full signal, the approximate equality will migrate from $s_{\tau\tau}(t)$ (\ref{stime}) to $s_{\nu\nu}(t)$ (\ref{sfreq}). See \ref{app:posfreq} for a detailed discussion of these issues and the quality of the approximation in (\ref{stime}).

\begin{remark}
Although in this work the original signal is assumed to be real, $s_{\tau\tau}(t)$ (\ref{sfreq}) represents an exception. This is because, without the phase-shift introduced by $e^{i\Delta\phi}$, the analogy between $s_{\tau\tau}(t)$ and $s_{\nu\nu}(t)$ would be incomplete while, intuitively, the time and frequency domains should have equal rights. In fact, since the WFT (\ref{wft}) and WT (\ref{wt}) are based on the positive frequency part of the signal, the resultant TFR will remain the same if instead of $s_{\nu\nu}(t)$ (\ref{sfreq}) one considers $s_{\nu\nu}^+(t)=[e^{i\nu t}+e^{i(\nu+\Delta\nu)t}e^{i\Delta\varphi}]/2$. The latter has the Fourier transform $\hat{s}_{\nu\nu}^+(\xi)=[\delta(\xi-\nu)+\delta(\xi-\nu-\Delta\nu)e^{i\Delta\varphi}]/2$, which is fully analogous to $s_{\tau\tau}(t)$. Nevertheless, it turns out that in all the following considerations nothing changes if one restricts $s_{\tau\tau}(t)$ to be real, which is satisfied for $\Delta\varphi=0,\pi$ in (\ref{stime}). For example, in Sec.\ \ref{sec:tfresD} below, averaging over these two values ($\Delta\varphi=0,\pi$) in the corresponding reconstruction errors will give the same result as averaging over the full range $\Delta\varphi\in[0,2\pi]$.
\end{remark}

It seems reasonable to define the time (frequency) resolution $\gamma_t$ ($\gamma_\omega$) of the transform as the reciprocal of the minimum time (frequency) difference $\Delta\tau$ ($\Delta\nu$) in $s_{\tau\tau}(t)$ ($s_{\nu\nu}(t)$) for which two delta-peaks (tones) can still be reliably resolved in the TFR:
\begin{equation}\label{tfres}
\gamma_t(\nu,\tau)=1/\Delta\tau_{\min}(\nu,\tau),\;\;\gamma_\omega(\nu,\tau)=1/\Delta\nu_{\min}(\nu,\tau).
\end{equation}
However, the meaning of ``reliably resolved'' still remains imprecise, and will be dealt with later. Note, that the definitions (\ref{tfres}) represent the most general case, where the resolutions and minimal differences are ``localized'', i.e.\ allowed to depend on both time and frequency. This is not only useful for the WT, which is characterized by frequency-dependent resolution properties, but also allows one to consider more complicated cases, e.g.\ the WFT with a time-dependent resolution parameter $f_0\rightarrow f_0(t)$.

The joint time-frequency resolution $\gamma_{\omega t}$ is most often defined as the reciprocal of the area of the minimal resolvable square $[\nu,\nu+\Delta\nu_{\min}(\nu,\tau)]$ $\times$ $[\tau,\tau+\Delta\tau_{\min}(\nu,t)]$, being equal to the product of the time and frequency resolutions. However, within such a square $\gamma_t(\nu,\tau)$ and $\gamma_\omega(\nu,\tau)$ can vary considerably, so that these variations should be taken into account to make $\gamma_{\omega t}$ more meaningful. The latter can therefore be defined as
\begin{equation}\label{jtfres}
\gamma_{\omega t}(\nu,\tau)
=\frac{\Delta\tau_{\min}(\nu,\tau)\Delta\nu_{\min}(\nu,\tau)}
{\int_{\nu}^{\nu+\Delta\nu_{\min}(\nu,\tau)}\Delta\tau_{\min}(\omega,\tau)d\omega
\int_{\tau}^{\tau+\Delta\tau_{\min}(\nu,\tau)}\Delta\nu_{\min}(\nu,t)dt}.
\end{equation}
Note that, if neither $\Delta\nu_{\min}$ nor $\Delta\tau_{\min}$ depends on time or frequency, then one has the traditional $\gamma_{\omega t}=\gamma_\omega\gamma_t$.

The definitions (\ref{tfres}) and (\ref{jtfres}) remain valid for any signal representation, not only the WFT and WT. For example, the time domain representation $s(t)$ has infinite time resolution, since, theoretically, two delta-peaks can be resolved for any $\Delta \tau$ in $s_{\tau\tau}(t)$, but zero frequency resolution, as both tones in $s_{\nu\nu}(t)$ have nonzero entries at almost all times, and thus cannot in principle be resolved, no matter how large $\Delta\nu$ is. On the other hand, the frequency domain representation (Fourier transform) $\hat{s}(\xi)$ has infinite frequency resolution, since (in theory) two tones can be perfectly resolved in $\hat{s}_{\nu\nu}(\xi)$ for however small frequency difference $\Delta\nu$; but, at the same time, it has zero time resolution, because the two delta-peaks are spread over all frequencies in $\hat{s}_{\tau\tau}(\xi)$.

TFRs, on the other hand, represent a ``mix'' of the time domain and frequency domain representations. As a result, both time and frequency resolutions are finite for them. This can clearly be seen from the WFT and WT for the two-peak and two-tone signals (\ref{stime}), (\ref{sfreq}). Thus, since the window/wavelet function has non-zero supports in both time and frequency, the peaks in $s_{\tau\tau}(t)$ (tones in $s_{\nu\nu}(t)$) will interfere in the WFT/WT for small enough $\Delta \tau$ ($\Delta\nu$), making their accurate separation and reconstruction impossible (e.g.\ see the components 1,2 in Fig.\ \ref{fig:tfrex}(a)).

From (\ref{sfreq}) and (\ref{stime}) one can see that the interference between two tones or between two delta-peaks in the WFT does not depend on $\nu$ or $\tau$. Therefore, for any meaningful definition of $\Delta \tau_{\min}$ and $\Delta\nu_{\min}$, they also should not depend on time or frequency. Furthermore, it is also clear that rescaling of the WFT window function $g(t)\rightarrow g(t/r)\Leftrightarrow \hat{g}(\xi)\rightarrow \hat{g}(r\xi)$ changes its resolution properties as $\{\Delta\tau_{\min},\Delta\nu_{\min}\}\rightarrow \{r\Delta\tau_{\min},\Delta\nu_{\min}/r\}$. This reflects the repeatedly-mentioned trade-off between time and frequency resolution, being a manifestation of the time-frequency uncertainty principle, which excludes the possibility of simultaneous sharp localization in time and frequency \cite{Kaiser:94,Mallat:08}. Thus, without changing the form of the window function, there is only the possibility of rescaling both resolutions, increasing one and decreasing the other, but not treating them separately. Their product, which in the case of the WFT is equal to the joint time-frequency resolution $\gamma_{\omega t}$ (\ref{jtfres}), remains fixed under such rescaling, thus representing an important characteristic of the window function: the higher it is, the better the trade-off that is possible.

Using the same arguments as for the WFT, from (\ref{stime}) and (\ref{sfreq}) it follows that for the WT both $\Delta \tau_{\min}$ and $\Delta\nu_{\min}$ depend on frequency, but not time, so that $\gamma_t=\gamma_t(\nu)$, $\gamma_\omega=\gamma_\omega(\nu)$. Next, one also has $\Delta_\tau(r\nu)=\Delta_\tau(\nu)/r$ and $\Delta\nu(r\nu)=r\Delta\nu(\nu)$, so that the time resolution of the WT increases with frequency, while the frequency resolution decreases. Nevertheless, as will be seen, the time-frequency resolution (\ref{jtfres}) does not depend on $\nu$ or $\tau$, being fixed for the specified wavelet parameters.

As mentioned above (see Remark \ref{rem:reschwt}), the direct rescaling $\{\psi(t)$, $\hat{\psi}(\xi)\} \rightarrow \{\psi(t/r), \hat{\psi}(r\xi)\}$ does not change anything for the WT, and to tune the resolutions at a particular frequency one needs to change the wavelet peak frequency $\omega_\psi\rightarrow r\omega_\psi$ while preserving the spreads and the forms of $\{\psi(t),\hat{\psi}(\xi)\}$ (or vice versa). If such a procedure was possible, it would lead to the same trade-off rule as for the WFT, i.e.\ $\gamma_t(\nu)\rightarrow \gamma_t(\nu)/r$, $\gamma_{\omega}(\nu)\rightarrow r\gamma_\omega(\nu)$, with $\gamma_{\omega t}$ remaining fixed. However, given restrictions such as the admissibility condition (\ref{adc}), changing the wavelet peak frequency will inevitably affect its form and/or its spread (on the linear scale). Thus, e.g.\ for the Morlet wavelet (\ref{mw}) a change of $\omega_\psi$ (by varying $f_0$) will be accompanied by a simultaneous change of its form because of the admissibility term; for the lognormal wavelet, changing $f_0$ will change its spread in frequency on the logarithmic scale, but on the linear scale $\hat{\psi}(\xi)$ will become more asymmetric, affecting the time domain form $\psi(t)$. Because of these issues, there does not seem to be a way of changing the resolution properties of the WT without altering its time-frequency resolution, in contrast to the WFT. Furthermore, $\gamma_{\omega t}$ (whichever way defined) progressively worsens with an increase of the wavelet time resolution (decrease of $f_0$), and it seems in principle impossible to reach a very sharp time-localization in the WT (see \ref{app:winwav}).

\subsection{Classical definitions and their flaws}\label{sec:tfresB}

Although $\Delta \tau_{\min}$ and $\Delta\nu_{\min}$ in (\ref{tfres}) were defined respectively as the minimum time and frequency difference which can reliably be resolved in the TFR, the meaning of ``reliably resolved'' remains mathematically unclear. Based on how it is defined, one can characterize the resolution properties of the TFRs in different ways.

The traditional approach \cite{Kaiser:94,Mallat:08} implicitly assumes two delta-peaks (tones) to be well-resolved if the time (frequency) distance between them exceeds some number of standard deviations of the squared window/wavelet function in the time (frequency) domain. Within this framework, for the WFT one has
\begin{equation}\label{clreswft}
\begin{gathered}
\begin{aligned}
&\Delta\nu_{\min}^{(cl)}=k_1\Delta_\omega,\quad
\Delta_\omega^2=E_g^{-1}\frac{1}{2\pi}\int (\omega-\omega_c)^2|\hat{g}(\omega)|^2d\omega,\\
&\Delta\tau_{\min}^{(cl)}=k_2\Delta_t,\;\quad
\Delta_t^2=E_g^{-1}\int (t-t_c)^2|g(t)|^2dt,
\end{aligned}\\
\omega_c \equiv E_g^{-1}\frac{1}{2\pi}\int\omega|\hat{g}(\omega)|^2d\omega,
\quad t_c \equiv E_g^{-1}\int t|g(t)|^2dt,\\
E_g \equiv \frac{1}{2\pi}\int |\hat{g}(\xi)|^2d\xi=\int |g(t)|^2dt\quad \mbox{(by Parseval's identity)},\\
\end{gathered}
\end{equation}
where $k_{1,2}$ are implicitly assumed to be the same for all window functions. Obviously, such a definition is far from universal, since for different windows different number of standard deviations are needed to resolve the two tones/delta-peaks. As an illustrative example, the window with asymptotics $|\hat{g}(\xi\rightarrow\pm\infty)|\sim[1+|\xi|^{5/4}]^{-1}$ has infinite variance of $|\hat{g}(\xi)|^2$, but can still be used and allows for an accurate resolution and reconstruction of the two tones for high enough frequency difference between them.

For the WT, the classic variance-based framework takes the form
\begin{equation}\label{clreswt}
\begin{gathered}
\begin{aligned}
&\Delta\nu_{\min}^{(cl)}(\nu)=k_1\frac{\nu}{\omega_\psi}\Delta_\omega,\quad
\Delta_\omega^2=E_\psi^{-1}\frac{1}{2\pi}\int (\omega-\omega_c)^2|\hat{\psi}(\omega)|^2d\omega,\\
&\Delta\tau_{\min}^{(cl)}(\nu)=k_2\frac{\omega_\psi}{\nu}\Delta_t,\;\quad
\Delta_t^2=E_\psi^{-1}\int (t-t_c)^2|\psi(t)|^2dt,
\end{aligned}\\
\omega_c \equiv E_\psi^{-1}\frac{1}{2\pi}\int\omega|\hat{\psi}(\omega)|^2d\omega,\quad
t_c \equiv E_\psi^{-1}\int t|\psi(t)|^2dt,\\
E_\psi \equiv \frac{1}{2\pi}\int |\hat{\psi}(\xi)|^2d\xi=\int |\psi(t)|^2dt\quad
\mbox{(by Parseval's identity)},\\
\end{gathered}
\end{equation}
where $k_{1,2}$ are again implicitly assumed to be the same for all wavelet functions. This approach has the same drawbacks as (\ref{clreswft}) for the WFT. However, in the case of the WT it is actually not appropriate at all, at least in terms of the frequency resolution. Thus, as can be seen from (\ref{Nwt}), the tones are represented in the WT as terms $\sim\hat{\psi}(\omega_\psi\nu/\omega)$, so that the decay of their contribution as $\omega\rightarrow\infty$, determined by the behavior of $\hat{\psi}(\xi)$ as $\xi\rightarrow 0$, will obviously have a big effect on the frequency resolution. At the same time, the usual variance $\Delta_\omega^2$ (\ref{clreswt}) takes no account of this fact, e.g.\ being invariant under $\hat{\psi}(\xi)\rightarrow\hat{\psi}(\xi+\omega_\psi)$, that makes the wavelet inadmissible (in which case tones that are infinitely distant in frequency still interfere, so that the frequency resolution becomes effectively zero). Therefore, for wavelets, it seems more appropriate to study at least the variance of $|\hat{\psi}(\omega_\psi/\xi)|^2$, but by no means that of $|\hat{\psi}(\xi)|^2$.

For both WFT and WT, the classic time-frequency resolution measure is taken as $\gamma_{\omega t}^{(cl)}=[\Delta_\omega\Delta_t]^{-1}$, with $\Delta_\omega,\Delta_t$ being given by (\ref{clreswft}) for the WFT and by (\ref{clreswt}) for the WT (note also the difference between $\gamma_{\omega t}^{(cl)}$ and (\ref{jtfres}) for the latter). It can be shown \cite{Mallat:08,Kaiser:94,Addison:10}, that this measure attains its maximum for the Gaussian window (\ref{gw}) and (up to the effect of the admissibility term $\sim e^{-(2\pi f_0)^2/2}$) for the Morlet wavelet (\ref{mw}). However, as follows from the discussion above, only in the case of the WFT does the classic $\gamma_{\omega t}^{(cl)}$ make some sense, though even in this case it remains highly non-universal.

\subsection{Notion of the window/wavelet $\epsilon$-support}\label{sec:tfresC}

Before proceeding to a reconsideration of the classic definitions, it is useful to introduce the notions of the window/wavelet $\epsilon$-supports in frequency $[\xi_1(\epsilon),\xi_2(\epsilon)]$ and time $[\tau_1(\epsilon),\tau_2(\epsilon)]$, which will be used frequently in what follows. These $\epsilon$-supports are defined as the widest intervals containing the $(1-\epsilon)$ part of the total integrals of the window/wavelet function which appear in $C_{g,\psi}$ and $\widetilde{C}_{g,\psi}$ (\ref{iwft}), (\ref{iwt}). As will be seen below, they are directly related to the accuracy with which the components can be recovered from the WFT/WT, and thus can be used effectively for quantifying it.

Considering first the WFT, for an arbitrary window function, including functions that are not always positive and can be oscillating or complex, the corresponding definitions are
\begin{equation}\label{eswft}
\begin{aligned}
&R_g(\omega)\equiv \frac{\int_{-\infty}^\omega\hat{g}(\xi)d\xi}{\int \hat{g}(\xi)d\xi}
=C_g^{-1}\frac{1}{2}\int_{-\infty}^\omega\hat{g}(\xi)d\xi,\\
&\xi_{1,2}(\epsilon):\;|R_g(\xi\leq\xi_1)|\leq\epsilon/2,\;|1-R_g(\xi\geq\xi_2)|\leq\epsilon/2,\\
&P_g(\tau)\equiv \frac{\int_{-\infty}^\tau g(t)dt}{\int g(t)dt}
=\widetilde{C}_g^{-1}\int_{-\infty}^\tau g(t)dt,\\
&\tau_{1,2}(\epsilon):\;|P_g(\tau\leq\tau_1)|\leq\epsilon/2,\;|1-P_g(\tau\geq\tau_2)|\leq\epsilon/2.\\
\end{aligned}
\end{equation}
Evidently, $|R_g(\omega)|$ and $|1-R_g(\omega)|$ quantify the relative parts of $\hat{g}(\xi)$ that are contained in the ranges $\xi<\omega$ and $\xi>\omega$, respectively, while the values $\xi_{1,2}(\epsilon)$ specify the limits within which the $(1-\epsilon)$ part of the window FT resides. In the same manner, $|P_g(\tau)|$ and $|1-P_g(\tau)|$ reflect the relative parts of $g(t)$ contained in the ranges $t<\tau$ and $t>\tau$, respectively, while $[\tau_1(\epsilon),\tau_2(\epsilon)]$ represents the region encompassing its $(1-\epsilon)$ part. The inequalities in the definitions of $\xi_{1,2}(\epsilon)$ ($\tau_{1,2}(\epsilon)$) are needed only if $\hat{g}(\xi)$ ($g(t)$) is not strictly positive, or complex, to ensure that the integral of the latter over any frequency (time) region containing the $\epsilon$-support $[\xi_1(\epsilon),\xi_2(\epsilon)]$ ($[\tau_1(\epsilon),\tau_2(\epsilon)]$) will always approximate the corresponding full integral with relative error not higher than $\epsilon$.

Considering the single tone $s(t)=\cos\nu t\Rightarrow G_s(\omega,t)=\hat{g}(\omega-\nu)e^{i\nu t}/2$, it is clear that its WFT at frequencies $\omega<\xi$ ($\omega>\xi$) will contain the $|R_g(\xi-\nu)|$ ($|1-R_g(\xi-\nu)|$) part of the signal. Furthermore, its $(1-\epsilon)$ part will be contained in the frequency range $[\nu+\xi_1(\epsilon),\nu+\xi_2(\epsilon)]$, so that e.g.\ for real $\hat{g}(\xi)$ the corresponding signal reconstructed by (\ref{iwft}) from the WFT in this range will be $(1-\epsilon)\cos\nu t$, representing exactly the $(1-\epsilon)$ part of the original tone.

Likewise, for the delta-peak $s(t)=\delta(t-\tau)\Rightarrow G_s(\omega,t)\cong g(\tau-t)e^{i\omega (t-\tau)}$ the WFT at times $t<t_0$ ($t>t_0$) will contain the $\cong|1-P_g(\tau-t_0)|$ ($\cong|P_g(\tau-t_0)|$) part of the signal, while its $(1-\epsilon)$ part will be contained in the time interval $\cong[\tau-\tau_2(\epsilon),\tau-\tau_1(\epsilon)]$. Thus, e.g.\ for real $g(t)$ the delta-peak's FT reconstructed by (\ref{iwft}) from this interval will be $\hat{s}(\xi)\cong(1-\epsilon)e^{-i\xi \tau}$.

Similarly to the case of the WFT, the $\epsilon$-supports for the WT are defined based on (\ref{iwt}) as
\begin{equation}\label{eswt}
\begin{aligned}
&R_\psi(\omega)\equiv \frac{\int_0^\omega\hat{\psi}^*(\xi)\frac{d\xi}{\xi}}{\int_0^\infty\hat{\psi}^*(\xi)\frac{d\xi}{\xi}}
=C_\psi^{-1}\frac{1}{2}\int_0^\omega\hat{\psi}^*(\xi)\frac{d\xi}{\xi},\\
&\xi_{1,2}(\epsilon):\;|R_\psi(\xi\leq\xi_1)|\leq\epsilon/2,\;|1-R_\psi(\xi\geq\xi_2)|\leq\epsilon/2,\\
&P_\psi(\tau)\equiv \frac{\int_{-\infty}^\tau \psi^*(t)e^{i\omega_\psi t}dt}{\int \psi^*(t)e^{i\omega_\psi t}dt}
=\widetilde{C}_\psi^{-1}\int_{-\infty}^\tau \psi^*(t)e^{i\omega_\psi t}dt,\\
&\tau_{1,2}(\epsilon):\;|P_\psi(\tau\leq\tau_1)|\leq\epsilon/2,\;|1-P\psi(\tau\geq\tau_2)|\leq\epsilon/2.\\
\end{aligned}
\end{equation}
Like $|P_g(\tau)|$ in (\ref{eswt}), $|P_\psi(\tau)|$ ($|1-P_\psi(\tau)|$) quantifies the relative part of $\psi(t)e^{-i\omega_\psi t}$ contained at $t<\tau$ ($t>\tau$), with $[\tau_1(\epsilon),\tau_2(\epsilon)]$ specifying the interval encompassing its $(1-\epsilon)$ part. In the same manner, $|R_\psi(\omega)|$ and $\xi_{1,2}(\epsilon)$ are related to the relative part of $\hat{\psi}(\xi)$, taken on a logarithmic scale.

However, due to the scaling nature of the WT, the relationships of (\ref{eswt}) to real quantities differ slightly from the case of the WFT. Thus, for the single tone $s(t)=\cos(\nu t)\Rightarrow W_s(\omega,t)=\hat{\psi}^*(\omega_\psi\nu/\omega)e^{i\nu t}/2$ the WT at frequencies $\omega<\xi$ ($\omega>\xi$) will contain the $|1-R_\psi(\omega_\psi\nu/\xi)|$ ($|R_\psi(\omega_\psi\nu/\xi)|$) part of the signal, while its $(1-\epsilon)$ part will lie in the band $[\omega_\psi\nu/\xi_2(\epsilon),\omega_\psi\nu/\xi_1(\epsilon)]$.

For the delta-function $s(t)=\delta(t-\tau)\Rightarrow W_s(\omega,t)\cong \frac{\omega}{\omega_\psi}\psi^*\Big(\frac{\omega(\tau-t)}{\omega_\psi}\Big)$ the WT spread in time will vary for different $\omega$. At each frequency the part of the delta-function's total FT contained in the WT at $t<t_0$ ($t>t_0$) will be $\cong |1-P_\psi(\omega(\tau-t_0)/\omega_\psi)|$ ($\cong|P_\psi(\omega(\tau-t_0)/\omega_\psi)|$), while its $(1-\epsilon)$ part will reside in the interval $\cong [\tau-\omega_\psi \tau_2(\epsilon)/\omega,\tau-\omega_\psi\tau_1(\epsilon)/\omega]$.

The quantities (\ref{eswft}), (\ref{eswt}) are very convenient and will be used extensively below,
not only in the present section. For simplicity, $\tau_{1,2}(\epsilon)$ and $\xi_{1,2}(\epsilon)$ denote the respective $\epsilon$-supports both for the window function in the WFT and for the wavelet function in the WT. The meaning will always be clear from the context. Note that the full supports of the window/wavelet in time ($g(t),\psi(t)$) and frequency ($\hat{g}(\xi),\hat{\psi}(\xi)$), whether finite or not, are $[\tau_1(0),\tau_2(0)]$ and $[\xi_1(0),\xi_2(0)]$, respectively.

\subsection{Reconsidered definitions}\label{sec:tfresD}

A more universal and appropriate approach (than the traditional variance-based one) is to regard two components as being reliably resolved if they can each be accurately identified and reconstructed from the signal's TFR. i.e.\ can be recovered with a relative error not exceeding some threshold. Consider the WFT of the two-tone signal (\ref{sfreq}), from which one wants to find the individual analytic signals $x_{\nu\nu;1}^a(t)=e^{i\nu t},x_{\nu\nu;2}^a(t)=e^{i(\nu+\Delta\nu) t}e^{i\Delta\varphi}$ for each of the two tones.
At any time $t$, this can be done by first dividing the frequency range at some $\omega=\omega_{\rm x}(t)$ into two parts, each responsible for a separate tone, and then integrating the WFT over the corresponding frequency ranges in the same way as in (\ref{iwft}). This will give the reconstructed analytic signals $\tilde{x}_{\nu\nu;1,2}^{a}(t)$ which, using (\ref{sfreq}) and (\ref{eswft}), can be represented as
\begin{equation}\label{tonerec}
\begin{aligned}
&\tilde{x}_{\nu\nu;1}^{a}(t)=\int_{-\infty}^{\omega_{\rm x}(t)}G_s(\omega,t)d\omega\\
&=\frac{C_g^{-1}}{2}e^{i\nu t}\left[\int_{-\infty}^{\omega_{\rm x}(t)}\hat{g}(\omega-\nu)d\omega
+e^{i(\Delta\nu t+\Delta\varphi)}\int_{-\infty}^{\omega_{\rm x}(t)}\hat{g}(\omega-\nu-\Delta\nu)d\omega\right]\\
&=e^{i\nu t}\left[\big(1-R_g(\nu-\omega_{\rm x}(t))\big)+R_g(\omega_{\rm x}(t)-\nu-\Delta\nu)e^{i(\Delta\nu t+\Delta\varphi)}\right],\\
&\tilde{x}_{\nu\nu;2}^{a}(t)=\int_{\omega_{\rm x}(t)}^{\infty}G_s(\omega,t)d\omega\\
&=\frac{C_g^{-1}}{2}e^{i\nu t}\left[\int_{\omega_{\rm x}(t)}^{\infty}\hat{g}(\omega-\nu)d\omega
+e^{i(\Delta\nu t+\Delta\varphi)}\int_{\omega_{\rm x}(t)}^{\infty}\hat{g}(\omega-\nu-\Delta\nu)d\omega\right]\\
&=e^{i\nu t}\left[R_g(\nu-\omega_{\rm x}(t))+\big(1-R_g(\omega_{\rm x}(t)-\nu-\Delta\nu)\big)e^{i(\Delta\nu t+\Delta\varphi)}\right],
\end{aligned}
\end{equation}
where $R_g(x)=1-R_g(-x)$ is defined in (\ref{eswft}).

Obviously, the reconstruction errors $x_{\nu\nu;1,2}^{a}(t)-\tilde{x}_{\nu\nu;1,2}^{a}(t)$ generally depend on the phase-shift $\Delta\varphi$. Therefore, in the corresponding expressions one should take the average over $\Delta\varphi$, which will be denoted as $\langle...\rangle_{\Delta\varphi}$. The relative errors of each tone's reconstruction $\varepsilon_{\nu\nu;1,2}(\nu,t,\Delta\nu)$ then become
\begin{equation}\label{tonerecerr}
\begin{aligned}
\varepsilon_{\nu\nu;1,2}^2(\nu,t,\Delta\nu)\equiv&
\frac{\langle |x_{\nu\nu;1,2}^{a}(t)-\widetilde{x}_{\nu\nu;1,2}^{a}(t)|^2\rangle_{\Delta\varphi}}
{\langle |x_{\nu\nu;1,2}^{a}(t)|^2\rangle_{\Delta\varphi}}\\
=&|R_g(\nu-\omega_{\rm x}(t))|^2+|R_g(\omega_{\rm x}(t)-\nu-\Delta\nu)|^2.\\
\end{aligned}
\end{equation}
Note that, in the present case, averaging over $\Delta\varphi$ and time-averaging will give the same results; however, in general the TFR resolution properties can depend on time, and taking the mean over phase-shifts allows one to localize these errors at each $t$.

The minimum resolvable frequency difference $\Delta\nu_{\min}(\nu,t)$ can be defined as the minimum $\Delta\nu$ in (\ref{tonerecerr}) for which the total error $\varepsilon_{\nu\nu}(\nu,t,\Delta\nu)$ is still smaller than some accuracy threshold $\epsilon_r$:
\begin{equation}\label{fdmin}
\begin{aligned}
&\Delta\nu\geq\Delta\nu_{\min}(\nu,t):\\
&\varepsilon_{\nu\nu}(\nu,t,\Delta\nu)\equiv
\left[\varepsilon_{\nu\nu;1}^2(\nu,t,\Delta\nu)+\varepsilon_{\nu\nu;2}^2(\nu,t,\Delta\nu)\right]^{1/2}\leq \epsilon_r.
\end{aligned}
\end{equation}
It can be expressed through the $\epsilon_r$-support of the window in frequency (\ref{eswft}). Thus, consider the WFT with real, positive and symmetric $\hat{g}(\omega)$, e.g.\ a Gaussian (\ref{gw}). Then it follows from (\ref{sfreq}) that the minimum WFT amplitude between the peaks corresponding to two tones will always appear at $\omega=\nu+\Delta\nu/2$ (unless these two peaks are merged into a single one at some times, which might happen if $\Delta\nu$ is too small). Therefore, in practice the respective frequency regions of the tones will be separated exactly at their average frequency (see Sec.\ \ref{sec:tfsrec} below), so that one should use $\omega_{\rm x}(t)=\nu+\Delta\nu/2$ when estimating the errors (\ref{tonerecerr}); it can also be shown that, in the present case, these errors are minimized by such a choice of $\omega_{\rm x}$. The overall reconstruction error is then
\begin{equation}\label{overr}
\varepsilon_{\nu\nu}(\Delta\nu)=\big[2|R_g(-\Delta\nu/2)|^2+2|R_g(-\Delta\nu/2)|^2\big]^{1/2}=2|R_g(-\Delta\nu/2)|,
\end{equation}
and, taking into account that $\xi_1(\epsilon)=-\xi_2(\epsilon)$ due to the assumed window symmetry, it follows from (\ref{fdmin}),(\ref{eswft}) that the frequency difference for which two tones are recovered with inaccuracy $\epsilon$ is exactly equal to the $\epsilon$-support of $\hat{g}(\xi)$. For other window forms (e.g. asymmetric $\hat{g}(\xi)$) all becomes more complicated, but one can still expect to get an overall error of around $\epsilon$ when $\Delta\nu=\xi_2(\epsilon)-\xi_1(\epsilon)$. Note that the above considerations hold for reasonably small $\epsilon$, so that $\Delta\nu$ is high enough and there are always two distinct peaks in the WFT amplitude; otherwise, if the peaks are merged at certain times, the actual reconstruction errors will be larger than (\ref{tonerecerr}).

The case of two delta-peaks (\ref{stime}) is closely similar to that of two tones, so the same considerations apply, with just $\xi_{1,2}(\epsilon)\rightarrow\tau_{1,2}(\epsilon)$. Hence, setting $\epsilon_r$ as the maximum allowable reconstruction error for which two tones/delta-peaks can still be regarded as resolved, the minimum resolvable time-delay $\Delta\tau_{\min}$ and frequency difference $\Delta\nu_{\min}$, and the other resolution parameters based on them, for the WFT take the forms
\begin{equation}\label{tfreswft}
\begin{gathered}
\Delta\nu_{\min}=\xi_2(\epsilon_r)-\xi_1(\epsilon_r),\quad
\Delta\tau_{\min}\cong\tau_2(\epsilon_r)-\tau_1(\epsilon_r),\\
\gamma_\omega\equiv\Delta\nu_{\min}^{-1},\quad \gamma_t=\Delta\tau_{\min}^{-1},\quad
\gamma_{\omega t}\equiv\gamma_\omega\gamma_t=[\Delta\nu_{\min}\Delta\tau_{\min}]^{-1},\\
\end{gathered}
\end{equation}
where $\tau_{1,2}(\epsilon)$ and $\xi_{1,2}(\epsilon)$ are defined in (\ref{eswft}). Generally, accurate reconstruction might reasonably be assumed as being at 95\% precision, so one can set $\epsilon_r=0.05$ in (\ref{tfreswft}). The resolution characteristics of different windows are listed in \ref{app:winwav}.

The same approach straightforwardly extends to the WT case, where one applies similar considerations in terms of (\ref{eswt}). Thus, it can be shown that two tones with frequency ratio $\frac{\nu+\Delta\nu}{\nu}=1+\Delta\nu/\nu=\frac{\xi_2(\epsilon)}{\xi_1(\epsilon)}$ are reconstructed from the WT with an overall relative error of around $\epsilon$. This estimate is exact if $\hat{\psi}(\xi)$ is real, positive and symmetric on a logarithmic scale (such as the lognormal wavelet (\ref{lw})), in which case the tones will always be separated at $\omega_{\rm x}(t)=\exp[\log\nu+\log(\nu+\Delta\nu)]=\sqrt{\nu(\nu+\Delta\nu)}$. For the resolution of two delta-peaks, it follows from (\ref{stime}) that one should consider the $\epsilon$-supports corresponding to $\psi^*(\omega t/\omega_\psi)$, so that the related error will be different at each frequency $\omega$, characterized by $\epsilon$ calculated from $\omega\Delta\tau/\omega_\psi=\tau_2(\epsilon)-\tau_1(\epsilon)$. Hence, the resolution parameters (\ref{tfres}), (\ref{jtfres}) for the WT are
\begin{equation}\label{tfreswt}
\begin{gathered}
\Delta\nu_{\min}(\nu)=\nu{\Big(}\frac{\xi_2(\epsilon_r)}{\xi_1(\epsilon_r)}-1{\Big)},\quad
\Delta\tau_{\min}(\nu)\cong\frac{\omega_\psi}{\nu}\big(\tau_2(\epsilon_r)-\tau_1(\epsilon_r)\big),\\
\gamma_\omega(\nu)\equiv[\Delta\nu_{\min}(\nu)]^{-1},\quad \gamma_t(\nu)=[\Delta\tau_{\min}(\nu)]^{-1},\\
\gamma_{\omega t}\cong
\Big[\omega_\psi\big(\tau_2(\epsilon_r)-\tau_1(\epsilon_r)\big)\log\frac{\xi_2(\epsilon_r)}{\xi_1(\epsilon_r)}\Big]^{-1},\\
\end{gathered}
\end{equation}
where $\tau_{1,2}(\epsilon)$ and $\xi_{1,2}(\epsilon)$ are defined in (\ref{eswt}), and $\epsilon_r$ denotes the maximum allowable reconstruction error, which can be set to $\epsilon_r=0.05$, similarly to that in (\ref{tfreswft}). The resolution characteristics of different wavelets are listed in \ref{app:winwav}.

Summarizing, in contrast to the classic resolution measures (\ref{clreswft}), (\ref{clreswt}), the quantities in (\ref{tfreswft}), (\ref{tfreswt}) are very universal and have clear physical meaning, being directly related to the accuracy with which two time or frequency events can be recovered from the resultant TFR.

\section{Extraction of components from the TFR}\label{sec:tfsrec}

The main purpose of time-frequency analysis can be formulated as the identification and quantification of the AM/FM components present in a signal. Being reliably represented in the TFR, they can be identified and reconstructed. In this respect the TFR can be used to decompose a signal into its constituent components, or just to recover some particular components of interest. The advantages of this approach are: (i) its noise robustness; and (ii) the possibility of adjusting the corresponding time and frequency resolutions (by varying the window/wavelet resolution parameter), i.e.\ specifying the trade-off between the restrictions on the frequency separation and time-modulation of the components (see Part II). Thus, it can easily be checked numerically that the TFR-based decomposition is less susceptible to noise than e.g.\ the (ensemble) empirical mode decomposition method \cite{Huang:98,Wu:09}; it is also known \cite{Feldman:09,Rilling:08} that the frequency resolution of the latter (which is logarithmic, like that of the WT) can only be adjusted slightly by varying the number of sifting iterations. We now discuss the two steps of component extraction from the TFR, namely its identification and reconstruction.

\begin{remark}
Note that the extraction of components from the ``continuous'' WT (\ref{wt}) represents a completely different (usually more accurate because of being adaptive) approach from subband filtering -- the other wavelet-based decomposition, corresponding to the ``discrete'' WT \cite{Vetterli:95,Akansu:00}.
\end{remark}

\subsection{Identification of the components}\label{sec:tfsrecA}

\begin{figure*}[t!]
\includegraphics[width=1.0\linewidth]{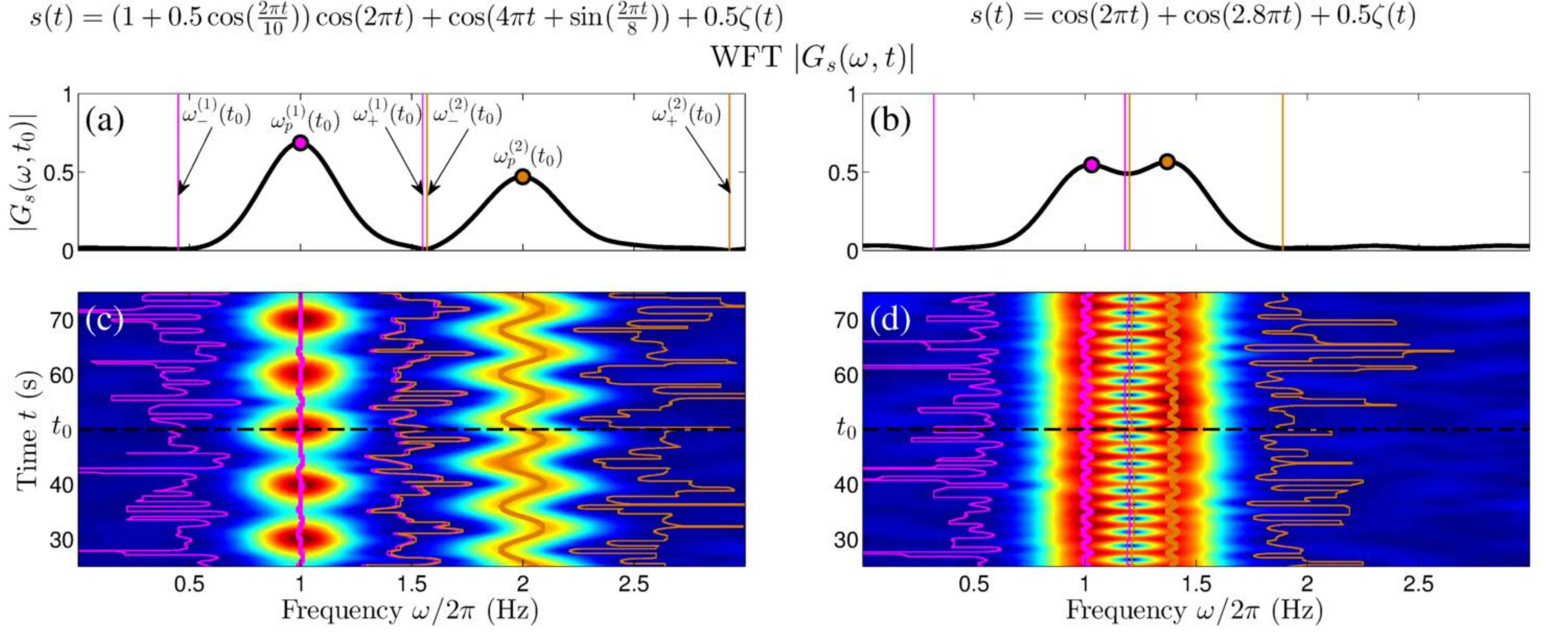}\\
\includegraphics[width=1.0\linewidth]{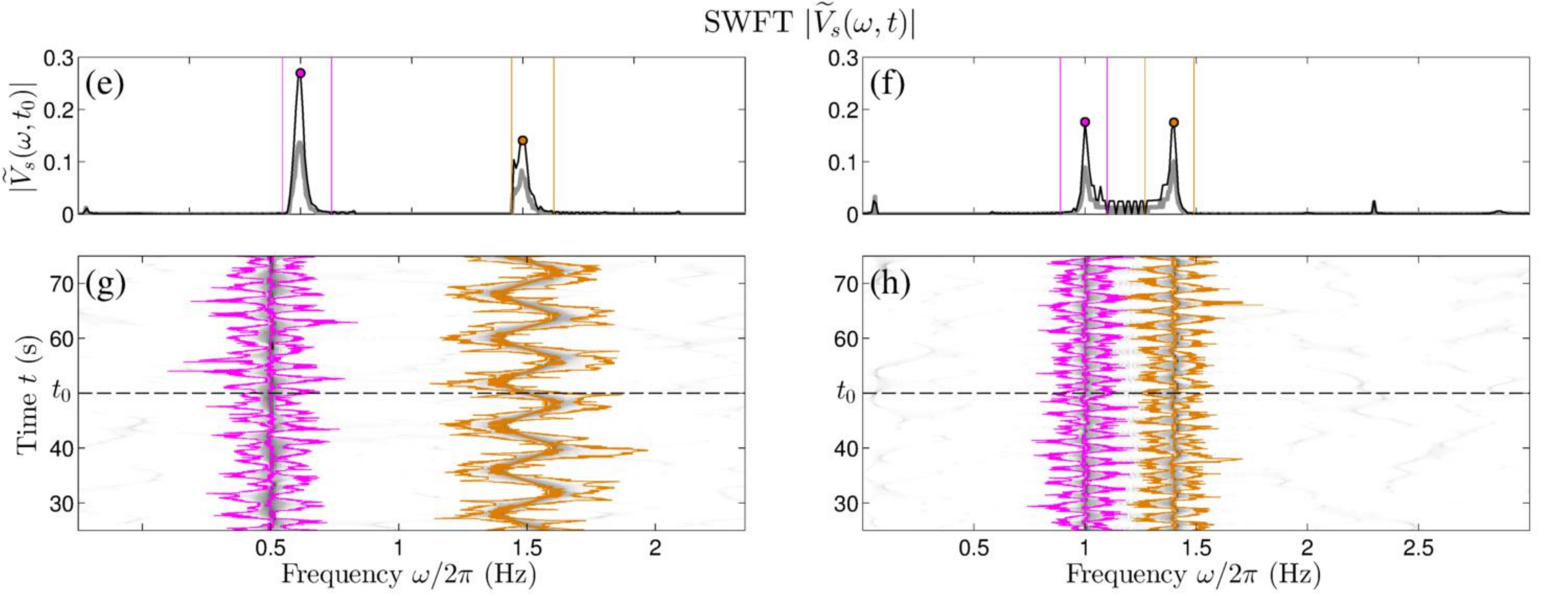}\\
\caption{Examples of the extracted components' ridge curves $\omega_p(t)$ and time-frequency supports (TFSs) $[\omega_-(t),\omega_+(t)]$ in the WFT and SWFT (the case of the WT and SWT is similar) for two different signals, which are defined by the equations above (a) and (b); in these equations, $\zeta(t)$ denotes unit-deviation Gaussian white noise. Each signal consists of two components, and in each panel the colors magenta and orange refer to the first and second components, respectively. (a,b): Snapshots of the WFT amplitudes (thick black lines) at time $t_0=50$ s, the borders of the components' TFSs at that time (thin colored lines) and the associated ridge points $\omega_p(t)$ (filled colored circles). (c,d): Full WFT amplitude time-evolutions ($t_0$, corresponding to the snapshots in (a,b), are indicated by the dashed black lines), with the components' TFSs $[\omega_-^{(1,2)}(t),\omega_+^{(1,2)}(t)]$ being indicated by thin colored lines, and their ridge curves $\omega_p^{(1,2)}(t)$ by thick solid lines of the same colors. (e,f,g,h): Same as (a,b,c,d), but for the SWFT amplitude. The gray background lines in (e,f) show the amplitudes of the SWFTs calculated with a frequency discretization step $\Delta\omega/2\pi=0.005$, which is half that of the original. The signals were sampled at 100 Hz for 100 s.}
\label{fig:tfsupp}
\end{figure*}

AM/FM components visually appear as ``curves'' in the TFR plots (see e.g.\ Fig.\ \ref{fig:tfrex}). Thus, if the construction of the (S)WFT or (S)WT is well-matched to the signal's structure, then each component will be represented by a unique sequence of TFR amplitude peaks -- the ridge curve $\omega_p(t)$, being mapped to a particular time-frequency region -- the time-frequency support (TFS) $[\omega_-(t),\omega_+(t)]$. The latter can be defined as the widest region of unimodal and non-zero (WFT/WT) or just non-zero (SWFT/SWT) TFR amplitude around $\omega_p(t)$:

\begin{definition}
The time-frequency support $[\omega_-(t),\omega_+(t)]$ of the component in the current TFR $Y_s(\omega,t)$ (WFT, WT, SWFT or SWT) is defined as the widest frequency region around the corresponding ridge curve $\omega_p(t)$ where:
\begin{enumerate}
\item The TFR amplitude is higher than zero:\newline
$|Y_s(\omega,t)|>0,\;\forall t,\omega\in[\omega_-(t),\omega_+(t)]$.
\item The TFR amplitude is unimodal:\newline
${\rm sign}(\omega-\omega_p(t))\partial_\omega |Y_s(\omega,t)|\leq0,\;\forall t,\omega\in[\omega_-(t),\omega_+(t)]$\newline
(smooth TFRs only, i.e.\ excluding SWFT and SWT).
\end{enumerate}
Conversely, the ridge curve is then the sequence of highest TFR amplitude peaks in the corresponding TFS at each time: $$\omega_p(t)=\underset{\omega\in[\omega_-(t),\omega_+(t)]}{\argmax}|Y_s(\omega,t)|.$$
Examples of ridge curves and TFSs for different components are shown in Fig.\ \ref{fig:tfsupp}.
\end{definition}

In the TFS definition, the first criterion establishes that the support is cut when the TFR amplitude drops to zero, usually implying that the region to which some AM/FM component is mapped has ended. Considering the SWFT/SWT, there will always be regions (usually large) where the TFR amplitude is exactly zero (see Fig.\ \ref{fig:tfsupp}(e,f)). Indeed, if one is using the same $N_f$ frequency bins for the SWFT/SWT as for the WFT/WT from which it is constructed then, in the process of synchrosqueezing, at each time one maps $N_f$ values of $\nu_{G,W}(\omega,t)$ onto $N_f$ bins. In practice many values will often be assigned to one bin, and in effect the components will be separated by regions of zero amplitude. Although harder to prove, this seems to be the case even when using more bins for the WFT/WT than for the SWFT/SWT. 

The second criterion in the TFS definition establishes that, if two supports are overlapping in the TFR (i.e.\ the corresponding components interfere with each other), then they will be nearly-optimally separated at frequency where the amplitude minimum between the two peaks occurs, indicating that one component's contribution has been overcome by that of the other. However, such treatment is suitable only for smooth TFRs, such as the WFT/WT. On the contrary, the synchrosqueezed TFRs are non-smooth, sometimes having saw-tooth patterns (see Fig.\ \ref{fig:tfsupp} (e,f)). Therefore, for the SWFT/SWT we cannot require strict unimodality, and the TFS is extracted as the widest region of non-zero amplitude around the given peaks. Neglecting the second criterion might seem to cause interfering components to be assigned to the same support (as would be the case for the WFT/WT). However, for the SWFT/SWT this does not appear to happen. We have found numerically that if the two interfering components can be distinguished in the WFT/WT (i.e.\ are not merged completely into one peak), then synchrosqueezing will map most of their power onto the two SWFT/SWT supports, separated by regions of zero amplitude (see Fig.\ \ref{fig:tfsupp}(f)).

\begin{remark}
It should be noted, that the TFS definition for the WFT/WT (but not SWFT/SWT) as presented here is inappropriate for windows and wavelets with multimodal (i.e.\ having multiple non-negligible peaks) $|\hat{g}(\xi)|$ and $|\hat{\psi}(\xi>0)|$, respectively. Thus, for such window/wavelet functions it might be very hard to understand where is the ``line'' separating the time-frequency regions corresponding to different components, because even for the single tone signal $s(t)\sim\cos\nu t$ there will be multiple peaks in the TFR amplitude. Due to this, windows/wavelets not satisfying approximate unimodality (in the frequency domain) enormously complicate the interpretation of the resultant WFT/WT, which is why they are rarely used.
\end{remark}

In general, the behavior of the synchrosqueezed TFRs is much more complex than that of the WFT/WT and, even in the noise-free case, the SWFT/SWT will contain a lot of spurious TFSs represented by small amplitude spikes. This can be seen in Fig.\ \ref{fig:tfsupp}(f), where two main supports are separated by many smaller ones. Thus, although most parts of the two well-defined TFSs in the WFT (Fig.\ \ref{fig:tfsupp}(b)) are mapped onto two dominant supports in the SWFT, interference causes some of the power to be dissipated over the other frequencies in the form of spiky TFSs. Such ``leakage'' is present also in the case of amplitude/frequency modulation (note the several small peaks around the first support in Fig.\ \ref{fig:tfsupp}(e)). This effect is usually very small and does not lead to large errors in signal reconstruction (see Part II), but the existence of so many TFSs greatly complicates the interpretation of the resultant SWFT/SWT and, without some prior knowledge, might lead to the identification of far more components than are really present in the signal. This represents a significant drawback of the SWFT/SWT in comparison with the WFT/WT. Moreover, in contrast to the latter, the amplitude of the numerically calculated synchrosqueezed TFRs (\ref{numstfr}) depends on the width of the frequency bins (in a quite sophisticated way), as illustrated in Fig.\ \ref{fig:tfsupp}(e,f). This reflects the important fact that, for the SWFT/SWT, only the integral over the whole support makes sense, so that e.g.\ the component's amplitude cannot be estimated based on heights of the corresponding peaks, as can be done in the case of the WFT/WT (see below).

\begin{remark}
In the case of the SWFT/SWT, it has also been proposed \cite{Thakur:13,Auger:13} to reconstruct the components from the time-frequency area of constant width $\Delta$ around the corresponding ridge curve, i.e.\ to define the TFS simply as $\omega_{\pm}(t)=\omega_p(t)\pm\Delta/2$. Such an approach, however, is highly non-universal, as well as non-adaptive. Thus, the resolution properties of the SWFT/SWT are determined by the window/wavelet characteristics, so that $\Delta$ should depend on these characteristics as well, with no fixed choice being suitable in all cases. Furthermore, the spread of the SWFT/SWT around each peak is often asymmetric and varies in time, being determined by many factors such as the strength of noise, the amplitude/frequency modulation of the considered component and its interference with the other components (see Fig.\ \ref{fig:tfsupp}(e-h)). The time-frequency support as defined in this work implicitly takes account of such effects and is therefore quite universal and adaptive.
\end{remark}

It should be noted, that the TFS definition assumes that we know the ridge curve $\omega_p(t)$ associated with the component of interest, which in reality is of course not so. Thus, in real cases there exist many peaks in the TFR amplitude at each time and, moreover, their number can vary in time (e.g.\ due to noise). In such circumstances it is often difficult to tell which peaks correspond to which component, and which ones are not physically meaningful at all. The problem of finding the appropriate sequence of peaks is widely discussed in the context of ridge analysis \cite{Delprat:92,Carmona:97,Carmona:99,Lilly:10}. However, it is really a separate topic, which we consider in detail in \cite{Iatsenko:ridge}.

\subsection{Estimation of the component's parameters}\label{sec:tfsrecB}

Having found the ridge curve $\omega_p(t)$ and TFS $[\omega_{-}(t),\omega_{+}(t)]$ corresponding to the chosen AM/FM component $A(t)\cos\phi(t)$, the question which immediately arises is how to reconstruct its instantaneous characteristics, namely the amplitude $A(t)$, phase $\phi(t)$ and frequency $\nu(t)$. These parameters can be reconstructed by two methods: direct and ridge. Both these approaches are considered below, while their relative performance and the corresponding errors are studied in detail in Part II. Note that, as previously mentioned (see Sec.\ \ref{sec:analyticsignal}) the TFR methods can, at best, reconstruct the analytic approximation (\ref{ap}) to the component's parameters. Therefore, in certain rare cases when the analytic estimates are inaccurate (see \ref{app:asigerr}), all TFR-based estimates will be inaccurate to the same or higher extent.

\subsubsection*{Direct reconstruction}

Given its TFS, the parameters of the component can be reconstructed using the inversion formulas (\ref{iwft}), (\ref{iwt}), (\ref{istfr}), where one should restrict the integration over $\omega$ to only the corresponding time-frequency region $[\omega_-(t),\omega_+(t)]$. Such a method will be referred to as \emph{direct reconstruction}. However, in this way one can estimate only the amplitude and phase of the component (through reconstructing its analytic signal), but not its frequency. Nevertheless, the expressions for the latter, as well as for any order time-derivatives of amplitude and phase, can be derived in a similar way to (\ref{iwft}) and (\ref{iwt}), as discussed in \ref{app:derdirrec}. The direct estimates of the component's parameters are then given as
\begin{equation}\label{dirrec}
\begin{aligned}
&\mbox{\textbf{direct[WFT]:}}\\
&A(t)e^{i\phi(t)}=C_g^{-1}\int_{\omega_-(t)}^{\omega_+(t)} G_s(\omega,t)d\omega,\\
&\nu(t)=\operatorname{Re}\Bigg[\frac{\int_{\omega_-(t)}^{\omega_+(t)} \omega G_s(\omega,t)d\omega}
{\int_{\omega_-(t)}^{\omega_+(t)} G_s(\omega,t)d\omega}-\overline{\omega}_g\Bigg],\;\;
\overline{\omega}_g\equiv\frac{C_g^{-1}}{2}\int \xi\hat{g}(\xi)d\xi,\\
&(\overline{\omega}_g=0\mbox{ for symmetric }\hat{g}(\xi),\mbox{ such as Gaussian window}),\\\midrule
&\mbox{\textbf{direct[WT]:}}\\
&A(t)e^{i\phi(t)}=C_\psi^{-1}\int_{\omega_-(t)}^{\omega_+(t)} W_s(\omega,t)\frac{d\omega}{\omega},\\
&\nu(t)=\operatorname{Re}\Bigg[\frac{D_\psi^{-1}\int_{\omega_-(t)}^{\omega_+(t)} \omega W_s(\omega,t)\frac{d\omega}{\omega}}
{C_\psi^{-1}\int_{\omega_-(t)}^{\omega_+(t)} W_s(\omega,t)\frac{d\omega}{\omega}}\Bigg],\;\;
D_\psi\equiv\frac{\omega_\psi}{2}\int_0^\infty \frac{1}{\xi}\hat{\psi}^*(\xi)\frac{d\xi}{\xi},\\\midrule
&\mbox{\textbf{direct[SWFT]:}}\\
&A(t)e^{i\phi(t)}=\int_{\omega_-(t)}^{\omega_+(t)} V_s(\omega,t)d\omega,\;\;
\nu(t)={\rm Re}\Bigg[\frac{\int_{\omega_-(t)}^{\omega_+(t)} \omega V_s(\omega,t)d\omega}{\int_{\omega_-(t)}^{\omega_+(t)} V_s(\omega,t)d\omega}\Bigg],\\\midrule
&\mbox{\textbf{direct[SWT]:}}\\
&A(t)e^{i\phi(t)}=\int_{\omega_-(t)}^{\omega_+(t)} T_s(\omega,t)d\omega,\;\;
\nu(t)={\rm Re}\Bigg[\frac{\int_{\omega_-(t)}^{\omega_+(t)} \omega T_s(\omega,t)d\omega}{\int_{\omega_-(t)}^{\omega_+(t)} T_s(\omega,t)d\omega}\Bigg],\\
\end{aligned}
\end{equation}
where the formula for reconstructing $\nu(t)$ from the SWFT/SWT is introduced simply by analogy to the WFT/WT case, thus being derived rather empirically. In practice, the frequency axis is partitioned onto bins, so the formulas (\ref{dirrec}) should be also discretized, and the correct way of doing this is discussed in Sec.\ \ref{sec:pract2B}. Note also, that the instantaneous frequency cannot usually be estimated reliably through the numerical differentiation of the reconstructed phase, because any noise in the latter will be greatly amplified by such a procedure.

Except for the SWFT/SWT-based $\nu(t)$, all direct estimates (\ref{dirrec}) by definition give exact values (up to the error of the analytic approximation (\ref{ap})) in the ``ideal'' case when TFS contains all the energy of the component and no other contributions. The SWFT/SWT-based instantaneous frequency estimate, though not rigorously derived, seems also to be exact (or at least very close to exact) in this case, as will be seen in Part II.

\begin{remark}
The direct WT-based estimate of instantaneous frequency (\ref{dirrec}) is inapplicable for wavelets characterized by infinite $D_\psi$ (\ref{dirrec}). For the latter to be finite, one needs $|\hat{\psi}(\xi)|$ to decay faster than $\xi$ when $\xi\rightarrow0$, a condition that is not satisfied e.g.\ for the Morlet wavelet (\ref{mw}). In such circumstances one can use some kind of hybrid reconstruction
\begin{equation}\label{hybridrec}
\begin{aligned}
\mbox{\textbf{hybrid[WFT]:}}\;&
\nu(t)={\rm Re}\Bigg[\frac{\int_{\omega_-(t)}^{\omega_+(t)} \nu_G(\omega,t)G_s(\omega,t)d\omega}
{\int_{\omega_-(t)}^{\omega_+(t)} G_s(\omega,t)d\omega}\Bigg],\\
\mbox{\textbf{ hybrid[WT]:}}\;&
\nu(t)={\rm Re}\Bigg[\frac{\int_{\omega_-(t)}^{\omega_+(t)} \nu_W(\omega,t)W_s(\omega,t)\frac{d\omega}{\omega}}
{\int_{\omega_-(t)}^{\omega_+(t)} W_s(\omega,t)\frac{d\omega}{\omega}}\Bigg].
\end{aligned}
\end{equation}
Thus, when there is one-to-one correspondence between the component's TFS in the WFT/WT and that in the SWFT/SWT (e.g.\ for an ``ideal'' case), it follows from (\ref{stfr}) that reconstruction (\ref{hybridrec}) will give exactly the same values of $\nu(t)$ as direct reconstruction from the SWFT/SWT (\ref{dirrec}). Even when this is not the case, e.g.\ when there is noise or interference, the difference between the two estimates is usually negligible (see Part II). However, when both direct and hybrid reconstructions are possible, direct is obviously the more accurate (as it was rigorously derived, see \ref{app:derdirrec}), although the difference between the two is almost unnoticeable in the majority of cases. Note that, in the case of the WFT, the hybrid reconstruction is rarely needed, as direct frequency estimation (\ref{dirrec}) is possible for most of the window functions (being inapplicable only for a very exotic ones with $|\int \xi\hat{g}(\xi)d\xi|=\infty$).
\end{remark}

\subsubsection*{Ridge reconstruction}

The other widely used possibility is to reconstruct the component's parameters using TFR values at the ridge points $\omega_p(t)$ \cite{Delprat:92,Mallat:08,Lilly:10}, which will be referred to as the \emph{ridge reconstruction}. Thus, for the single tone signal $s(t)=A\cos(\nu t+\varphi)$ it follows from (\ref{Nwft}), (\ref{Nwt}) that the TFR amplitude at each time will be peaked at the tone frequency, and that the tone amplitude and phase can be perfectly reconstructed from the TFR value at any frequency as $Ae^{i(\nu t+\varphi)}=2G_s(\omega,t)/\hat{g}(\omega-\nu)$ (WFT) or $Ae^{i(\nu t+\varphi)}=2W_s(\omega,t)/\hat{\psi}^*(\omega_\psi\nu/\omega)$ (WT). Generalizing such an approach to the case of any AM/FM component, one obtains the ridge reconstruction formulas:
\begin{equation}\label{ridgerec}
\begin{aligned}
&\mbox{\textbf{ridge[WFT]:}}\\
&\nu(t)=\omega_p(t)+\delta\nu_d(t),\;\;
A(t)e^{i\phi(t)}=\frac{2G_s(\omega_p(t),t)}{\hat{g}(\omega_p(t)-\nu(t))},\\
\midrule
&\mbox{\textbf{ridge[WT]:}}\\
&\nu(t)=\omega_p(t)e^{\delta\log\nu_d(t)},\;\;
A(t)e^{i\phi(t)}=\frac{2W_s(\omega_p(t),t)}{\hat{\psi}^*(\omega_\psi\nu(t)/\omega_p(t))},\\
\midrule
&\mbox{\textbf{ridge[SWFT]:}}\\
&\nu(t)=\omega_p(t)+\delta\nu_d(t),\;\;
\phi(t)=\operatorname{arg}{\big[}V_s(\omega_p(t),t){\big]},\\
\midrule
&\mbox{\textbf{ridge[SWT]:}}\\
&\nu(t)=\omega_p(t)e^{\delta\log\nu_d(t)},\;\;
\phi(t)=\operatorname{arg}{\big[}T_s(\omega_p(t),t){\big]},\\
\end{aligned}
\end{equation}
where $\delta\nu_d(t)$ and $\delta\log\nu_d(t)$ are the corrections for discretization of the frequency scale (see below). The expressions for the amplitude estimates for the SWFT and SWT are not given because there is no possibility of reconstructing amplitude from the ridges of synchrosqueezed TFRs. This is because single points in the latter do not reflect the component's amplitude and will generally depend on the widths of frequency bins, as illustrated in Fig.\ \ref{fig:tfsupp} (see also \cite{Iatsenko:ridge,Montejo:12}). Note that, for the WT, the expressions (\ref{ridgerec}) should be modified if any other normalization than (\ref{wt}) is used (see Remark \ref{rem:wtrem}).

\begin{remark}\label{rem:phaseridge}
For the WFT/WT, instead of the amplitude ridge points $\omega_p(t)$, corresponding to maxima in the TFR amplitude, one can use phase ridge points, defined as the points where the frequency is equal to the TFR phase velocity $\omega_p(t)=\nu_{G,W}(\omega_p(t),t)$ and $[\partial_\omega\nu_{G,W}(\omega,t)]_{\omega=\omega_p(t)}<1$. Nevertheless, for the WFT with a Gaussian window it can be shown that both ridges coincide, representing one more advantage of this window function. In general, however, two ridge points will be different (even for the WT with lognormal/Morlet wavelet), but the difference between related estimates will be very small compared to the inherent errors of ridge reconstruction \cite{Lilly:10}. Therefore, the amplitude and phase ridge points can both be used, but the amplitude ridges seem to be the more convenient in practice (e.g. are easier to search for). Note also that phase ridges cannot be defined for the SWFT/SWT, because the corresponding phase velocities $\nu_{V,T}(\omega,t)$ (defined in the same way as those for the WFT/WT (\ref{iftfr})) are by construction very much the same as the frequencies for which they are calculated: $\nu_{V,T}(\omega,t)\approx\omega$.
\end{remark}

The corrections $\delta\nu_d(t)$ in (\ref{ridgerec}) arise because, in practice, the TFRs are calculated at the discrete frequencies $\omega_k=(k-k_0)\Delta\omega$ for the (S)WFT and $\omega_k/2\pi=2^{\frac{k-k_0}{n_v}}$ for the (S)WT. Therefore, the peak positions $\omega_p(t)$ are determined only up to the half-width of a frequency bin, leading to inaccuracies in all estimates. For the WFT/WT, these discretization errors can greatly be reduced by using quadratic interpolation to better locate the peak, so that in (\ref{ridgerec}) one sets
\begin{equation}\label{deltanu}
\begin{aligned}
\mbox{\textbf{WFT:}}\;\;&\delta\nu_d(t)=\frac{\Delta\omega}{2}\frac{a_3-a_1}{2a_2-a_1-a_3},\\
&a_{\{1,2,3\}}\equiv |G_s(\omega_{\{k_p(t)-1,k_p(t),k_p(t)+1\}},t)|,\\
\mbox{\textbf{ WT:}}\;\;&\delta\log\nu_d(t)=\frac{n_v^{-1}\log 2}{2}\frac{a_3-a_1}{2a_2-a_1-a_3},\\
&a_{\{1,2,3\}}\equiv |W_s(\omega_{\{k_p(t)-1,k_p(t),k_p(t)+1\}},t)|,\\
\end{aligned}
\end{equation}
where $k_p(t)$ denotes the discrete index of the peak at each time: $\omega_p(t)\equiv\omega_{k_p(t)}$. The corrections (\ref{deltanu}) then propagate to the amplitude and phase estimates (\ref{ridgerec}), making them more accurate as well. The interpolation cannot, however, be performed for the SWFT/SWT due to their non-smoothness, so in this case one uses the ``uncorrected'' estimates $\nu(t)=\omega_{k_p(t)}$ (susceptible to the discretization effects):
\begin{equation}\label{ssdeltanu}
\mbox{\textbf{SWFT:}}\;\delta\nu_d(t)=0,\;\;\;\;\;\;\mbox{\textbf{SWT:}}\;\delta\log\nu_d(t)=0.
\end{equation}
The discretization errors of ridge reconstruction for each TFR will be discussed in detail in Sec.\ \ref{sec:pract2C} below.

In contrast to direct reconstruction (\ref{dirrec}), ridge estimates (\ref{ridgerec}) are not exact, having errors proportional to the strengths of the amplitude and frequency modulations of the component \cite{Delprat:92,Lilly:10}. This is true even in the ``ideal'' case when the TFS contains the full energy of the component and no other influences. However, as will be seen in Part II, ridge reconstruction is much more robust in the face of noise and interference between the components than direct methods, so the choice between the two depends on the situation. The adaptive criteria for selecting the optimal estimates will be devised in Part II.

\begin{remark}
For the WFT/WT one can alternatively estimate the instantaneous frequency in (\ref{ridgerec}) through the TFR phase velocities at the peaks \cite{Lilly:10}: $\nu(t)=\nu_{G,W}(\omega_p(t),t)$. In fact, for the Gaussian window (only) one can show that $\nu_{G,W}(\omega_p(t),t)=\omega_p(t)$. However, for more exotic window/wavelet functions, such as the exponential window (see \ref{app:winwav}), one might encounter the situation where, for some peaks, $\nu_{G,W}(\omega_p(t),t)$ is completely different from $\omega_p(t)$. Hence, the original ridge reconstruction (\ref{ridgerec}) in conjunction with quadratic interpolation (\ref{deltanu}) appears to be a more universal approach than that based on the WFT/WT phase velocities.
\end{remark}

\section{Practical issues}\label{sec:pract}

In theory, one has infinite time and frequency scales, and both of these variables are continuous. In practice, however, everything is finite and discrete, which has specific consequences in terms of the resultant TFRs. In this section, the issues that arise while dealing with real signals are reviewed and studied.

\subsection{Signal preprocessing}\label{sec:pract1}

\begin{figure*}[t]
\includegraphics[width=1.0\linewidth]{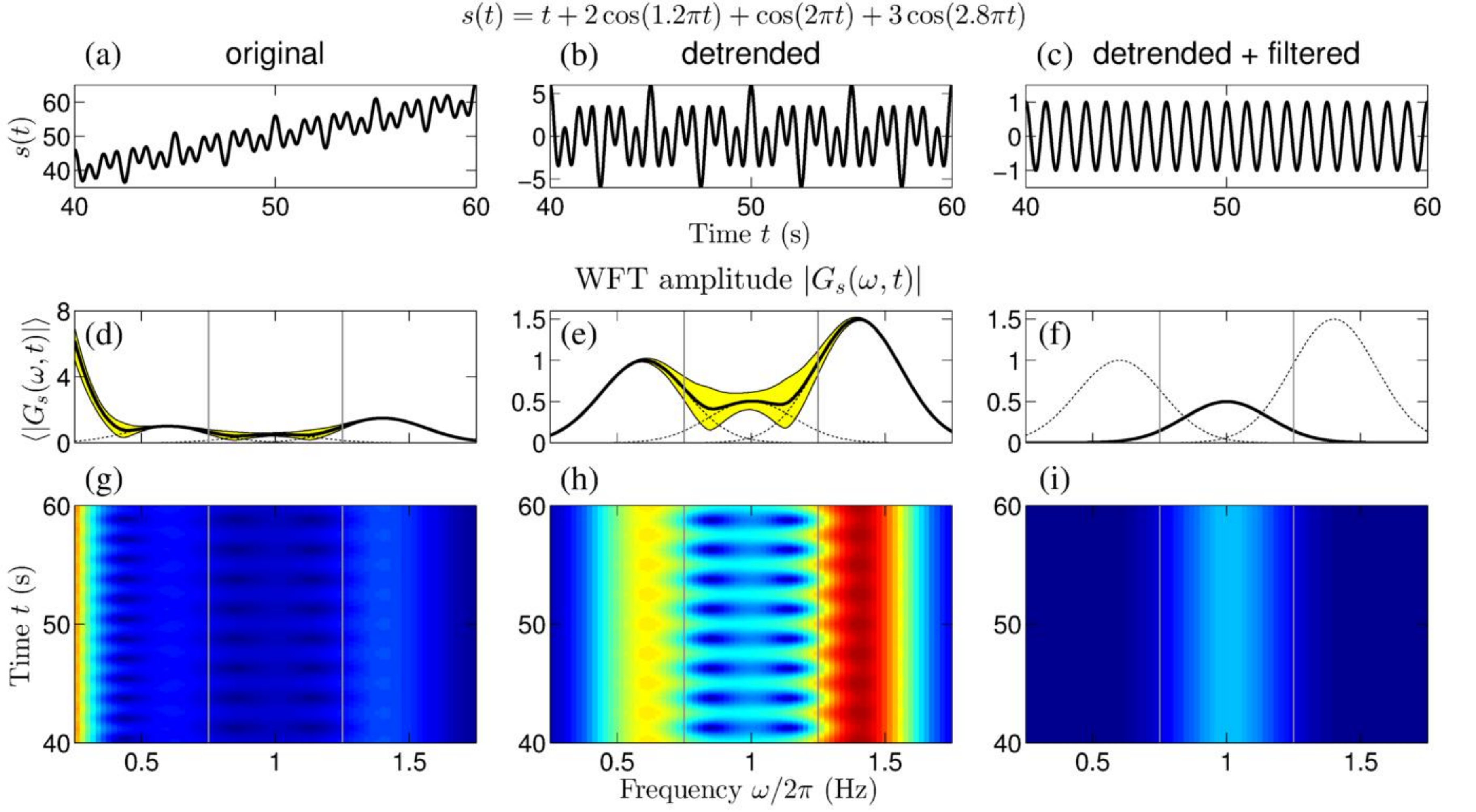}
\caption{Preprocessing and its effects on the resultant WFT (the case of the WT is qualitatively similar). (a): The original signal, consisting of the trend and three tones (shown at the top). (b): The same signal after detrending. (c): The same signal after detrending plus filtering in the band $[0.75,1.25]$ Hz. (d,e,f): The corresponding mean WFT amplitudes, with respective $\pm\sqrt{2}\std()$ being indicated by yellow regions; dotted lines show the mean amplitudes for the WFTs of each tone (out of 3) separately. (g,h,i) The corresponding WFT amplitudes in the time-frequency plane, with the respective color ranges being the same as $y$-axis limits in (d,e,f). It is assumed that one is interested in calculating the TFR for the frequency range [0.75,1.25] Hz, which is indicated by gray vertical lines in panels (d-i), but for completeness the resultant transforms are shown within wider ranges. The signal was sampled at 100 Hz for 100 s.}
\label{fig:preprocess}
\end{figure*}

To obtain a reliable TFR, an initial preprocessing of the signal should be performed. It consists of eliminating trends, followed by bandpass filtering in the frequency band of interest (for which the TFR is to be calculated). These two steps are considered below, with their effects being illustrated in Fig.\ \ref{fig:preprocess}.

\subsubsection{Removing the trends}\label{sec:pract1A}

If the signal contains a trend-like term (i.e.\ a term of the form $Kt$ or more generally $Kt^\alpha$), it can seriously corrupt the resultant TFR and complicate its interpretation. This is because trends make non-negligible contributions to the signal's spectral power in a wide frequency band. Thus, the FT of a trend existing for time $T$ will be $\int_0^T Kte^{-i\xi t}dt=\frac{KT}{\xi}e^{-i\xi T}+\frac{K}{\xi^2}(e^{-i\xi T}-1)$, which is proportional to the overall time of the trend and decays slowly with $\xi$. In general, any order contributions $Kt^n$, $n\in\mathds{N}$ will have FT $\sim iT^n e^{-i\xi T}/\xi$ in the first order over $\xi^{-1}$. Therefore, it is clear from the frequency domain form of (\ref{wft}) and (\ref{wt}) that trends might seriously affect the representation of the other components in the TFR, as illustrated in Fig.\ \ref{fig:preprocess}(d,g).

To avoid this, one should remove trends before doing any time-frequency analysis. Usually the best way to do so is to subtract a simple linear fit of the data, which will eliminate the $\sim\xi^{-1}$ spectral contribution of any term $\sim t^n,n\geq1$, changing it to $\sim\xi^{-2}$ for $n>1$. However, in general it seems even better to subtract a third order polynomial fit, which will in addition fully remove the trends $\sim t^2,t^3$ and reduce the spectral power of the higher order terms to $\sim \xi^{-4}$; it will also eliminate to a large extent the step-increases in mean value that are sometimes present in real data. This approach is very simple and introduces minimal undesirable distortions to the signal's spectrum, as can be seen by comparing (d,g) and (e,h) in Fig.\ \ref{fig:preprocess}. Third order is suggested because such a polynomial has at most 3 zero crossings and thus cannot model more than one oscillation during the whole time-series; and since oscillations having less than one cycle over the signal time length cannot in principle be reliably studied in the TFR (see Sec.\ \ref{sec:pract4} below), one does not lose anything by filtering them out.

\begin{remark}
Sometimes detrending is performed by subtracting a moving average from the signal, which additionally filters out the low-frequency spectral content. However, due both to the properties of such a filter and its associated boundary effects, this procedure usually introduces more unwanted spectral distortions than subtraction of a simple polynomial fit followed by bandpass filtering (see below).
\end{remark}

\subsubsection{Bandpass filtering}\label{sec:pract1B}

Using time-frequency analysis for a given signal, one is usually interested only in a particular frequency range (e.g.\ containing a chosen AM/FM component). At the same time, due to the peak broadening resulting from the TFR's finite frequency resolution, the frequency content of the signal {\it outside} the given range can significantly influence the TFR {\it inside} that range. This effect can be especially prominent when the spread of the window/wavelet function in the frequency domain is large. Therefore, one should always filter the signal in the frequency range of interest $[\omega_{\min},\omega_{\max}]$ before or during the application of the TFR, i.e.\ set $\hat{s}(\xi<\omega_{\min})=\hat{s}(\xi>\omega_{\max})=0$ in (\ref{wft}) and (\ref{wt}) (or use some other filter to remove spectral content lying outside the considered frequency range).

This issue is illustrated in Fig.\ \ref{fig:preprocess}(e,h), where the WFT of the tone in a frequency range of interest (indicated by gray vertical lines) is seriously corrupted by the other two tones that are near in frequency. However, as seen from Fig.\ \ref{fig:preprocess}(f,i), filtration of the signal only within the band considered solves this problem, allowing for an accurate representation of the corresponding component. Note that, performed alone, bandpass filtering does not fully remove the influence of a trend in the frequency range under consideration. That is why trends should be removed first, as described in the previous subsection.

\begin{remark}
It should be noted, that even if the instantaneous frequency of the AM/FM component lies within the considered frequency range $[\omega_{\min},\omega_{\max}]$, some of the related tones responsible for its amplitude/frequency modulation (see Eq.\ \ref{amfmft} and its discussion) might lie outside this range. In this case they will be filtered out, which will spoil the representation of the corresponding component to some extent. However, the probability that the related spectral content lies outside the considered band is not higher than the probability that there exists some unrelated components which, if unfiltered, might affect the components of interest to the same or a greater extent. Additionally, real signals are usually corrupted by noise, which can be viewed as many undesirable tones. Therefore, despite the possible drawbacks, it is generally preferable to bandpass filter than not to do so.
\end{remark}

\subsection{Frequency discretization}\label{sec:pract2}

In theory one has a continuous frequency variable $\omega$, for which TFRs are calculated. In practice, however, the frequency axis is discretized, being partitioned into bins centered at the chosen discrete values $\omega=\omega_k$, and the (S)WFT/(S)WT is calculated only at these frequencies. As discussed previously (see Sec.\ \ref{sec:tfr}), the (S)WFT and (S)WT have linear and logarithmic frequency resolutions, so that the discretization should also be performed on linear and logarithmic scales, respectively. It is convenient to take $\omega_k=(k-k_0)\Delta\omega$ (frequency bins $[(k-k_0-1/2)\Delta\omega,(k-k_0+1/2)\Delta\omega]$) for the (S)WFT and $\omega_k/2\pi=2^{(k-k_0)/n_v}$ (frequency bins $2\pi[2^{(k-k_0-1/2)/n_v},2^{(k-k_0+1/2)/n_v}]$) for the (S)WT. We now discuss the appropriate choice of $\Delta\omega$ and $n_v$, which determine the widths of the frequency bins, and how frequency discretization affects the resultant TFR and the estimation of the components' parameters from it.

\subsubsection{Choice of discretization parameters $\Delta\omega$ and $n_v$}\label{sec:pract2A}

\begin{figure*}[t]
\includegraphics[width=1.0\linewidth]{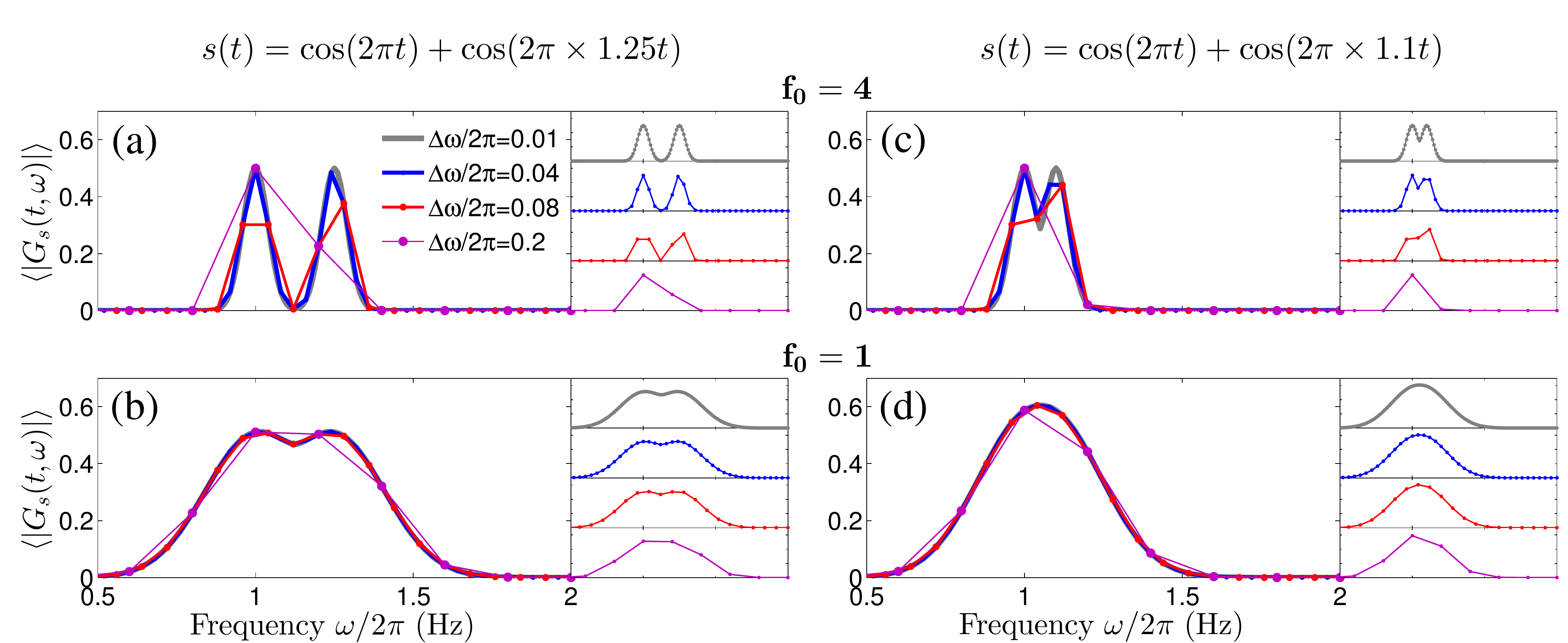}\\
\caption{Time-averaged WFT amplitudes for the signals (a,b) $s(t)=\cos(2\pi t)+\cos(2\pi\cdot 1.25t)$ and (c,d) $s(t)=\cos(2\pi t)+\cos(2\pi\cdot 1.1t)$, calculated using different resolution parameters $f_0$ and frequency bin widths $\Delta\omega$; the small insets on the right show the respective picture for each $\Delta\omega$ separately. The situation for the WT is qualitatively the same.}
\label{fig:fbins}
\end{figure*}

Theoretically the frequency resolution of the TFR is determined only by the window/wavelet parameters (see Sec.\ \ref{sec:tfres}). However, due to frequency discretization there also appears numerical frequency resolution, determined by the widths of the frequency bins, i.e.\ the choice of $\Delta\omega$ for the (S)WFT or the number-of-voices $n_v$ for the (S)WT. It imposes an upper bound on the effective frequency resolution of the transform, which is equal to the minimum between theoretical and numerical resolutions. Thus, if $\Delta\omega$ or $1/n_v$ is chosen to be too large, a few peaks in the TFR amplitude might be merged into one frequency bin, leading to an inability to distinguish between them and lowering the effective frequency resolution, as illustrated in Fig.\ \ref{fig:fbins}. On the other hand, if $\Delta\omega$ or $1/n_v$ is selected to be too small, it cannot improve frequency resolution beyond the theoretical maximum, but will increase computational cost due to requiring the calculation of the TFR at more frequencies. Therefore, $\Delta\omega$ and $n_v$ should be selected as such to retain the original, theoretical frequency resolution, predicted from the chosen window/wavelet parameters. This maximizes the effective resolution and occurs when the numerical frequency resolution is the same as or better than the theoretical one. Obviously, the optimal values of $\Delta\omega$ and $n_v$ depend on the chosen form of window/wavelet and its parameters (see Fig.\ \ref{fig:fbins}). Thus, while for $f_0=1$ a frequency step of $\Delta\omega/2\pi=0.08$ is sufficient to distinguish between two nearby but theoretically resolved components (Fig.\ \ref{fig:fbins}(b)), for $f_0=4$ it is already insufficient (Fig.\ \ref{fig:fbins}(c)).

Evidently, to avoid reducing the original frequency resolution by discretization, for the WFT the width of a frequency bin $\Delta\omega$ should be smaller than the theoretical $\Delta\nu_{\min}$ (see Sec.\ \ref{sec:tfres}); while for the WT one has the same criterion but on a logarithmic frequency scale. In other words, one should break the minimum resolvable frequency difference (WFT) or ratio (WT) into a large enough number of segments. However, one is now interested in preserving the TFR's appearance more generally, i.e.\ not only for a reliably resolved components, but in relation to {\it any} components that can be distinguished in the TFR (even where they are substantially corrupted by interference). Thus, to determine the optimal frequency bin widths, the $\epsilon_r$ in $\Delta\nu_{\min}$ (\ref{tfreswft}) and (\ref{tfreswt}) should be set quite high, and we consider $\epsilon_r=0.5$ as an extreme case. Based on these considerations, one can select an appropriate $\Delta\omega$ or $n_v$ as
\begin{equation}\label{fbins}
\begin{aligned}
\mbox{\textbf{(S)WFT:}}\;\Delta\omega&=\frac{\xi_2(0.5)-\xi_1(0.5)}{N_b}\\
&\approx 1.35/f_0N_b\mbox{ for Gaussian window (\ref{gw})},\\
\mbox{\textbf{ (S)WT:}}\;\;\;n_v&=\frac{N_b\log 2}{\log\xi_2(0.5)-\log\xi_1(0.5)}\\
&\approx 3.23f_0N_b\mbox{ for lognormal wavelet (\ref{lw})},\\
\end{aligned}
\end{equation}
where $N_b>1$ denotes the chosen number of bins into which we divide the 50\% window/wavelet FT support (we use $N_b=10$ in calculations, and $N_b=20$ for figures), and $\xi_{1,2}(\epsilon)$ are as defined in (\ref{eswft}) and (\ref{eswt}) for the WFT and WT, respectively. Since synchrosqueezing does not seem to change frequency resolution (see Sec.\ \ref{sec:dwftwt} and Part II), the same frequency discretization (\ref{fbins}) appears to be appropriate for the SWFT/SWT as well (unless one wants to estimate the instantaneous frequency from it, see below).

The criteria (\ref{fbins}) are based on the results of Sec.\ \ref{sec:tfres}, where we studied the resolution of two equal-amplitude tones. However, it is clear that the best separation occurs for tones of the same amplitude; while in the case of different amplitudes (e.g.\ $s(t)=\cos\nu t+0.1\cos(\nu+\Delta\nu) t$), the difference $\Delta\nu$ between the frequencies of the tones should be larger in order for them to be resolved (see also Part II). Consequently, the optimal $\Delta\omega$ and $n_v$ (\ref{fbins}) are sufficient to resolve \emph{any} components if they are separated theoretically, as desired.

\begin{remark}
It should be noted, that even with the choice (\ref{fbins}) it is still possible that at certain times close peaks will appear in the WFT/WT amplitude that are separated theoretically, but not practically, being merged due to insufficient numerical resolution. However, as mentioned previously, if the bin widths (\ref{fbins}) are not small enough to resolve two tones, most of the time they will behave as a single component in the TFR, and it will be impossible to investigate them separately. Hence, attributing the corresponding peaks to one or few components will not make much difference.
\end{remark}

\subsubsection{Discretization errors of direct parameter estimation}\label{sec:pract2B}

While reconstructing signal parameters by direct methods (\ref{dirrec}), in theory one integrates over the continuous TFS at each time; but in practice we discretize integrals to the finite sum over the frequency bins. This discretization obviously introduces some numerical errors into the direct schemes. Nevertheless, it appears that these errors are very small (except in the case of SWFT/SWT frequency reconstruction, see below), at least if the frequency discretization is not too coarse and one performs the numerical integration correctly, i.e.\ using the mid-point rule. For example, even in the case of $\Delta\omega/2\pi=0.08$ at $f_0=4$, which is obviously insufficient for an accurate representation (see Fig.\ \ref{fig:fbins}(c)), one can reconstruct the single tone $s(t)=\cos\nu t$ from its WFT with relative accuracy $<2\%$ while, for $\Delta\omega$ determined by (\ref{fbins}) ($\approx 2\pi\cdot 0.005$ for Gaussian window with $f_0=4$, and using $N_b=10$), the relative inaccuracy is $<10^{-10}\%$.

However, the midpoint rule generally needs equidistant points so, in the case of the equilogspaced frequencies used for (S)WT, one should also discretize the integrals on a logarithmic scale. The correct discretization rules are summarized below
\begin{equation}\label{drules}
\begin{aligned}
&\mbox{\textbf{Linear frequency scale (SWFT/WFT)}}:\\
&\int_{\omega_-}^{\omega_+}[...(\omega)]d\omega\rightarrow \sum_{k:\omega_k\in[\omega_-,\omega_+]}[...(\omega_k)]\frac{\omega_{k+1}-\omega_{k-1}}{2},\\
&\mbox{\textbf{Logarithmic frequency scale (SWT/WT)}}:\\
&\int_{\omega_-}^{\omega_+}[...(\omega)]d\omega\rightarrow
\left[\begin{array}{l}
\mbox{RIGHT}: \sum_{k:\omega_k\in[\omega_-,\omega_+]}[...(\omega_k)]\frac{\omega_k}{2}\log \frac{\omega_{k+1}}{\omega_{k-1}},\\
\mbox{NOT RIGHT}: \sum_{k:\omega_k\in[\omega_-,\omega_+]}[...(\omega_k)]\frac{\omega_{k+1}-\omega_{k-1}}{2},\\
\end{array}\right.\\
&\mbox{\textbf{Numerical SWFT/SWT (any frequency scale)}}:\\
&\int_{\omega_-}^{\omega_+}V_s(\omega,t)[...(\omega)]d\omega\rightarrow\sum_{k:\omega_k\in[\omega_-,\omega_+]}\widetilde{V}_s(\omega_k,t)[...(\omega_k)],\\
&\int_{\omega_-}^{\omega_+}T_s(\omega,t)[...(\omega)]d\omega\rightarrow\sum_{k:\omega_k\in[\omega_-,\omega_+]}\widetilde{T}_s(\omega_k,t)[...(\omega_k)],\\
\end{aligned}
\end{equation}
where for the usual frequency discretization one has $(\omega_{k+1}-\omega_{k-1})/2=\Delta\omega$ for (S)WFT and $\frac{1}{2}\log \frac{\omega_{k+1}}{\omega_{k-1}}=\frac{\log 2}{n_v}$ for (S)WT. In (\ref{drules}), $[...(\omega)]$ means any functional of frequency (e.g.\ for WT-based direct signal reconstruction (\ref{iwt}) one has $[...(\omega)]=\frac{W_s(\omega,t)}{\omega}$). The SWFT/SWT calculated in practice $\widetilde{V}_s,\widetilde{T}_s$ (\ref{numstfr}) represents the theoretical SWFT/SWT $V_s,T_s$ already integrated over the one frequency bin (using the correct rule), so the discretization procedures for them are very simple. Note that the rules outlined here apply generally, and not only for direct reconstruction methods, e.g.\ discretization while performing numerical synchrosqueezing (\ref{numstfr}) is also done according to (\ref{drules}).

The relative discretization error $\epsilon_d$ in the direct reconstruction (\ref{dirrec}) of the signal $s^{a}(t)$ (and so also in its amplitude and phase) can be estimated through the error of discretizing the integrals determining $C_{g,\psi}$ (\ref{iwft}), (\ref{iwt}). Using mid-point rule error estimates, one has
\begin{equation}\label{derr}
\begin{aligned}
\mbox{\textbf{WFT:}}\;
\epsilon_d
\lesssim& \max_\nu \frac{C_g^{-1}}{2}{\Big|}\int \hat{g}
(\omega-\nu)d\omega-\sum_{n=-\infty}^\infty \hat{g}(n\Delta\omega-\nu)\Delta\omega{\Big|}\\
\leq& \max_\nu \frac{C_g^{-1}}{2}\frac{\Delta\omega^3}{24}\sum_{n=-\infty}^\infty \Big[\big|\hat{g}''(n\Delta\omega-\nu)\big|+O(\Delta\omega)\Big]\\
\approx& \Bigg[\frac{C_g^{-1}}{2}\int_{-\infty}^{\infty} |\hat{g}''(\omega)|d\omega\Bigg]\frac{\Delta\omega^2}{24}+O(\Delta\omega^3),\\
\mbox{\textbf{WT:}}
\;\epsilon_d
\lesssim& \max_\nu \frac{C_\psi^{-1}}{2}{\Big|}\int_{-\infty}^\infty \hat{\psi}(e^x\nu)dx-\sum_{n=-\infty}^\infty \hat{\psi}(2^{n/n_v}\nu)\frac{\log 2}{n_v}{\Big|}\\
\leq& \max_\nu \frac{C_\psi^{-1}}{2}\frac{(\frac{\log 2}{n_v})^3}{24}
\sum_{n=-\infty}^\infty \Big[\big|2^{n/n_v}\nu\hat{\psi}'(2^{n/n_v}\nu)\\
&+(2^{n/n_v}\nu)^2\hat{\psi}''(2^{n/n_v}\nu)\big|+O(n_v^{-1})\Big]\\
\approx& \Bigg[\frac{C_\psi^{-1}}{2}\int_{0}^{\infty} |\hat{\psi}'(\omega)+\omega\hat{\psi}''(\omega)|d\omega\Bigg]\frac{\log^2 2}{24n_v^2}+O(n_v^{-3}).\\
\end{aligned}
\end{equation}
The true error, however, is usually much smaller than the upper bound given by (\ref{derr}). It is hard to derive error bounds in the case of the SWFT/SWT, but it appears that same situation arises there as well. Note that the expressions (\ref{derr}) are derived on the assumption that the first and second derivatives of the window/wavelet FT are finite everywhere; if this is not so (e.g.\ $\hat{g}(\xi)\sim e^{-\alpha|\xi|}\Rightarrow\hat{g}''(0)=\infty$), then the error increases to $O(\Delta\omega)$, and discretization effects become considerable (which is a significant drawback of such windows/wavelets).

For direct reconstruction of the frequency (\ref{dirrec}), the situation is slightly different. For the WFT and WT the absolute error, apart from $\epsilon_d$ (\ref{derr}), will include the discretization error for the integral determining $\overline{\omega}_g$ (WFT) or $D_\psi$ (WFT) in (\ref{dirrec}), respectively. They can be derived in the same way as (\ref{derr}), and so the total error is $O(\Delta\omega^2)$ (WFT) and $O(n_v^{-2})$ (WT), which is quite small in theory and even smaller in practice. On the other hand, the SWFT/SWT-based direct frequency estimates are considerably affected by discretization effects. Thus, for a single tone signal $s(t)=\cos\nu t$ the SWFT/SWT has only one nonzero coefficient at each time, so the discretization error of the direct estimate of the tone's frequency (\ref{dirrec}) is $\leq\Delta\omega/2$ for the SWFT and $\leq(1-2^{-1/2n_v})\nu$ for the SWT. These seem to be the maximal error bounds: although the situation becomes more complicated when amplitude/frequency modulation is present, the error is expected to be smaller (and it is indeed so in practice) because the SWFT/SWT now contains few non-zero coefficients at each time and the numerical integration is no longer limited to a single point.

To conclude, with the use of frequency bin widths determined by (\ref{fbins}) and the correct integration rules (\ref{drules}), the discretization errors of direct parameter estimation (\ref{dirrec}) become very small and are in practice usually negligible, so one need not worry about them. The situation is different, however, for the SWFT/SWT-based frequency estimates, which suffer greatly from discretization effects. Note, that one can use more sophisticated integration rules, e.g.\ Gauss quadrature, but the simple and intuitive mid-point method (\ref{drules}) is already extremely accurate for most purposes.

\subsubsection{Discretization errors of ridge-based parameter estimation}\label{sec:pract2C}

Ridge reconstruction methods estimate all from the peaks in the TFR amplitude, the positions of which are not precisely determined due to frequency discretization. However, in the case of the WFT/WT, the related discretization errors are nearly eliminated by using the corrections $\delta\nu_d(t),\delta\log\nu_d(t)$ (\ref{deltanu}) in the formulas (\ref{ridgerec}). Without this modification, the ridge estimates can suffer greatly from discretization effects.

Thus, consider the paradigmatic case of a single tone signal $s(t)=A\cos\nu t\rightarrow G_s(\omega_k,t)=\frac{A}{2}\hat{g}(\omega_k-\nu)e^{i\nu t},\;W_s(\omega_k,t)=\frac{A}{2}\hat{\psi}^*(\omega_\psi\frac{\nu}{\omega_k})$, with $\omega_k=(k-k_0)\Delta\omega$ for the WFT and $\omega_k/2\pi=2^{(k-k_0)/n_v}$ for the WT. For this signal, the peak in the WFT/WT amplitude will occur at the discrete frequency $\omega_p(t)=\omega_{k_p(t)}$ that is closest to $\nu$, so that the maximum difference between the two is equal to the half-width of a frequency bin. Therefore, if using the conventional (unmodified) ridge reconstruction formulas, i.e.\ (\ref{ridgerec}) with $\delta\nu_d(t)=\delta\log\nu_d(t)=0$, the discretization-related inaccuracies can be shown to be
\begin{equation}\label{rerr}
\begin{aligned}
&\mbox{\textbf{WFT:}}\;\;|\Delta\nu_d(t)| \leq\frac{\Delta\omega}{2},\;\;\hat{g}(\xi)\equiv a_g(\xi)e^{i\varphi_g(\xi)}\\
&\Big|\frac{\Delta A_d(t)}{A}\Big|\leq \max\Big|1-\frac{a_g(\pm\Delta\nu_d)}{a_g(0)}\Big|
\leq\Big|\frac{a_g''(0)}{2a_g(0)}\frac{\Delta\omega^2}{4}+O(\Delta\omega^3)\Big|\\
&\big|\Delta\phi_d(t)\big| \leq\big|\varphi_g(\pm\Delta\nu_d(t))-\varphi_g(0)\big|\\
&\hphantom{\big|\Delta\phi_d(t)\big|}\leq\Big|\varphi_g'(0)\frac{\Delta\omega}{2}+\frac{\varphi_g''(0)}{2}\frac{\Delta\omega^2}{4}+O(\Delta\omega^3)\Big|\\
\midrule
&\mbox{\textbf{WT:}}\;\;|\Delta\log\nu_d(t)|\leq\frac{n_v^{-1}\log 2}{2},\;\;\hat{\psi}(\xi)\equiv a_\psi(\xi)e^{i\varphi_\psi(\xi)}\\
&\Big|\frac{\Delta A_d(t)}{A}\Big| \leq\max\Big|1-\frac{a_\psi(e^{\pm\Delta\log\nu_d(t)}\omega_\psi)}{a_\psi(\omega_\psi)}\Big|\\
&\hphantom{\Big|\frac{\Delta A_d(t)}{A}\Big|} \leq\Big|\frac{\omega_\psi^2a_\psi''(\omega_\psi)}{2a_\psi(\omega_\psi)}\frac{n_v^{-2}\log^2 2}{4}+O(n_v^{-3})\Big|\\
&\big|\Delta\phi_d(t)\big| \leq\big|\varphi_g(e^{\pm\Delta\log\nu_d(t)}\omega_\psi)-\varphi_\psi(\omega_\psi)\big|\\
&\hphantom{|\Delta\phi_d(t)|} \leq\Big|\omega_\psi\varphi_\psi'(\omega_\psi)\frac{n_v^{-1}\log 2}{2}
+\frac{\omega_\psi^2\varphi_\psi''(\omega_\psi)}{2}\frac{n_v^{-2}\log^2 2}{4}+O(n_v^{-3})\Big|\\
\end{aligned}
\end{equation}
where we have taken into account that, by definition, $a_g'(0)=0$ and $a_\psi'(\omega_\psi)=0$. Note, that for the present case of the single tone signal the phase error $|\Delta\phi_d(t)|$ vanishes when $\varphi_{g,\psi}(\xi)=\operatorname{const}$ (e.g.\ for real $\hat{g}(\xi)$ or $\hat{\psi}(\xi)$), which is satisfied by most of the windows/wavelets in use. However, when amplitude/frequency modulation is present, then $\partial_\omega \operatorname{arg}[G(\omega,t)]$ will often not be equal to zero at the TFR amplitude peak, even for windows/wavelets with real FTs, so that in general one has $\Delta\phi_d(t)=O(\Delta\omega),\,O(n_v^{-1})$.

Nevertheless, by using the quadratic interpolation (\ref{deltanu}) in (\ref{ridgerec}), one can greatly reduce the discretization errors, as compared to (\ref{rerr}). Thus, one can show that for the general case of any AM/FM component the related errors become $\Delta\nu_d(t)=O(\Delta\omega^2)$, $\Delta A_d(t)/A=O(\Delta\omega^3)$, $\Delta\phi_d(t)=O(\Delta\omega^2)$ for the WFT and respectively ($\Delta\omega\rightarrow n_v^{-1}\log 2$) for the WT. These inaccuracies are of the same (frequency, phase) or higher (amplitude) order in frequency bin width as the inaccuracies of the direct estimates (see Sec.\ \ref{sec:pract2B}), which are negligible in practice. Therefore, selecting discretization parameters by (\ref{fbins}) and using correction terms (\ref{deltanu}) in (\ref{ridgerec}), one need not worry about the discretization errors of ridge-based parameter estimation from the WFT and WT. Note, that the above considerations apply only if $|\hat{g}'(0)|<\infty,\,|\hat{g}''(0)|<\infty$ and $|\hat{\psi}'(\omega_\psi)|<\infty,\,|\hat{\psi}''(\omega_\psi)|<\infty$; otherwise, e.g.\ if $a_g''(0)=\infty$ (as for $\hat{g}(\xi)\sim e^{-|\xi|}$), all inaccuracies raise to $O(\Delta\omega)$ and become considerable.

For the SWFT/SWT the situation is different, and ridge reconstruction from synchrosqueezed TFRs will generally suffer from discretization effects. Thus, as discussed in Sec.\ \ref{sec:tfsrecB}, there are no straightforward ways of reducing the related errors of ridge reconstruction. As a result, the frequency inaccuracies will be $|\Delta\nu_d(t)|\leq\Delta\omega/2$ for the SWFT and $|\Delta\nu_d(t)|\leq n_v^{-1}\log 2/2$ for the SWT, being the same as for the SWFT/SWT-based direct reconstruction (see Sec.\ \ref{sec:pract2B}). Hence, for accurate estimation of the instantaneous frequency from the SWFT/SWT, by any method, one needs to use very fine frequency binning even for large spreads of $\hat{g}(\xi)$ and $\hat{\psi}(\xi)$. This increases (sometimes enormously) the computational cost and thus represents a significant drawback of synchrosqueezing. Furthermore, in contrast to the corresponding direct estimates, the SWFT/SWT-based ridge estimates of phase also appear to contain discretization errors of order $O(\Delta\omega)$ and $O(n_v^{-1})$ (this claim is based on numerical evidence and is hard to derive rigorously). As mentioned previously, there is no possibility of amplitude reconstruction from SWFT/SWT ridges.

\subsection{Boundary effects and padding}\label{sec:pract3}

Theoretically, one integrates over an infinite time or frequency while computing the WFT (\ref{wft}) and WT (\ref{wt}), but in practice the signal has a finite time duration and sampling frequency. Consequently, the resultant TFR becomes ill-defined near the signal's time borders (when $t$ is close either to $0$ or to the overall time-length $T$). Irrespectively of how this problem is tackled, it often leads to distortions of the WFT/WT near both signal ends -- \emph{boundary effects}. Obviously, the same feature propagates to the SWFT/SWT.

\subsubsection{Padding and its schemes}\label{sec:pract3A}

\begin{figure*}[t]
\includegraphics[width=1.0\linewidth]{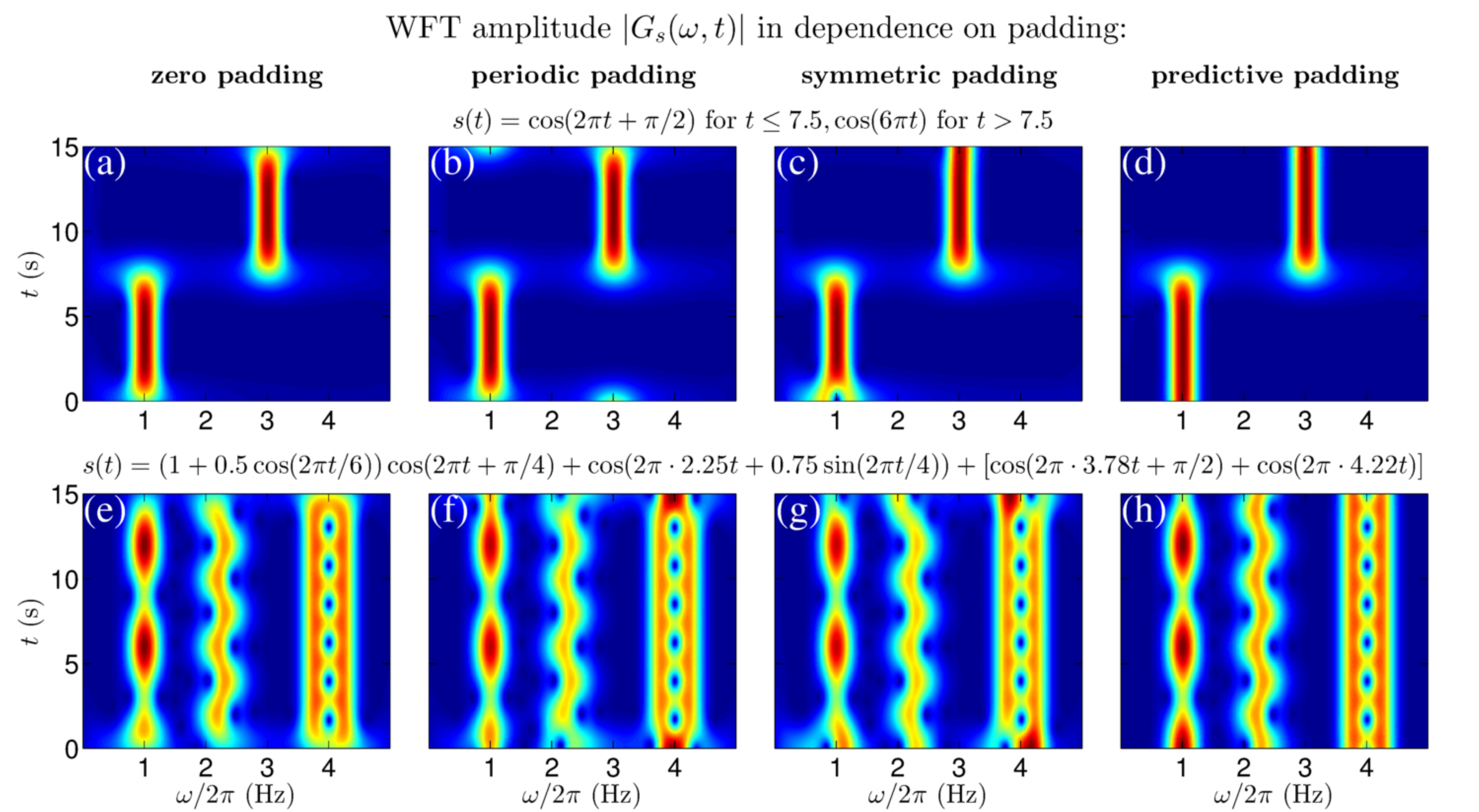}\\
\caption{Examples of the WFT amplitudes calculated using each of the four padding schemes discussed in the text (columns) for two different signals (rows): (a,e) zero padding; (b,f) periodic padding; (c,g) symmetric padding; (d,h) predictive padding. The effects of different padding strategies on the SWFT and (S)WT are qualitatively the same. The signals are indicated in the figure, and they were sampled at $f_s=100$ Hz for $T=15$ s.}
\label{fig:padding}
\end{figure*}

Although this is not really the case in practice (see below), suppose for now that the WFT and WT are calculated using the convolution in the time domain in (\ref{wft}) and (\ref{wt}), respectively. The signal's positive part $s^+(\tau)$ is thus multiplied by $g(\tau-t)e^{-i\omega(\tau-t)}$ (WFT) or $\psi^*(\omega(\tau-t)/\omega_\psi)$ (WT) and integrated over $\tau$. Therefore, one needs to devise a rule by which the integration can be performed when $\tau$ lies outside the signal's time limits, i.e.\ when $\tau<0$ or $\tau>T$. Generally, this can be done by continuing the signal beyond its original time interval $[0,T]$ in some way, e.g.\ setting $s^+(\tau\notin[0,T])=0$ in (\ref{wft}), (\ref{wt}). In practice, the signal is padded at both ends using a particular convention, then the convolution is calculated with this padded signal, and after this only the part lying within the original time limits is taken. An appropriate padding can be constructed by many different schemes, and the most common/useful ones are listed below, with their effect on the resultant TFR being illustrated in Fig.\ \ref{fig:padding}. Note, that padding should be performed \emph{after} the initial signal preprocessing, discussed in Sec.\ \ref{sec:pract1}, has been completed.

\textbf{Zero padding}, when one pads the signal with zeros at both ends, represents the simplest and most predictable form of padding. This scheme effectively sets $s(\tau<0)=s(\tau>T)=0$ in (\ref{wft}), (\ref{wt}), ``cutting'' the convolution when $\tau$ goes beyond the signal's time limits. Its effect is illustrated in Fig.\ \ref{fig:padding}(a,e): as can be seen, this type of padding does not introduce any new behavior (in contrast to a few other types, see below), always leading to a gradual decay of the TFR amplitude towards the time ends. This allows one to derive the expressions for boundary errors that arise in the case of zero padding, which will be done in the following.

\textbf{Periodic padding} is constructed by periodic continuation of the signal. It leads to quite unpredictable boundary effects, with the TFR near the time borders being affected by the signal's behavior at both ends. For example, as illustrated in Fig.\ \ref{fig:padding}(b), the tone occurring at the beginning (end) has its ``phantom'' tail near the time end (beginning) in the WFT. Consequently, this type of padding does not represent the preferred choice, unless one knows that the signal is fully periodic over the whole interval, which almost never occurs in practice.

\textbf{Symmetric padding}, whose effects are shown in Fig.\ \ref{fig:padding}(c,g), is performed by reflecting the signal along $t=0$ and $t=T$. It localizes boundary inaccuracies at each edge, thus solving the problem of the end-to-end influence that occurs for periodic padding. Indeed, comparing (c,g) and (b,f) of Fig.\ \ref{fig:padding}, it can be seen that the ``phantom tails'' near the boundaries have been removed. However, the effect of symmetric padding depends to a large extent on the phases of oscillations at $t=0,T$, which are quite unpredictable. Thus, if the initial/end phase is not equal to zero or $\pi$, symmetric padding introduces a phase jump and leads to a splitting of the single peak in the TFR amplitude into several peaks near the time borders. As an example, in Fig.\ \ref{fig:padding}(b) the first oscillation has a phase shift of $\pi/2$ at $t=0$ and thus is doubled, while the second one has zero shift at $t=T$ and thus is well represented.

\textbf{Predictive padding}, as implied by its name, continues the signal beyond its time limits by inferring/forecasting its past/future behavior. Thus, one tries to predict the values of the signal by assuming some model of the process generating the time-series and fitting data to this model to find its parameters. In the present case, one is mainly interested in continuing the signal in such a way as to best represent in the TFR its existing characteristics near the boundaries, and not in finding the unknown signal behavior beyond the available data, which is generally impossible. Because time-frequency analysis is devoted to studying the oscillatory properties of the data, it seems that the most relevant approach would be to forecast the signal based on its local spectral characteristics at both ends. A scheme for doing this is discussed in \ref{app:forecasting}, while the effect of such padding is shown in Fig.\ \ref{fig:padding}(d,h). As can be seen, for both signals in Fig.\ \ref{fig:padding} predictive padding almost completely eliminates the boundary distortions. However, it does not represent an ``ultimate cure'', and for a complicated signal some errors might remain; they depend on the signal's structure and thus are hard to estimate in general. Nevertheless, in terms of its ability to reduce boundary effects, predictive padding usually outperforms all other schemes, and it is therefore used by default in this work.

The remaining question is how many values to pad. At the beginning of this section it was assumed for simplicity that the WFT/WT are calculated using the time domain forms of (\ref{wft}), (\ref{wt}). Such an approach, however, is computationally very expensive, requiring O($N$) convolutions for each $t_{n=1,...,N}$ and thus being of O($N^2$) in cost. In practice, therefore, the WFT and WT are calculated using the frequency domain convolution in (\ref{wft}) and (\ref{wt}) which, utilizing the FFT algorithm, can be performed in O($N\log N$) computations (see Sec.\ \ref{sec:stepalg} below). However, the discrete FT of the signal $\hat{s}(\xi_n)$, which is used in numerical convolution, represents the periodic spectrum estimate, being an exact FT of the periodically continued signal: $\hat{s}(\xi_n)=\int s({\rm mod}_T t)e^{-i\xi_n t}dt$. Thus, no padding is in practice completely equivalent to periodic padding.

One therefore needs to add enough points at both signal ends to guarantee that, within the original signal time-length, the effects of implicit periodic continuation due to FFT-based convolution are small. This minimum number of padded values will obviously depend on the spread of the window/wavelet in time: the larger it is, the more values one should pad so that most of the window/wavelet is convoluted with padded values and not with the periodically-continued other end of the signal. Since the effective window/wavelet length in time (up to a predefined tolerance $\epsilon$) can be expressed through its $\epsilon$-supports $\tau_{1,2}(\epsilon)$ (\ref{eswft}), (\ref{eswt}), the minimum number of points $n_1^{(\min)}$ ($n_2^{(\min)}$) with which one should pad the signal for $t<0$ ($t>T$) can be determined as
\begin{equation}\label{minpad}
\begin{aligned}
\mbox{\textbf{(S)WFT:}}\;&n^{(\min)}_{1,2}(\epsilon)=f_s|\tau_{1,2}(\epsilon)|,\\
\mbox{\textbf{ (S)WT:}}\;&n^{(\min)}_{1,2}(\epsilon)=\frac{\omega_\psi}{\omega_{\min}}f_s|\tau_{1,2}(\epsilon)|,\\
\end{aligned}
\end{equation}
where $f_s=1/\Delta t$ is the sampling frequency of the signal, and $\omega_{\min}$ denotes the minimum frequency for which the TFR is calculated. Thus, for the WT the wavelet is rescaled at each frequency, so the number of points needed to assure the specified precision $\epsilon$ also scales; to guarantee the accuracy for all frequencies, the maximum number of points is taken, which corresponds to the lowest frequency.

Additionally, the FFT algorithm requires the total length of the signal in samples to be a power of two. Therefore, in practice one should simultaneously assure both the power-of-two points and criterion (\ref{minpad}), which can be done by padding the signal from originally $N$ to $N_p$ points, with $n_1$ values to the left ($t<0$) and $n_2$ to the right ($t>T$), given by
\begin{equation}\label{numpad}
\begin{gathered}
N_p={\rm NextPowerOfTwo}[N+n_1^{(min)}(\epsilon)+n_2^{(min)}(\epsilon)]=N+n_1+n_2,\\
n_{1,2}=\frac{n_{1,2}^{(\min)}(\epsilon)}{n_1^{(\min)}(\epsilon)+n_2^{(\min)}(\epsilon)}(N_p-N),\\
\end{gathered}
\end{equation}
where $n_{1,2}^{(\min)}(\epsilon)$ are as defined in (\ref{minpad}); by default we use $\epsilon=0.001$.

\subsubsection{Error estimates and the cone-of-influence}\label{sec:pract3B}

As discussed above, for periodic, symmetric and predictive padding the boundary effects very much depend on the signal's structure. For zero padding, however, they are quite universal and can thus be estimated. A simple and straightforward way to do so is to quantify the relative boundary inaccuracies $\epsilon_b(\omega,t)$ through the difference between the theoretical (\ref{Nwft}), (\ref{Nwt}) and ``practical'' WFT/WT of the single-tone signal $s(t)=\cos(\nu t+\varphi)$, calculated at the tone frequency $\nu$. Then in the case of zero padding one obtains
\begin{equation}\label{berrwft}
\begin{aligned}
\mbox{\textbf{WFT:}}\;\;\;\;&\\
\epsilon_b(\nu,t)
\equiv&\frac{\big|G_s(\nu,t)-\widetilde{G}_s(\nu,t)\big|}{\left|G_s(\nu,t)\right|}\\
=&\frac{\left|\int g(\tau-t)d\tau
-\int_0^T g(\tau-t)d\tau\right|}
{\left|\int g(\tau-t)d\tau\right|}\\
\leq& \frac{\big|\int_{-\infty}^{0}g(\tau-t)d\tau\big|+\big|\int_{T}^{\infty}g(\tau-t)d\tau\big|}{\big|\int g(\tau)d\tau\big|}\\
=&|P_g(-t)|+|P_g(t-T)|,\\
\end{aligned}
\end{equation}
\begin{equation}\label{berrwt}
\begin{aligned}
\mbox{\textbf{WT:}}\;\;\;\;\;&\\
\epsilon_b(\nu,t)
\equiv&\frac{\big|W_s(\nu,t)-\widetilde{W}_s(\nu,t)\big|}{\left|W_s(\nu,t)\right|}\\
=&\frac{\left|\int e^{i\nu\tau}\psi^*\Big(\frac{\nu(\tau-t)}{\omega_\psi}\Big)\frac{\nu d\tau}{\omega_\psi}
-\int_0^T e^{i\nu\tau}\psi^*\Big(\frac{\nu(\tau-t)}{\omega_\psi}\Big)\frac{\nu d\tau}{\omega_\psi}\right|}
{\left|\int e^{i\nu\tau}\psi^*\Big(\frac{\nu(\tau-t)}{\omega_\psi}\Big)\frac{\nu d\tau}{\omega_\psi}\right|}\\
\leq& \frac{\big|\int_{-\infty}^{0}\psi^*\Big(u-\frac{\nu t}{\omega_\psi}\Big)e^{i\omega_\psi u}du\big|
+\big|\int_{T}^{\infty}\psi^*\Big(u-\frac{\nu t}{\omega_\psi}\Big)e^{i\omega_\psi u}du\big|}
{\big|\int \psi^*\Big(u-\frac{\nu t}{\omega_\psi}\Big)e^{i\omega_\psi u}du\big|}\\
=&|P_\psi(-\nu t/\omega_\psi)|+|P_\psi(\nu (t-T)/\omega_\psi)|,\\
\end{aligned}
\end{equation}
where $\widetilde{G}_s(\omega,t)$ and $\widetilde{W}_s(\omega,t)$ denote respectively the WFT and WT obtained in practice, while $P_{g,\psi}$ are as defined in (\ref{eswft}), (\ref{eswt}). Note that $\epsilon_b(\omega,t)$ (\ref{berrwft}), (\ref{berrwt}) does not depend on the phase-shift $\varphi$ of the tone, as would be the case for the symmetric and periodic padding schemes.

\begin{figure*}[t!]
\includegraphics[width=1.0\linewidth]{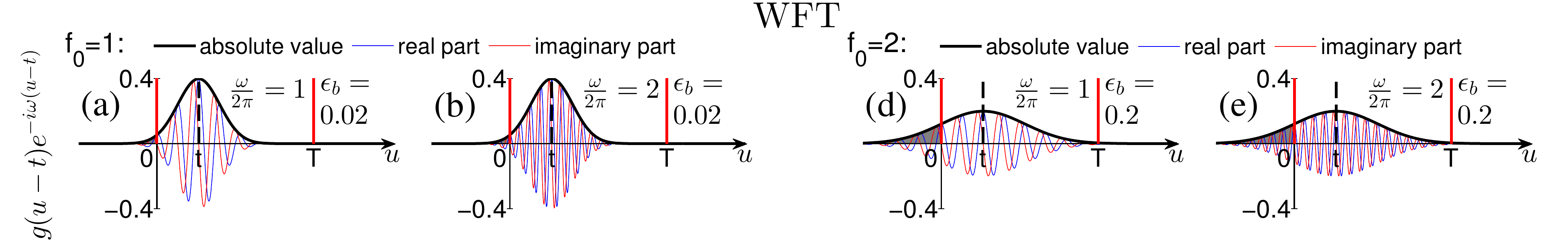}\\
\includegraphics[width=1.0\linewidth]{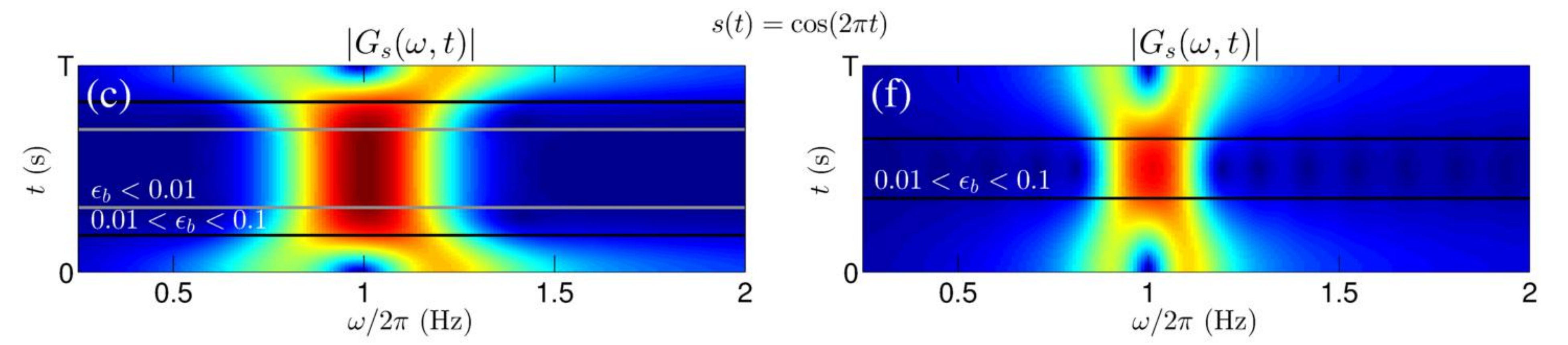}\\
\caption{(a,b,d,e): Examples of the time domain Gaussian window functions (\ref{gw}) (multiplied by the corresponding oscillatory part), with which the signal is convoluted while constructing the WFT (\ref{wft}), for different resolution parameters $f_0$ and frequencies $\omega$; thick vertical lines indicate the time limits of the signal (solid red) and the time $t$ at which the window is centered (dashed black); the area under the window modulus outside the time limits is filled with gray (note that, for a Gaussian window (\ref{gw}), $|g(t)|=g(t)$). (c,f): WFTs for the single tone signal $s(t)=\cos(2\pi t)$ in the same time limits $[0,T]$ as in (a,b,d,e), with black and gray lines enclosing the time-frequency regions where boundary errors (\ref{berrwft}) are small: $\epsilon_b<0.1$ and $\epsilon_b<0.01$, respectively. The signal was sampled at $f_s=100$ Hz for $T=7.5$ s, and in (c,f) padding with zeros was used.}
\label{fig:coiwft}
\end{figure*}

\begin{figure*}[t!]
\includegraphics[width=1.0\linewidth]{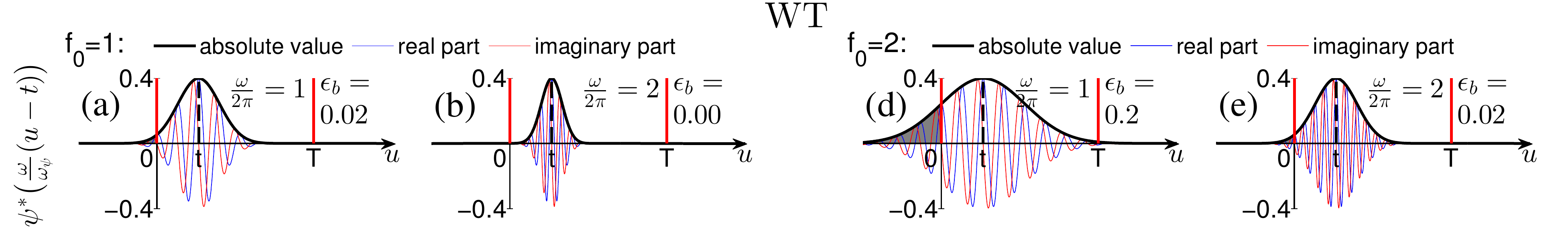}\\
\includegraphics[width=1.0\linewidth]{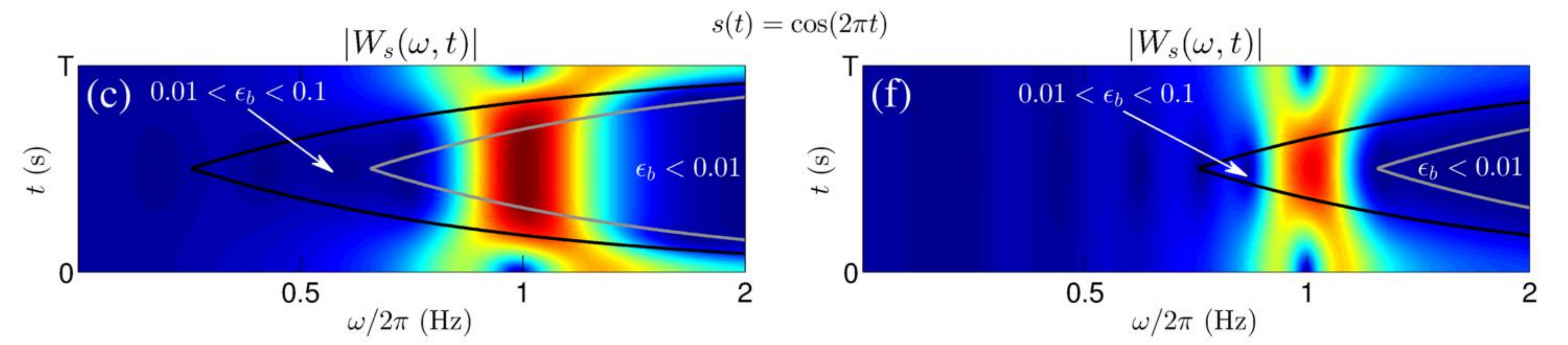}\\
\caption{Same as Fig.\ \ref{fig:coiwft}, but for the WT (\ref{wt}) with Morlet wavelet (\ref{mw}). (a,b,d,e): the wavelet function with which the signal is convoluted while constructing the WT (\ref{wt}); the area under the wavelet modulus outside the time limits is filled with gray (note that $|\psi(t)|\approx \psi(t)e^{-i\omega_\psi t}$ in the present case). (c,f): WTs for the single tone signal $s(t)=\cos(2\pi t)$ within the same time limits $[0,T]$ as in (a,b,d,e).}
\label{fig:coiwt}
\end{figure*}

It is clear that the boundary distortions in the TFR manifest themselves when a non-negligible portion of the window $g(\tau-t)e^{-i\omega(\tau-t)}$ or wavelet $\psi^*(\omega(\tau-t)/\omega_\psi)$, by which the signal $s^+(\tau)$ is multiplied and integrated in (\ref{wft}), (\ref{wt}), lies outside the time limits. This is true for any padding but, in the case of padding with zeros, the related errors take the simple forms (\ref{berrwft}), (\ref{berrwt}). From the latter it can be seen, that the dependence of the boundary inaccuracies on time and frequency is different for the WFT and WT, as illustrated in Figs.\ \ref{fig:coiwft} and \ref{fig:coiwt}.

In the case of the WFT (Fig.\ \ref{fig:coiwft}), the boundary inaccuracies are independent of the frequency $\omega$ and depend only on time $t$: $\epsilon_b(\omega,t)=\epsilon_b(t)$. This is because $\omega$ controls the effective number of oscillations within the window, but not its spread in time, which is determined by the resolution parameter $f_0$. Hence, for each frequency the relative part of $g(\tau-t)$ lying outside the time limits (gray-shaded areas in Fig.\ \ref{fig:coiwft}(a,b,d,e)) is the same, so that for some (central) time range the WFT is well-behaved at all frequencies, as shown in Fig. \ref{fig:coiwft}(c,f). For the WT (Fig. \ref{fig:coiwt}), on the other hand, the wavelet is rescaled at each frequency $\omega$, and so is the part of the wavelet outside the time limits (gray-shaded areas in Fig.\ \ref{fig:coiwt}(a,b,d,e)). As a result, the boundary errors depend on both frequency and time, and the region where the WT coefficients are determined with some predefined accuracy in terms of boundary effects takes the form of a cone, as shown in Fig.\ \ref{fig:coiwt}(c,f).

The time-frequency region where the boundary errors are small in the current TFR is called the \emph{cone-of-influence} (although this region has a truly conic form only for the WT):
\begin{equation}\label{coidef}
\mbox{\textbf{Cone-of-influence}}\;\{\omega,t\}_{coi}^{(\epsilon)}:\;
(\omega,t)\in\{\omega,t\}_{coi}^{(\epsilon)}\Leftrightarrow \epsilon_b(\omega,t)\leq\epsilon,\\
\end{equation}
where $\epsilon$ is the chosen accuracy threshold. Examples of the cones-of-influence for $\epsilon=0.1$ and $\epsilon=0.01$ have already been presented in Fig.\ \ref{fig:coiwft}(c,f) and Fig.\ \ref{fig:coiwt}(c,f), where they are enclosed by black and gray lines, respectively. Based on the error estimates (\ref{berrwft}), (\ref{berrwt}), the cones-of-influence for the WFT and WT can be expressed through the corresponding window/wavelet $\epsilon$-supports in time (\ref{eswft}), (\ref{eswt}) as
\begin{equation}\label{coi}
\begin{aligned}
\mbox{\textbf{WFT:}}\;&\{\omega,t\}_{coi}^{(\epsilon)}
=\Big\{\big[\omega_{\min},\omega_{\max}\big],\big[\delta t_1(\epsilon),T-\delta t_2(\epsilon)\big]\Big\},\\
&\delta t_1(\epsilon)\leq -\tau_1(\epsilon),\quad
\delta t_2(\epsilon)\leq \tau_2(\epsilon),\\
\mbox{\textbf{ WT:}}\;&\{\omega,t\}_{coi}^{(\epsilon)}
=\Big\{\big[\omega_{\min},\omega_{\max}\big],\big[\delta t_1(\epsilon,\omega),T-\delta t_2(\epsilon,\omega)\big]\Big\},\\
&\delta t_1(\epsilon,\omega)\leq -\omega_\psi\tau_1(\epsilon)/\omega,\quad
\delta t_2(\epsilon,\omega)\leq \omega_\psi\tau_2(\epsilon)/\omega,\\
\end{aligned}
\end{equation}
where $[\omega_{\min},\omega_{\max}]$ is the frequency range in which the WFT/WT is calculated (the restrictions on this range and related issues are discussed in Sec.\ \ref{sec:pract4} below).

The cone-of-influence for the SWFT/SWT differs from that for the underlying WFT/WT, being a ``trimmed'' version of the latter as well as shifted up in frequency. This is because, when performing the synchrosqueezing, one sums many WFT/WT coefficients into one SWFT/SWT coefficient, and so needs to account for most of the peak supports in frequency. For example, in the case of a single tone signal $s(t)=\cos(\nu t)$, the WFT/WT calculated in the frequency range $[\omega_{\min},\omega_{\max}]$ will contain the $|R_g(\omega_{\max}-\nu)-R_g(\omega_{\min}-\nu)|$ (WFT) or $|R_\psi(\omega_\psi\nu/\omega_{\min})-R_\psi(\omega_\psi\nu/\omega_{\max})|$ (WT) part of the tone's total amplitude, where $R_{g,\psi}(\xi)$ are defined in (\ref{eswft}), (\ref{eswt}). Hence, it is exactly this part that will be accounted for in the process of synchrosqueezing. Consequently, for any SWFT/SWT coefficient to be accurate (in terms of time and frequency boundary effects) its associated WFT/WT should be well-defined not only at the corresponding time-frequency point, but also over a sufficiently broad frequency range around it. Taking into account the cones-of-influence (\ref{coi}) for the WFT and WT, and using the notation of window/wavelet $\epsilon$-supports (\ref{eswft}), (\ref{eswt}), the cones-of-influence for synchrosqueezed transforms can be shown to be
\begin{equation}\label{sscoi}
\begin{aligned}
\mbox{\textbf{SWFT:}}\;&
\{\omega,t\}_{coi}^{(\epsilon)}=
\Big\{
\big[\omega_{\min}+\delta\omega_1(\epsilon),\omega_{\max}-\delta\omega_2(\epsilon)\big],\\
&\hphantom{\{\omega,t\}_{coi}^{(\epsilon)}=\Big\{}\quad
\big[\delta t_1(\epsilon),T-\delta t_2(\epsilon)\big]
\Big\},\\
&\delta\omega_1(\epsilon)\leq-\xi_1(\epsilon),\quad
\delta\omega_2(\epsilon)\leq\xi_2(\epsilon),\\
\mbox{\textbf{ SWT:}}\;&
\{\omega,t\}_{coi}^{(\epsilon)}=\Big\{
\big[\omega_{\min}r_2(\epsilon),\omega_{\max}/r_1(\epsilon)\big],\\
&\hphantom{\{\omega,t\}_{coi}^{(\epsilon)}=\Big\{}\quad
\big[\delta t_1(\epsilon,\omega/r_2(\epsilon)),T-\delta t_2(\epsilon,\omega/r_2(\epsilon))\big]
\Big\},\\
&r_1(\epsilon)\leq\xi_2(\epsilon)/\omega_\psi,\quad
r_2(\epsilon)\leq\omega_\psi/\xi_1(\epsilon),
\end{aligned}
\end{equation}
where $\delta t_{1,2}(\epsilon)$ are the same as in (\ref{coi}). From (\ref{sscoi}) it is also clear that, to calculate the SWFT/SWT within some frequency range, the underlying WFT/WT should be calculated over a wider range, an issue that will be discussed in more detail in Remark \ref{rem:omom} below.

Because boundary effects can greatly influence the TFR behavior, especially for the (S)WT, it is recommended that all characteristics (e.g.\ the WFT/WT mean amplitudes) be calculated using only TFR coefficients inside the cone-of-influence. The most appropriate padding scheme in this case is zero padding, as the corresponding boundary effects have a universal and well-defined form. In fact, the boundary errors (\ref{berrwft}), (\ref{berrwt}) and the cones-of-influence (\ref{coi}), (\ref{sscoi}) were rigorously estimated exclusively for this type of padding.

However, if one wants to extract some AM/FM component from the signal's TFR (see Sec.\ \ref{sec:tfsrec}), then it should obviously be extracted for all time, and the consideration cannot be restricted to the cone-of-influence only. If all the TFR coefficients are to be used, then zero padding does not represent a good choice, because the (S)WFT/(S)WT near the boundaries will \emph{surely} contain considerable errors. In this case predictive padding would be the most suitable, as it usually has the best performance in terms of reducing boundary effects.

\subsection{TFR frequency range}\label{sec:pract4}

For completeness, the restrictions on the frequency range $[\omega_{\min},\omega_{\max}]$ over which to calculate the TFR, i.e.\ how high/low in frequency one can in principle go, should be also discussed. Consider a signal $s(t)$ sampled at $f_s$ Hz for $T$ seconds. Then one has
\begin{equation}\label{mrange}
\omega_{\min}/2\pi\geq1/T,\quad \omega_{\max}/2\pi\leq f_s/2.
\end{equation}
The restriction on $\omega_{\max}$ follows from the Nyquist theorem, which states that oscillations with frequencies higher than half of the sampling frequency cannot be represented in a discrete signal. The constraints on $\omega_{\min}$, on the other hand, are based on simple logic: it is clear that for a particular oscillation to be reliably studied (generally by any method), there should be at least one of its cycles within the signal.

However, from the statistical viewpoint, to reach any meaningful conclusions about the characteristics of the oscillatory process, such as its typical amplitude (as calculated e.g.\ from the time-averaged TFR amplitude), there should be at least 5-6 corresponding cycles within the signal \cite{Keselbrener:96}, so that
\begin{equation}\label{omstat}
\omega_{\min}^{\rm (stat)}/2\pi=5/T.
\end{equation}
Thus, although one can in principle estimate the properties of the components at lower frequencies, the resultant estimates might be highly untypical. Indeed, the number of cycles can be associated with a number of trials. For example, when testing the effects of some drug, one cannot base conclusions on only one subject (since the probability that the subject tested is an outlier is quite high); but if the same effects appear in 5-6 subjects, this suggests that they are quite common. Note also that, for statistical comparison of oscillatory properties between different data (e.g.\ as is done in \cite{Iatsenko:cardio,Shiogai:10}), all these properties should be statistically meaningful in themselves, so that the related oscillations have frequencies $\geq\omega_{\min}^{\rm (stat)}$.

Finally, since the cone-of-influence for the (S)WT contracts towards the lower frequencies (see Sec.\ \ref{sec:pract3}), there exists also the minimal frequency $\omega_{\min}^{(\epsilon)}$ from the viewpoint of precision (while there is no such restriction in the case of the (S)WFT). The latter can be defined as the minimum frequency for which at least one (S)WT coefficient is determined with accuracy $\epsilon$ in terms of the boundary effects. Using the estimates (\ref{coi}), (\ref{sscoi}) one obtains
\begin{equation}\label{omeps}
\begin{aligned}
\mbox{\textbf{ WFT:}}\;&\omega_{\min}^{(\epsilon)}=-\infty;\quad\;\;
\mbox{\textbf{  WT:}}\;\omega_{\min}^{(\epsilon)}\leq \omega_\psi\big[\tau_2(\epsilon)-\tau_1(\epsilon)\big]/T;\\
\mbox{\textbf{SWFT:}}\;&\omega_{\min}^{(\epsilon)}=-\infty;\quad
\mbox{\textbf{ SWT:}}\;\omega_{\min}^{(\epsilon)}\leq \xi_2(\epsilon)\big[\tau_2(\epsilon)-\tau_1(\epsilon)\big]/T.
\end{aligned}
\end{equation}

Summarizing, it is quite safe to calculate and analyse the TFR in the range $[\omega_{\min},\omega_{\max}]$ with $\omega_{\min}>\max[2\pi/T,\omega_{\min}^{(\epsilon)}]$ and $\omega_{\max}<2\pi f_s/2$. However, if one needs to draw conclusions about the typical oscillatory parameters of the data, they should be based only on the TFR within the cone-of-influence (see Sec.\ \ref{sec:pract3}) and on at least 5 cycles of the related oscillations, so that $\omega_{\min}\geq\max[2\pi\times 5/T,\omega_{\min}^{(\epsilon)}]$.

\begin{remark}\label{rem:omom}
Even if one is interested in the frequency range $[\omega_{\min},\omega_{\max}]$ (so that the signal is filtered in this range during the preprocessing stage as discussed in Sec.\ \ref{sec:pract1}), it may be useful to calculate the TFR over a larger interval (up to $(-\infty,\infty)$ for the WFT and $(0,\infty)$ for the WT). This allows one to trace the ``tails'' of the TFR amplitude peaks located in the original range, which might be needed e.g.\ for reliable parameter reconstruction by the direct method (see Sec.\ \ref{sec:tfsrec}). For example, in Fig.\ \ref{fig:preprocess}(f,i) the considered frequency band -- between the gray vertical lines -- contains only the main part of the corresponding peak, with its tails occupying a wider range. If one wants to encompass the supports of all possible peaks contained in $[\omega_{\min},\omega_{\max}]$ then, while bandpass filtering the signal in the original frequency band, its WFT/WT should be calculated within a slightly larger region $[\wt{\omega}_{\min},\wt{\omega}_{\max}]$. To achieve relative accuracy $\epsilon$, so that the new range contains the $(1-\epsilon)$ part of the total area under any peak in the range of interest, one should use
\begin{equation}\label{omom}
\begin{aligned}
\mbox{\rm\textbf{WFT:}}\;&[\wt{\omega}_{\min},\wt{\omega}_{\max}]
=[\omega_{\min}+\xi_1(\epsilon),\omega_{\max}+\xi_2(\epsilon)],\\
\mbox{\rm\textbf{WT:}}\;&[\wt{\omega}_{\min},\wt{\omega}_{\max}]
=\Big[\omega_{\min}\frac{\xi_1(\epsilon)}{\xi_2(\epsilon)},\omega_{\max}\frac{\xi_2(\epsilon)}{\xi_1(\epsilon)}\Big],\\
\end{aligned}
\end{equation}
where the $\epsilon$-supports $\xi_{1,2}(\epsilon)$ for the window and wavelet are defined in (\ref{eswft}) and (\ref{eswt}), respectively. This issue applies especially to the SWFT/SWT (\ref{stfr}) which, by construction, tries to collect and join all WFT/WT coefficients associated with the same component. Thus, if some coefficients are not accounted for, e.g.\ are lying outside the frequency range considered in the WFT/WT, such a ``missing'' part is subtracted from the corresponding SWFT/SWT entries (see also the related discussion in Sec.\ \ref{sec:pract3B}). Therefore, to calculate the SWFT/SWT in the range $[\omega_{\min},\omega_{\max}]$, one needs to synchrosqueeze the WFT/WT calculated for a wider range $[\wt{\omega}_{\min},\wt{\omega}_{\max}]$ given by (\ref{omom}).
\end{remark}

\section{Summary: step-by-step algorithms and codes}\label{sec:stepalg}

MatLab codes for computing the different TFRs considered in this work (as well as for extracting the components from them) are available at \cite{freecodes} together with detailed documentation and video-instructions. They allow for the use of any window/wavelet, both standard and user-defined (either in time or frequency domain, or both), and include implementation of all the aspects discussed in previous sections, e.g.\ signal preprocessing, automatic determination of the frequency bin widths and the cone-of-influence, different padding schemes (with predictive padding being default) etc. Below we consider each algorithm in detail, summarizing all related issues and steps.

In the present section it is assumed that the original signal $s(t_n)$ is sampled at $t_n=(n-1)\Delta t$, $n=1,..,N$ (so its sampling frequency and the overall time-duration are $f_s=1/\Delta t$ and $T=(N-1)\Delta t$, respectively), and the signal's TFR needs to be calculated for frequencies $\omega\in[\omega_{\min},\omega_{\max}]$. In what follows, $\epsilon$ will denote the fixed relative precision for determination of the cone-of-influence and padding (see Sec.\ \ref{sec:pract3}), while $[x]^{\uparrow}$ and $[x]^{\downarrow}$ will stand for rounding of $x$ up and down, respectively (e.g.\ $[1.3]^{\uparrow}=2$, $[1.3]^{\downarrow}=1$).

The discrete FT of the signal calculated using the FFT algorithm will be denoted as
\begin{equation}\label{sigft}
\hat{s}(\xi_n):\;[\hat{s}(\xi_1),...,\hat{s}(\xi_N)]={\rm FFT}[s(t_1),...,s(t_N)],
\end{equation}
where usually (e.g.\ in the MatLab FFT implementation) the FT is returned at frequencies $\xi_n$ in the order
\begin{equation}\label{ftfreq}
\begin{aligned}
\big[\xi_1,...,\xi_N\big]=\big[&0,2\pi f_s/N,...,2\pi[(N-1)/2]^{\uparrow}f_s/N,\\
&-2\pi[(N-1)/2]^{\downarrow}f_s/N,...,-2\pi f_s/N\big],
\end{aligned}
\end{equation}
so that the positive frequencies come first, and only then the negative ones. The original signal can be recovered by applying the inverse fast Fourier transform (IFFT) algorithm to the signal's discrete FT: $[s(t_1),...,s(t_N)]={\rm IFFT}[\hat{s}(\xi_1),...,\hat{s}(\xi_N)]$. It should also be noted, that throughout this work the literals $\omega$ and $\xi$ (e.g.\ $\omega_{\min,\max}$) denote the circular frequency (in rad/s), with only the sampling frequency $f_s$, where present, being taken in Hz; in the codes \cite{freecodes}, however, all frequencies are in Hz, so that e.g.\ the frequencies for which the TFR is calculated are returned as $f_k=\omega_k/2\pi$.

\begin{remark}
In some steps of the TFR algorithms one will need to calculate the $\epsilon$-supports (\ref{eswft}), (\ref{eswt}), as well as $\omega_\psi$ (\ref{wtsc}) for the WT. The corresponding formulas assume that window/wavelet functions are known both in time and frequency, but this is not always so, e.g.\ the explicit form might be available for $g(t)$ only, and not for $\hat{g}(\xi)$. In such cases, one will need to calculate numerically the window/wavelet FT or inverse FT, and use it to obtain the values needed. The implementation of this approach for the general case is quite cumbersome and brings some numerical issues. It will not be further discussed here, but it is included in the codes \cite{freecodes}.
\end{remark}

\subsection{WFT $G_s(\omega_k,t_n)$}\label{sec:stepalgA}

\begin{enumerate}[leftmargin=*] \parskip4pt
\item Preprocess a signal as described in Sec.\ \ref{sec:pract1}: first subtract a third-order polynomial fit from it, and then bandpass-filter what is left in the frequency band $[\omega_{\min},\omega_{\max}]$.
\item Select a padding scheme (see Sec.\ \ref{sec:pract3}) and accordingly pad the preprocessed signal at both ends with $n_1$ values to the left and $n_2$ to the right, where $n_{1,2}$ are determined by (\ref{numpad}). Denote the padded signal as $s_p(\tilde{t}_j)$ with $\tilde{t}_j=(j-1-n_1)\Delta t$ and $j=1,...,N_p=N+n_1+n_2$.
\item Break the frequency interval $[\omega_{\min},\omega_{\max}]$ into bins $\omega_k=(k-k_0)\Delta\omega$, where $k_0=1-[\frac{\omega_{\min}}{\Delta\omega}]^{\uparrow}$ and $k=1,...,k_0+[\frac{\omega_{\max}}{\Delta\omega}]^{\downarrow}$; the optimal frequency step $\Delta\omega$ can be determined by (\ref{fbins}) (in the latter, $N_b=10$ is used by default). Note that the positions of $\omega_k$ do not depend on any signal parameters (e.g.\ $f_s$ or $T$), as is sometimes the case, being determined only by $\Delta\omega$: this is convenient since it allows WFTs of different signals to be calculated at the same frequencies, so that their characteristics (e.g.\ mean amplitudes at each frequency) can be compared.
\item Calculate the FT of a padded signal: $[\hat{s}_p(\tilde{\xi}_1),...,\hat{s}_p(\tilde{\xi}_{N_p})]={\rm FFT}[s_p(\tilde{t}_1),...,s_p(\tilde{t}_{N_p})]$, where frequencies $\tilde{\xi}_j$ are given by (\ref{ftfreq}) with $N\rightarrow N_p,n\rightarrow j$. Set the FT at negative frequencies to zero: $\hat{s}(\tilde{\xi}_j\leq0)=0$.
\item For each frequency $\omega_k$:

\textbf{(a)} If the explicit form of the window function FT is available (e.g.\ for the Gaussian window (\ref{gw})), calculate $\hat{g}(\omega_k-\tilde{\xi}_j)$ by direct substitution of the arguments $\omega_k-\tilde{\xi}_j$ into the known $\hat{g}(\xi)$. Otherwise, if only $g(t)$ is available, compute $[\hat{g}(\omega_k-\tilde{\xi}_1),...,\hat{g}(\omega_k-\tilde{\xi}_{N_p})]=\Delta t\,{\rm FFT}[g(\tau_1)e^{-i\omega_k\tau_1},...,g(\tau_{N_p})e^{-i\omega_k\tau_{N_p}}]$, where $\tau_j=-(j-1)\Delta t$ for $j=1,..,[(N_p-1)/2]^{\uparrow}$ and $\tau_j=(N_p-j+1)\Delta t$ for other $j$.

\textbf{(b)} According to the frequency domain form of (\ref{wft}), calculate the convolutions $c(\tilde{t}_j)$ of the padded signal $s_p^+(u)$ with $g(u-\tilde{t}_j)e^{-i\omega_k(u-\tilde{t}_j)}$ as: $[c(\tilde{t}_1),...,$ $c(\tilde{t}_{N_p})]={\rm IFFT}[\hat{s}_p(\tilde{\xi}_1)\hat{g}(\omega_k-\tilde{\xi}_1),...,\hat{s}_p(\tilde{\xi}_{N_p})\hat{g}(\omega_k-\tilde{\xi}_{N_p})]$, where $\hat{g}(\omega_k-\tilde{\xi}_j)$ were obtained in the previous substep (a), while $\hat{s}_p(\tilde{\xi}_j)$ were calculated in step 4 (and are zero at negative frequencies).

\textbf{(c)} The WFT at frequency $\omega_k$ is then equal to $c(\tilde{t}_j)$ for times within the original signal's time limits: $G_s(\omega_k,t_{n=1,...,N})=c(\tilde{t}_{j=1+n_1,...,N_p-n_2})$.
\end{enumerate}

\subsection{WT $W_s(\omega_k,t_n)$}\label{sec:stepalgB}

\begin{enumerate}[leftmargin=*] \parskip4pt
\item Preprocess the signal as described in Sec.\ \ref{sec:pract1}: first subtract a third-order polynomial fit from it, and then bandpass-filter what is left in the frequency band $[\omega_{\min},\omega_{\max}]$.
\item Select a padding scheme (see Sec.\ \ref{sec:pract3}) and accordingly pad the preprocessed signal at both ends with $n_1$ values to the left and $n_2$ to the right, where $n_{1,2}$ are determined by (\ref{numpad}). Denote the padded signal as $s_p(\tilde{t}_j)$ with $\tilde{t}_j=(j-1-n_1)\Delta t$ and $j=1,...,N_p=N+n_1+n_2$.
\item Break the frequency interval $[\omega_{\min},\omega_{\max}]$ into bins $\omega_k/2\pi=2^{(k-k_0)/n_v}$, where $k_0=1-[n_v\log_2\frac{\omega_{\min}}{2\pi}]^{\uparrow}$ and $k=1,...,k_0+[n_v\log_2\frac{\omega_{\max}}{2\pi}]^{\downarrow}$; the optimal number-of-voices $n_v$ can be determined by (\ref{fbins}) (in the latter, $N_b=10$ is used by default). Note that the positions of $\omega_k$ do not depend on any signal parameters (e.g.\ $f_s$ or $T$), as is sometimes the case, being determined only by $n_v$: this is convenient since it allows WTs of different signals to be calculated at the same frequencies, so that their characteristics (e.g.\ mean amplitudes at each frequency) can be compared.
\item Calculate the FT of a padded signal: $[\hat{s}_p(\tilde{\xi}_1),...,\hat{s}_p(\tilde{\xi}_{N_p})]={\rm FFT}[s_p(\tilde{t}_1),...,s_p(\tilde{t}_{N_p})]$, where frequencies $\tilde{\xi}_j$ are given by (\ref{ftfreq}) with $N\rightarrow N_p,n\rightarrow j$. Set the FT at negative frequencies to zero: $\hat{s}(\tilde{\xi}_j\leq0)=0$.
\item For each frequency $\omega_k$:

\textbf{(a)} If the explicit form of the wavelet function FT is available (e.g.\ for the lognormal wavelet (\ref{lw})), calculate $\hat{\psi}^*{\Big(}\frac{\omega_\psi\tilde{\xi}_j}{\omega_k}{\Big)}$ by direct substitution of the arguments $\omega_\psi\tilde{\xi}_j/\omega_k$ into the known $\hat{\psi}^*(\xi)$. Otherwise, if only $\psi(t)$ is available, compute $\Big[\hat{\psi}^*{\Big(}\frac{\omega_\psi\tilde{\xi}_1}{\omega_k}{\Big)},$ $...,
\hat{\psi}^*{\Big(}\frac{\omega_\psi\tilde{\xi}_{N_p}}{\omega_k}{\Big)}\Big]
=\Delta t\,{\rm FFT}\Big[\psi^*{\Big(}\frac{\omega_k\tau_1}{\omega_\psi}{\Big)},...,
\psi^*{\Big(}\frac{\omega_k\tau_{N_p}}{\omega_\psi}{\Big)}\Big]$,
where $\tau_j=-(j-1)\Delta t$ for $j=1,..,[(N_p-1)/2]^{\uparrow}$ and $\tau_j=(N_p-j+1)\Delta t$ for other $j$.

\textbf{(b)} According to the frequency domain form of (\ref{wt}), calculate the convolutions $c(\tilde{t}_j)$ of the padded signal $s_p^+(u)$ with $\frac{\omega}{\omega_\psi}\psi^*\Big(\frac{\omega_k(u-\tilde{t}_j)}{\omega_\psi}\Big)$ as: $[c(\tilde{t}_1),...,c(\tilde{t}_{N_p})]={\rm IFFT}\big[\hat{s}_p(\tilde{\xi}_1)\hat{\psi}^*(\omega_\psi\tilde{\xi}_1/\omega_k),$ $...,
\hat{s}_p(\tilde{\xi}_{N_p})\hat{\psi}^*(\omega_\psi\tilde{\xi}_{N_p}/\omega_k)\big]$, where $\hat{\psi}^*(\omega_\psi\tilde{\xi}_j/\omega_k)$ were obtained in the previous substep (a), while $\hat{s}_p(\tilde{\xi}_j)$ were calculated in step 4 (and are zero at negative frequencies).

\textbf{(c)} The WT at frequency $\omega_k$ is then equal to $c(\tilde{t}_j)$ for times within the original signal's time limits: $W_s(\omega_k,t_{n=1,...,N})=c(\tilde{t}_{j=1+n_1,...,N_p-n_2})$.
\end{enumerate}

\subsection{Synchrosqueezing (SWFT $\widetilde{V}_s(\omega_k,t_n)$ and SWT $\widetilde{T}_s(\omega_k,t_n)$)}\label{sec:stepalgC}

Here we present the fast algorithm for synchrosqueezing, which is based on a slightly modified idea of \cite{Thakur:13}, complemented with the utilization of sparse matrices. The latter are used because synchrosqueezed transforms usually have only a very small percentage of nonzero entries, so the sparse matrix representation is very suitable for them. It should be noted, that synchrosqueezing (at least performed according to a given algorithm) is not computationally expensive, as is sometimes thought: its computational cost is of the same order as for the underlying WFT/WT transforms. The algorithm is as follows:

\begin{enumerate}[leftmargin=*] \parskip4pt
\item Preprocess the signal as described in Sec.\ \ref{sec:pract1}: first subtract third-order polynomial fit from it, and then bandpass-filter what is left in the frequency band $[\omega_{\min},\omega_{\max}]$.
\item Calculate WFT $G_s(\wt{\omega}_i,t_n)$ or WT $W_s(\wt{\omega}_i,t_n)$ as described previously, but for a wider frequency interval (see Sec.\ \ref{sec:pract4}): $\wt{\omega}_i\in[\omega_{\min}+\xi_1(\epsilon),\omega_{\max}+\xi_2(\epsilon)]$ for the WFT or $\wt{\omega}_i\in\Big[\frac{\omega_{\min}\omega_\psi}{\xi_2(\epsilon)}$, $\frac{\omega_{\max}\omega_\psi}{\xi_1(\epsilon)}\Big]$ for the WT, where $\xi_{1,2}(\epsilon)$ are defined in (\ref{eswft}), (\ref{eswt}).
\item Compute the phases of the WFT $\phi_{G}(\wt{\omega}_i,t_n)\equiv{\rm arg}[G_s(\wt{\omega}_i,t_n)]$ or WT $\phi_{W}(\wt{\omega}_i,t_n)\equiv{\rm arg}[W_s(\wt{\omega}_i,t_n)]$, and at each $\wt{\omega}_i,t_n$ unwrap them over time. Then estimate the TFR phase velocities (\ref{iftfr}) as
\begin{equation*}
\begin{aligned}
\nu_{G,W}(\wt{\omega}_i,t_n)=&[\phi_{G,W}(\wt{\omega}_i,t_{n+1})-\phi_{G,W}(\wt{\omega}_i,t_{n-1})]/(2\Delta t)\\
&\mbox{ for }n=2,...,N-1,\\
\nu_{G,W}(\wt{\omega}_i,t_1)=&[\phi_{G,W}(\wt{\omega}_i,t_2)-\phi_{G,W}(\wt{\omega}_i,t_1)]/\Delta t,\\
\nu_{G,W}(\wt{\omega}_i,t_N)=&[\phi_{G,W}(\wt{\omega}_i,t_N)-\phi_{G,W}(\wt{\omega}_i,t_{N-1})]/\Delta t.\\
\end{aligned}
\end{equation*}
Alternatively, one can estimate $\nu_{G,W}(\wt{\omega}_i,t_n)$ as
\begin{equation*}
\nu_{G}(\wt{\omega}_i,t_n)=\operatorname{Im}\Bigg[\frac{\partial_t G_s(\wt{\omega}_i,t_n)}{G_s(\wt{\omega}_i,t_n)}\Bigg],\;
\nu_{W}(\wt{\omega}_i,t_n)=\operatorname{Im}\Bigg[\frac{\partial_t W_s(\wt{\omega}_i,t_n)}{W_s(\wt{\omega}_i,t_n)}\Bigg],
\end{equation*}
with $\partial_t G_s(\wt{\omega}_i,t_n)$ and $\partial_t W_s(\wt{\omega}_i,t_n)$ being calculated as the usual WFT and WT, but using in the final steps 5(a-c) of the previously outlined algorithms $\hat{g}(\omega-\xi)\rightarrow i\xi\hat{g}(\omega-\xi)$ and $\hat{\psi}^*(\omega_\psi\xi/\omega)\rightarrow i\xi \hat{\psi}^*(\omega_\psi\xi/\omega)$, respectively; such a way of calculating them follows from the expressions for time derivatives of the frequency-domain forms of (\ref{wft}) and (\ref{wt}). This latter approach is more robust to time discretization effects, but is also more computationally expensive than estimation of phase velocities by direct numerical differentiation.
\item Break the frequency interval $[\omega_{\min},\omega_{\max}]$ into bins $\omega_k$ constructed as discussed previously (for the SWFT/SWT --- as for the WFT/WT). Note, that one might want to use finer frequency binning for synchrosqueezed transforms (i.e.\ lower $\Delta\omega$ or $n_v$ than given by (\ref{fbins})) to better reflect their non-smooth behavior and/or improve accuracy of the SWFT/SWT-based frequency estimation (see Sec.\ \ref{sec:pract2}).
\item According to (\ref{numstfr}), all WFT/WT coefficients corresponding to $\nu_{G,W}(\wt{\omega}_i,t_n)$ located within the frequency bin centered at $\omega_k$ are joined into $\widetilde{V}_s(\omega_k,t_n)$ or $\widetilde{T}_s(\omega_k,t_n)$. Therefore, calculate the numbers $k(\wt{\omega}_i,t_n)$ of the SWFT/SWT frequency bins into which the corresponding WFT/WT coefficients are joined:
\begin{equation*}
\begin{aligned}
\mbox{\textbf{SWFT:}}\;&k(\wt{\omega}_i,t_n)=1+\left[\frac{1}{2}+\frac{\nu_G(\wt{\omega}_i,t_n)-\omega_1}{\Delta\omega}\right]^{\downarrow},\\
\mbox{\textbf{SWT:}}\;&k(\wt{\omega}_i,t_n)=1+\left[\frac{1}{2}+\frac{\log\nu_W(\wt{\omega}_i,t_n)-\log\omega_1}{n_v^{-1}\log 2}\right]^{\downarrow}.\\
\end{aligned}
\end{equation*}
\item Find all pairs of indices $i$ and $n$ so that $1\leq k(\wt{\omega}_i,t_n)\leq N_f$, where $N_f$ is the number of SWFT/SWT frequencies $\omega_k$. After this, for each relevant time-frequency point $p$ one will have a triple of indices $[i_p,n_p,k_p\equiv k(\wt{\omega}_i,t_n)]$.
\item Initialize $N_f\times N$ sparse matrix $S[k,j]=0$ and for all $p$ do
\begin{equation*}
\begin{aligned}
\mbox{\textbf{SWFT:}}\;&S[k_p,n_p]=S[k_p,n_p]+G_s(\wt{\omega}_{i_p},t_{n_p})C_g^{-1}\Delta\omega,\\
\mbox{\textbf{SWT:}}\;&S[k_p,n_p]=S[k_p,n_p]+W_s(\wt{\omega}_{i_p},t_{n_p})C_\psi^{-1}\frac{\log 2}{n_v}.\\
\end{aligned}
\end{equation*}
In MatLab, the repeated entries are automatically added up in the process of construction of a sparse matrix, so for the SWFT the full procedure of this step can be implemented simply as $S=sparse(k,n,Q*H_s(l),N_f,N)$, where $i=[i(1);i(2),...]$, $n=[n(1);n(2);...]$, $k=[k(1);k(2);...]$ correspond to vectors of $i_p,n_p,k_p$, while $Q$ is equal to $C_{g}^{-1}\Delta\omega$ for the WFT and to $C_\psi^{-1}\frac{\log 2}{n_v}$ for the WT; $l=sub2ind([\tilde{N}_f,N],n,j)$ denote linear indices associated with the pairs $[i_p,n_p]$ for $\wt{N}_f\times N$ matrix of the WFT or WT: $H_s[i,n]=\big\{G_s[i,n]\;\mbox{or}\;W_s[i,n]\big\}$ ($\wt{N}_f$ is the number of $\wt{\omega}_i$).
\item Finally, transform $S[k,n]$ into an ordinary matrix form:\newline $\widetilde{V}_s(\omega_k,t_n)={\rm full}(S[k,n])$ or $\widetilde{T}_s(\omega_k,t_n)={\rm full}(S[k,n])$. Alternatively, it can be left in the sparse form if one needs only to save it for further usage, or if the operations needed for the analysis are defined for sparse matrices.
\end{enumerate}

\appendix

\section{Nomenclature}\label{app:nomenclature}

\subsection{Abbreviations}

\begin{center}
\begin{tabular}{@{}p{1.5cm}@{} @{}p{6.5cm}@{}}
ECG & Electrocardiogram
\end{tabular}\\
\vspace{0.1cm}
\begin{tabular}{@{}p{1.5cm}@{} @{}p{6.5cm}@{}}
FT & Fourier Transform
\end{tabular}\\
\vspace{0.1cm}
\begin{tabular}{@{}p{1.5cm}@{} @{}p{6.5cm}@{}}
FFT & Fast Fourier Transform (algorithm)
\end{tabular}\\
\vspace{0.1cm}
\begin{tabular}{@{}p{1.5cm}@{} @{}p{6.5cm}@{}}
IFFT & Inverse Fast Fourier Transform (algorithm)
\end{tabular}\\
\vspace{0.1cm}
\begin{tabular}{@{}p{1.5cm}@{} @{}p{6.5cm}@{}}
WFT & Windowed Fourier Transform (see Sec.\ \ref{sec:tfrWFT})
\end{tabular}\\
\vspace{0.1cm}
\begin{tabular}{@{}p{1.5cm}@{} @{}p{6.5cm}@{}}
WT & Wavelet Transform (see Sec.\ \ref{sec:tfrWT})
\end{tabular}\\
\vspace{0.1cm}
\begin{tabular}{@{}p{1.5cm}@{} @{}p{6.5cm}@{}}
SWFT & Synchrosqueezed WFT (see Sec.\ \ref{sec:tfrSS})
\end{tabular}\\
\vspace{0.1cm}
\begin{tabular}{@{}p{1.5cm}@{} @{}p{6.5cm}@{}}
SWT & Synchrosqueezed WT (see Sec.\ \ref{sec:tfrSS})
\end{tabular}\\
\vspace{0.1cm}
\begin{tabular}{@{}p{1.5cm}@{} @{}p{6.5cm}@{}}
TFR & Time-Frequency Representation (includes WFT, WT, SWFT, SWT and many others, but in this work only the former four are considered)
\end{tabular}\\
\vspace{0.1cm}
\begin{tabular}{@{}p{1.5cm}@{} @{}p{6.5cm}@{}}
TFS & Time-Frequency Support (see Sec.\ \ref{sec:tfsrec})
\end{tabular}
\end{center}

\subsection{Terminology}

\begin{center}
\begin{tabular}{>{\centering}p{3cm}<{\centering} p{5cm}}
AM/FM component, or simply component
&
Sinusoidal oscillation with amplitude and/or frequency modulation, i.e.\ the function of time $t$ of the form $x(t)=A(t)\cos \phi(t)$ with $A(t)>0,\phi'(t)>0,\forall t$ (see Sec.\ \ref{sec:analyticsignal} for a more detailed discussion).
\end{tabular}\\
\vspace{0.1cm}
\begin{tabular}{>{\centering}p{3cm}<{\centering} p{5cm}}
Tone
&
AM/FM component of constant amplitude and frequency, i.e.\ a simple sine $x(t)=A\cos(\nu t+\varphi)$.
\end{tabular}\\
\vspace{0.1cm}
\begin{tabular}{>{\centering}p{3cm}<{\centering} p{5cm}}
Gaussian window
&
Window function of the form (\ref{gw}).\\
\end{tabular}\\
\vspace{0.1cm}
\begin{tabular}{>{\centering}p{3cm}<{\centering} p{5cm}}
Morlet wavelet
&
Wavelet function of the form (\ref{mw}).\\
\end{tabular}\\
\vspace{0.1cm}
\begin{tabular}{>{\centering}p{3cm}<{\centering} p{5cm}}
Lognormal wavelet
&
Wavelet function of the form (\ref{lw}).\\
\end{tabular}\\
\vspace{0.1cm}
\begin{tabular}{>{\centering}p{3cm}<{\centering} p{5cm}}
Time and frequency resolution of the TFR
&
The reciprocal of the minimal time lag and frequency difference between the two time events (delta-peaks) and two frequency events (tones) for which they can both be represented reliably in the TFR; see Sec.\ \ref{sec:tfres}.
\end{tabular}\\
\vspace{0.1cm}
\begin{tabular}{>{\centering}p{3cm}<{\centering} p{5cm}}
Joint time-frequency resolution of the TFR
&
The reciprocal of the minimal time-frequency area where the interference between two time events and between two frequency events is simultaneously small in the TFR; see Sec.\ \ref{sec:tfres}.
\end{tabular}
\end{center}

\subsection{Main notation}

\begin{center}
\begin{tabular}{>{\centering}p{2.2cm}<{\centering} p{6cm}}
$\hat{f}(\xi)$
&
Fourier Transform of $f(t)$, see (\ref{nt}). 
\end{tabular}\\
\vspace{0.1cm}
\begin{tabular}{>{\centering}p{2.2cm}<{\centering} p{6cm}}
$f^{\pm}(t)$
&
Positive/negative frequency part of $f(t)$, see (\ref{nt}). 
\end{tabular}\\
\vspace{0.1cm}
\begin{tabular}{>{\centering}p{2.2cm}<{\centering} p{6cm}}
$\langle f(t)\rangle$
&
Time average of $f(t)$, see (\ref{nt}). 
\end{tabular}\\
\vspace{0.1cm}
\begin{tabular}{>{\centering}p{2.2cm}<{\centering} p{6cm}}
$\std[f(t)]$
&
Standard deviation of $f(t)$, see (\ref{nt}). 
\end{tabular}\\
\vspace{0.1cm}
\begin{tabular}{>{\centering}p{2.2cm}<{\centering} p{6cm}}
$f^*(t)$
&
Complex conjugate of $f(t)$.
\end{tabular}\\
\vspace{0.1cm}
\begin{tabular}{>{\centering}p{2.2cm}<{\centering} p{6cm}}
$c.c.$
&
Complex conjugate of the preceding expression.
\end{tabular}\\
\vspace{0.1cm}
\begin{tabular}{>{\centering}p{2.2cm}<{\centering} p{6cm}}
${\rm Re}[x],{\rm Im}[x]$
&
Real and imaginary parts of $x$.
\end{tabular}\\
\vspace{0.1cm}
\begin{tabular}{>{\centering}p{2.2cm}<{\centering} p{6cm}}
$s^a(t)$
&
Analytic signal (\ref{asig}), which is twice the signal's positive frequency part $s^{a}(t)=2s^{+}(t)$. If the original signal is real, as is assumed in this work, then $s(t)={\rm Re}[s^{a}(t)]$.
\end{tabular}\\
\vspace{0.1cm}
\begin{tabular}{>{\centering}p{2.2cm}<{\centering} p{6cm}}
${\rm sign}(x)$
&
Sign function: ${\rm sign}(x>0)=1$, ${\rm sign}(x<0)=-1$, ${\rm sign}(0)=0$.
\end{tabular}\\
\vspace{0.1cm}
\begin{tabular}{>{\centering}p{2.2cm}<{\centering} p{6cm}}
$J_n(x)$
&
$n^{\rm th}$ order Bessel function of the first kind.
\end{tabular}\\
\vspace{0.1cm}
\begin{tabular}{>{\centering}p{2.2cm}<{\centering} p{6cm}}
$I_n(x)$
&
$n^{\rm th}$ order modified Bessel function of the first kind.
\end{tabular}\\
\vspace{0.1cm}
\begin{tabular}{>{\centering}p{2.2cm}<{\centering} p{6cm}}
$\Gamma(a)$
&
Gamma function $\Gamma(n)\equiv\int_0^\infty x^{a-1}e^{-x}dx$ ($=(a-1)!$ if $a\in\mathds{N}$).
\end{tabular}\\
\vspace{0.1cm}
\begin{tabular}{>{\centering}p{2.2cm}<{\centering} p{6cm}}
${\rm erf}(x)$
&
Gauss error function ${\rm erf}(x)\equiv\frac{2}{\sqrt{\pi}}\int_0^xe^{-u^2}du$.
\end{tabular}\\
\vspace{0.1cm}
\begin{tabular}{>{\centering}p{2.2cm}<{\centering} p{6cm}}
$n_G(\epsilon)$
&
Number of standard deviations of the Gaussian distribution within which its $(1-\epsilon)$ part is contained, i.e.\ $\frac{\int_{-n_G}^{n_G}e^{-u^2/2}du}{\int_{-\infty}^{\infty}e^{-u^2/2}du}={\rm erf}(n_G/\sqrt{2})=1-\epsilon$. For example, $n_G(0.05)\approx 2$, $n_G(0.01)\approx 2.5$, $n_G(0.001)\approx3.3$.\\
\end{tabular}\\
\vspace{0.1cm}
\begin{tabular}{>{\centering}p{2.2cm}<{\centering} p{6cm}}
$G_s(\omega,t),W_s(\omega,t)$
&
WFT (\ref{wft}) and WT (\ref{wt}) of the considered signal $s(t)$, respectively.
\end{tabular}\\
\vspace{0.1cm}
\begin{tabular}{>{\centering}p{2.2cm}<{\centering} p{6cm}}
$g(t),\hat{g}(\xi)$
&
Time domain and frequency domain forms of the window function used for computation of the WFT (\ref{wft}), respectively.
\end{tabular}\\
\vspace{0.1cm}
\begin{tabular}{>{\centering}p{2.2cm}<{\centering} p{6cm}}
$\psi(t),\hat{\psi}(\xi)$
&
Time domain and frequency domain forms of the wavelet function used for computation of the WT (\ref{wt}), respectively.
\end{tabular}\\
\vspace{0.1cm}
\begin{tabular}{>{\centering}p{2.2cm}<{\centering} p{6cm}}
$\omega_g,\omega_\psi$
&
Window and wavelet peak frequencies: $\omega_g\equiv\underset{\omega}{\rm argmax}|\hat{g}(\omega)|$, $\omega_\psi\equiv\underset{\omega}{\rm argmax}|\hat{\psi}(\omega)|$. It is assumed that $\omega_g=0$.
\end{tabular}\\
\vspace{0.1cm}
\begin{tabular}{>{\centering}p{2.2cm}<{\centering} p{6cm}}
$C_g,C_\psi$
&
Integration constants defined in (\ref{iwft}) and (\ref{iwt}), respectively.
\end{tabular}\\
\vspace{0.1cm}
\begin{tabular}{>{\centering}p{2.2cm}<{\centering} p{6cm}}
$\overline{\omega}_g,D_\psi$
&
Integration constants defined in (\ref{D1wft}) and (\ref{D1wt}), respectively.
\end{tabular}\\
\vspace{0.1cm}
\begin{tabular}{>{\centering}p{2.2cm}<{\centering} p{6cm}}
$f_0$
&
Window/wavelet resolution parameter which determines the trade-off between its time and frequency resolutions (in the sense that increasing $f_0$ increases frequency resolution but decreases time resolution). The way it is introduced for different window and wavelet functions is described in \ref{app:winwav}.
\end{tabular}\\
\vspace{0.1cm}
\begin{tabular}{>{\centering}p{2.2cm}<{\centering} p{6cm}}
$R_g(\omega),P_g(\tau)$
&
Quantitative measures of the area below $\hat{g}(\xi<\omega)$ and $g(t<\tau)$, respectively, as defined in (\ref{eswft}).
\end{tabular}\\
\vspace{0.1cm}
\begin{tabular}{>{\centering}p{2.2cm}<{\centering} p{6cm}}
$R_\psi(\omega),P_\psi(\tau)$
&
Quantitative measures of the area below $\hat{\psi}(0<\xi<\omega)$ and $g(t<\tau)$, respectively, as defined in (\ref{eswt}).
\end{tabular}\\
\vspace{0.1cm}
\begin{tabular}{>{\centering}p{2.2cm}<{\centering} p{6cm}}
$\xi_{1,2}(\epsilon)$, $\tau_{1,2}(\epsilon)$
&
$\epsilon$-supports of the window/wavelet function in time and frequency, defined in (\ref{eswft}) for windows and in (\ref{eswt}) for wavelets; whether it refers to the former or to the latter is always clear from the context.
\end{tabular}\\
\vspace{0.1cm}
\begin{tabular}{>{\centering}p{2.2cm}<{\centering} p{6cm}}
$V_s(\omega,t),T_s(\omega,t)$
&
SWFT and SWT (\ref{stfr}) of the considered signal $s(t)$, respectively.
\end{tabular}\\
\vspace{0.1cm}
\begin{tabular}{>{\centering}p{2.2cm}<{\centering} p{6cm}}
$\widetilde{V}_s(\omega,t),\widetilde{T}_s(\omega,t)$
&
Numerical analogs (\ref{numstfr}) of $V_s(\omega,t)$ and $T_s(\omega,t)$, respectively.
\end{tabular}\\
\vspace{0.1cm}
\begin{tabular}{>{\centering}p{2.2cm}<{\centering} p{6cm}}
$[\omega_{\min},\omega_{\max}]$
&
Frequency range for which the currently considered TFR is calculated.
\end{tabular}\\
\vspace{0.1cm}
\begin{tabular}{>{\centering}p{2.2cm}<{\centering} p{6cm}}
$\omega_p(t)$
&
Ridge curve of some component in the currently considered TFR (see Sec.\ \ref{sec:tfsrec}).
\end{tabular}\\
\vspace{0.1cm}
\begin{tabular}{>{\centering}p{2.2cm}<{\centering} p{6cm}}
$[\omega_-(t),\omega_+(t)]$
&
Time-frequency support of some component in the currently considered TFR (see Sec.\ \ref{sec:tfsrec}).
\end{tabular}
\end{center}

\subsection{Assumptions and conventions}

\begin{itemize}[leftmargin=*]

\item Where undefined, the integrals are taken over $(-\infty,\infty)$ or, in practice, over the full range of the corresponding variable. For example, given a signal $s(t)$, one has $\int s(t)dt\equiv\int_{-\infty}^{\infty}s(t)dt$ theoretically and $\int s(t)dt\equiv\int_{0}^{T}s(t)dt$ practically, with $T$ denoting the overall time duration of the signal.

\item All TFRs are assumed to be computed for a real signal $s(t)$. Additionally, without loss of generality it is always assumed that $\underset{\xi}{\rm argmax}|\hat{g}(\xi)|=0$, i.e.\ that the window function used for calculating the (S)WFT (\ref{wft}) is peaked at zero in the frequency domain (see Sec.\ \ref{sec:tfrWFT}).

\item In all examples, unless specified otherwise, the (S)WFT and (S)WT are calculated using the Gaussian window (\ref{gw}) and Morlet wavelet (\ref{mw}), respectively, with the resolution parameter set to $f_0=1$. The frequency axis is discretized as $\omega_k=(k-k_0)\Delta\omega$ for the WFT, and as $\omega_k/2\pi=2^{(k-k_0)/n_v}$ for the WT, with the discretization parameters $\Delta\omega$ and $n_v$ being selected according to the criteria (\ref{fbins}). To reduce boundary effects in the TFRs (see Sec.\ \ref{sec:pract3}), predictive padding is used by default, with the number of padded values being calculated according to (\ref{numpad}); in the majority of cases, however, only the ``distortion-free'' TFR parts (i.e.\ those lying within the corresponding cones-of-influence \ref{coidef} with $\epsilon=0.001$) are presented.

\item To simplify all expressions, the circular frequencies -- in rad/s -- are mainly used, denoted by $\omega,\nu,\xi$ (with additional subscripts or superscripts); to convert these frequencies to Hz, one needs to divide them by $2\pi$. Note, however, that the signal's sampling frequency $f_s$ is always taken to be in Hz.

\item Throughout this work we use three colorcodes (e.g.\ see Figs.\ \ref{fig:tfrex} and \ref{fig:aserr}): blue-red for WFT/WT transforms, white-black for SWFT/SWT (to make their structure, which sometimes contains very sharp lines, well visible) and white-red for the other quantities (e.g.\ different errors and alike).

\end{itemize}

\section{Errors of the analytic estimates}\label{app:asigerr}

For a single AM/FM component (\ref{amfm}) the error of the analytic estimates $A^{a}(t),\phi^{a}(t)$ (\ref{ap}) of its amplitude and phase $A(t),\phi(t)$ can be quantified by $\varepsilon^a(t)\equiv A(t)e^{i\phi(t)}-A^a(t)e^{i\phi^a(t)}$. Substituting $s(t)=(1/2)[A(t)e^{i\phi(t)}+(A(t)e^{i\phi(t)})^*]$ into (\ref{asig}), and taking into account that for real signals $\widehat{s^*}(\xi)=[\hat{s}(-\xi)]^*$, one obtains
\begin{equation}\label{aphiamp}
s^a(t)\equiv A^{a}(t)e^{i\phi^{a}(t)}=[A(t)e^{i\phi(t)}]^{+}+{\big(}[A(t)e^{i\phi(t)}]^{-}{\big)}^*,
\end{equation}
Then, subtracting from $A(t)e^{i\phi(t)}=\langle A(t)e^{i\phi(t)}\rangle+[A(t)e^{i\phi(t)}]^++[A(t)e^{i\phi(t)}]^-$ the expression for $A^{a}(t)e^{i\phi^{a}(t)}$ (\ref{aphiamp}), the error of the analytic approximation (\ref{ap}) can be brought to a form
\begin{equation}\label{aserr}
\varepsilon^{a}(t)\equiv A(t)e^{i\phi(t)}-A^{a}(t)e^{i\phi^{a}(t)}=\langle Ae^{i\phi(t)}\rangle+2i{\rm Im}\left[{\big[}A(t)e^{i\phi(t)}{\big]}^{-}\right].\\
\end{equation}
It is usually dominated by the second (purely imaginary) term, while $\langle Ae^{i\phi(t)}\rangle$ is often negligible or is exactly zero (see below). Note, that by differentiating (\ref{aserr}) one can also find the inaccuracies of the analytic estimates for the higher amplitude/phase time-derivatives, e.g.\ that for the analytic frequency $\partial_t\phi(t)-\partial_t\phi^{a}(t)=O(\partial_t\varepsilon^{a}(t),\varepsilon^{a}(t))$.

To better understand the properties of $\varepsilon^a(t)$ (\ref{aserr}), consider an AM/FM component with a simple sinusoidal law of amplitude and frequency modulation:
\begin{equation}\label{amfmsimple}
\begin{aligned}
s(t)\;=\;&(1+r_a\cos[\nu_a t+\varphi_a])\cos[\nu t+\varphi+r_b\sin(\nu_b t+\varphi_b)]\\
=\;&\frac{1}{2}\Bigg[{\Big(}1+\frac{1}{2}r_ae^{i(\nu_a t+\varphi_a)}+\frac{1}{2}r_ae^{-i(\nu_a t+\varphi_a)}{\Big)}\\
&\hphantom{\frac{1}{2}\Bigg[}\times\sum_{n=-\infty}^\infty J_n(r_b)e^{i[(\nu+n\nu_b)t+(\varphi+n\varphi_b)]}\Bigg]+c.c.,\\
\end{aligned}
\end{equation}
where the expansion $e^{ia\sin \phi}=\sum_{n=-\infty}^\infty J_n(a)e^{in\phi}$ was used, $J_n(x)=(-1)^nJ_{-n}(x)$ denote Bessel functions of the first kind, and $c.c.$ stands for the complex conjugate of the preceding expression. Note that, because the definition of phase implies $\phi'(t)\geq0$, one has a restriction $r_b\nu_b\leq\nu$ in (\ref{aaex}), whereas according to the amplitude definition $A(t)\geq0$ one has $r_a\leq1$.

Substituting (\ref{amfmsimple}) into (\ref{aserr}), one obtains
\begin{equation}\label{aaex}
\begin{gathered}
\varepsilon^{a}(t)\;=\;\frac{1}{2}\langle\tilde{\varepsilon}(t)\rangle+i{\rm Im}[\tilde{\varepsilon}(t)],\\
\begin{aligned}
\tilde{\varepsilon}(t)\;\equiv\;&\sum_{n:\nu+n\nu_b\leq0}J_n(r_b)e^{i(\varphi+n\varphi_b)}e^{i(\nu+n\nu_b)t}\\
&+\frac{r_a}{2} \sum_{n:\nu+n\nu_b+\nu_a\leq0}J_n(r_b)e^{i(\varphi+n\varphi_b+\varphi_a)}e^{i(\nu+n\nu_b+\nu_a)t}\\
&+\frac{r_a}{2} \sum_{n:\nu+n\nu_b-\nu_a\leq0}J_n(r_b)e^{i(\varphi+n\varphi_b-\varphi_a)}e^{i(\nu+n\nu_b-\nu_a)t}.
\end{aligned}
\end{gathered}
\end{equation}
From (\ref{aaex}) it can be seen, that the average $\langle A(t)e^{i\phi(t)}\rangle$ is not zero only when $\nu+n\nu_b=0$ or $\nu+n\nu_b\pm\nu_a=0$ for some $n$, i.e.\ when the frequency content of $A(t)$ intersects with the frequency content of $e^{i\phi(t)}$. Obviously, this is rather a special case, which requires particular relationships between $\nu,\nu_a,\nu_b$.

\begin{figure*}[t]
\includegraphics[width=1.0\linewidth]{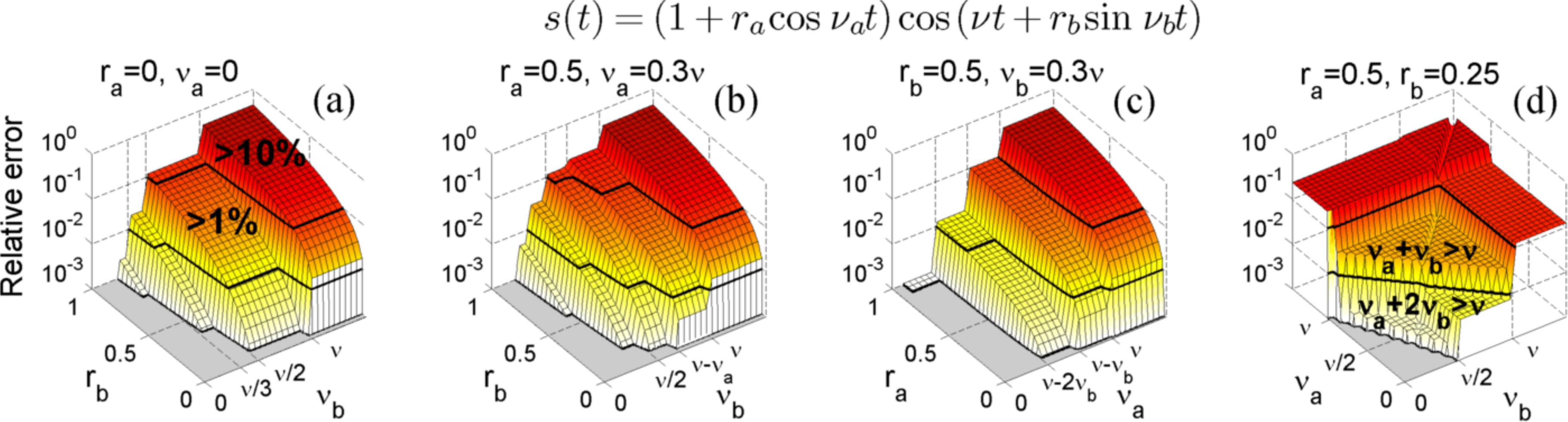}\\
\caption{Relative error of the analytic approximation $\langle |\varepsilon^a(t)|^2\rangle/\langle A^2(t)\rangle$, where $\varepsilon^a(t)$ is given by (\ref{aserr}), in dependence on signal parameters (assuming infinite time-length and sampling frequency): (a) its dependence on the parameters of frequency modulation $r_b,\nu_b$ for the FM component $s(t)=\cos(\nu t+r_b\sin\nu_b t)$; (b) its dependence on the parameters of frequency modulation $r_b,\nu_b$ for the AM/FM component $s(t)=(1+0.5\cos(0.3\nu t))\cos(\nu t+r_b\sin\nu_b t)$, with parameters of the amplitude modulation being fixed; (c) its dependence on the parameters of amplitude modulation $r_a,\nu_a$ for the AM/FM component $s(t)=(1+r_a\cos(\nu_a t))\cos(\nu t+0.5\sin0.3\nu t)$, with parameters of the frequency modulation being fixed; (d) its dependence on the frequencies of amplitude and frequency modulations $\nu_a,\nu_b$ for the AM/FM component $s(t)=(1+0.5\cos(\nu_a t))\cos(\nu t+0.25\sin\nu_b t)$, with other parameters remaining fixed. Solid black lines indicate levels of $0.001$, $0.01$ and $0.1$; for simplicity, the behavior of the relative error below $0.001$ is not shown.}
\label{fig:aserr}
\end{figure*}

As can be seen from (\ref{aaex}), the error $\varepsilon^{a}(t)$ is proportional to $r_a,r_b$, but the coefficients of this proportionality are determined by the values of $n$ for which $\nu+n\nu_b\leq0$ and/or $\nu+n\nu_b\pm\nu_a\leq0$. Hence, the behavior of $\varepsilon^{a}(t)$ can be partitioned into different regimes, separated by discontinuous step increases in error. When only amplitude, but not frequency, modulation is present ($r_b=0\Rightarrow J_{n\neq0}(r_b)=0$, so that one can set $n=0$ in all terms of (\ref{aaex})), it follows that the analytic approximation is exact ($\varepsilon^{a}(t)=0$) when $\nu_a<\nu$; otherwise the error is linearly proportional to $r_a$. When there is no amplitude modulation ($r_a=0$), but frequency modulation exists, the $\varepsilon^{a}(t)$ depends on $r_b$ and $\nu_b$ in a complicated fashion, as the relationships between $J_n(r_b)$ in (\ref{aaex}) change with $r_b$. However, when $r_b\in[0,1]$, one has $|J_{n_1}(r_b)|\gg|J_{n_2}(r_b)|$ if $|n_1|<|n_2|$. Therefore, in this case the error is largely determined by the smallest $n$ for which $\nu-n\nu_b\leq0$, being negligible for large $n$ but often considerable for $n\leq2$; it is also proportional to $r_b$, but the dependence is not simple.

The latter case is illustrated in Fig.\ \ref{fig:aserr}(a), where it can be seen that the relative error of the analytic approximation increases in steps when $\nu_b$ passes the levels $\nu/n$. Because for $r_b\leq1$ one has $J_{n>3}(r_b)<0.0025$, it becomes non-negligible only for $\nu_b\geq\nu/3$. In the case when both amplitude and frequency modulations are present the situation becomes more complicated, on account of an additional contribution from the mixing terms $\nu+n\nu_b\pm\nu_a\leq0$ in (\ref{aaex}). Thus, there appear additional step increases in error when $\nu_b$ crosses not only the levels $\nu/n$, but also $(\nu\pm\nu_a)/n$, as shown in Fig.\ \ref{fig:aserr}(b-d). Nevertheless, for $r_b\in[0,1]$, not-small $r_a$ and $\nu_a<\nu$, the error is largely determined by the smallest $n$ for which $\nu-\nu_a-n\nu_b\leq0$, as is clearly seen in Fig.\ \ref{fig:aserr}(d).

Amplitude/frequency modulations more complex than (\ref{amfmsimple}) can always be expanded in Fourier series (see (\ref{amfmft})), so that the expression for $\varepsilon^a(t)$ (\ref{aaex}) in this case will include additional terms
$$
\sim r_a^{(i)}\sum_{\{n^{(j)}\}:\nu\pm\nu_a^{(i)}+\sum_{j}n^{(j)}\nu_b^{(j)}\leq0}e^{i(\nu\pm\nu_a^{(i)}
+\sum_{j}n^{(j)}\nu_b^{(j)})t}\prod_{j}J_{n^{(j)}}(r_b^{(j)})
$$
for each Fourier term $i$ in the amplitude modulation (including the DC value $\nu_a^{(0)}=0$) and each Fourier term $j$ in the frequency modulation; the behavior of the error will therefore be partitioned into additional possible regimes.

\section{Importance of taking only positive frequency axis in the (S)WFT/(S)WT computation}\label{app:posfreq}

\begin{figure*}[t!]
\includegraphics[width=1.0\linewidth]{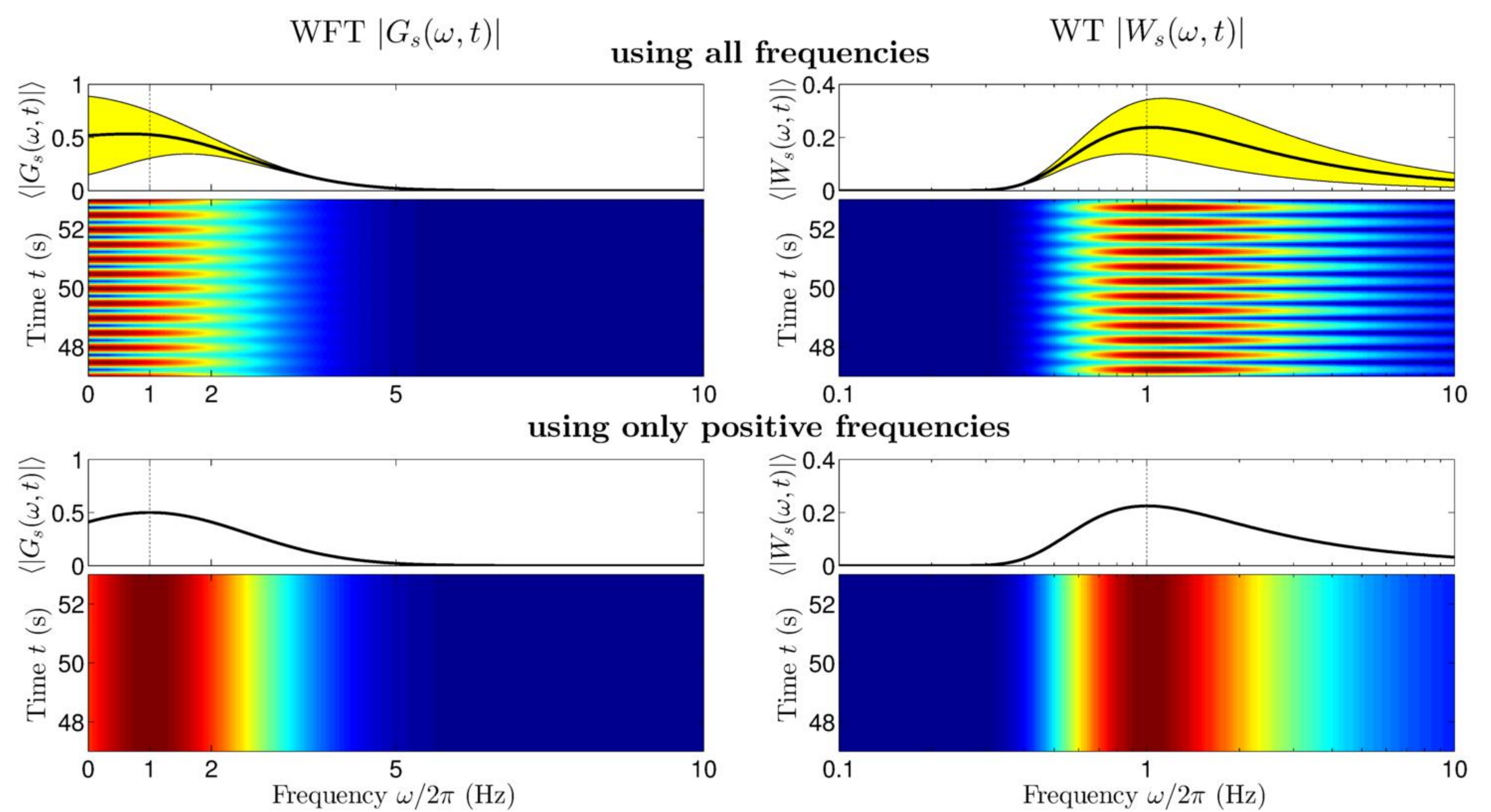}\\
\caption{WFT (left) and WT (right) amplitudes for a single-tone signal $s(t)=\cos2\pi t$ calculated using unfiltered signal (upper) and only its positive frequency part $s^{+}(t)=\frac{1}{2}e^{i2\pi t}$ (bottom). No other preprocessing, e.g.\ that described in Sec.\ \ref{sec:pract1}, was performed here. Yellow regions show mean $\pm\sqrt{2}\std()$ of the TFR amplitudes. WFT (WT) is based on a Gaussian window (Morlet wavelet) with $f_0=0.1$. SWFT and SWT are not shown, but similar situation is observed there as well.}
\label{fig:posfreq}
\end{figure*}

Often one computes the (S)WFT/(S)WT using a full signal, but here we show that it is essential to use only its positive frequency part $s^{+}(t)$ (\ref{nt}), as we do in this work (see (\ref{wft}), (\ref{wt})). This is because, otherwise, for a large spread of $\hat{g}(\xi)$ or $\hat{\psi}(\xi)$ the negative frequency part, in which one is not interested, will interfere with positive frequencies and thus corrupt both the resultant representation and the components reconstructed from it (e.g.\ by (\ref{iwft}), (\ref{iwt}), (\ref{dirrec}) or (\ref{ridgerec})). This is illustrated in Fig.\ \ref{fig:posfreq}. As can be seen, for large window/wavelet spread in frequency ($f_0=0.2$), even a simple tone's TFR is badly spoiled if one takes signal as it is; but by filtering its negative frequencies one can avoid this problem. The same obviously applies to the SWFT and SWT, because they are based on the original WFT/WT. These issues express themselves slightly differently in the (S)WFT and (S)WT, so we now consider these two representations separately.

\subsection{WFT and SWFT}

Here we consider the problems that arise if in the WFT (\ref{wft}), instead of the signal's positive frequency part $s^{+}(t)$, one uses the full signal $s(t)$, which we will assume to be the case for this section. Note, that in all formulas (e.g.\ reconstruction (\ref{iwft}), (\ref{dirrec})) one should then integrate over $\omega\in(0,\infty)$ instead of the full range. To investigate the interference effects coming from negative frequencies, consider for simplicity a single tone signal, for which such a WFT can be brought to the form:
\begin{equation}\label{negwft0}
\begin{aligned}
s(t)=&A\cos(\nu t+\varphi)=\frac{1}{2}(e^{i(\nu t+\varphi)}+e^{-i(\nu t+\varphi)})\Rightarrow\\
G_s(\omega,t)=&\frac{A}{2}[\hat{g}(\omega-\nu)e^{i(\nu t+\varphi)}+\hat{g}(\omega+\nu)e^{-i(\nu t+\varphi)}]\\
\equiv& \frac{A}{2}\hat{g}(\omega-\nu)e^{i(\nu t+\varphi)}[1+\epsilon_{neg}(\omega)e^{-i(2\nu t+2\varphi)}],\\
\end{aligned}
\end{equation}
where we have defined the interference term as
\begin{equation}\label{nerrwft}
\epsilon_{neg}(\omega)\equiv \frac{\hat{g}(\omega+\nu)}{\hat{g}(\omega-\nu)}\mbox{   }\left(=e^{-2f_0^2\nu\omega}\mbox{ for Gaussian window}\right).
\end{equation}

From (\ref{negwft0}) one can also calculate the WFT amplitude, phase and frequency as:
\begin{equation}\label{negwft}
\begin{aligned}
|G_s(\omega,t)|=&\frac{A}{2}|\hat{g}(\omega-\nu)|\sqrt{1+\epsilon_{neg}^2-2\epsilon_{neg}\cos(2\nu t+2\varphi)}\\
=& \frac{A}{2}|\hat{g}(\omega-\nu)|[1-\epsilon_{neg}\cos(2\nu t+2\phi)+O(\epsilon_{neg}^2)],\\
\operatorname{arg}[G_s(\omega,t)]=&\arctan\left(\frac{1-\epsilon_{neg}}{1+\epsilon_{neg}}\tan(\nu t+\varphi)\right)\\
=& \nu t+\varphi-\epsilon_{neg}\sin(2\nu t+2\varphi)+O(\epsilon_{neg}^2),\\
\nu_G(\omega,t)=&\frac{\nu}{\frac{1+\epsilon_{neg}}{1-\epsilon_{neg}}\cos^2(\nu t+\varphi)+\frac{1-\epsilon_{neg}}{1+\epsilon_{neg}}\sin^2(\nu t+\varphi)}\\
=& \nu[1-2\epsilon_{neg} \cos(2\nu t+2\varphi)+O(\epsilon_{neg}^2)].\\
\end{aligned}
\end{equation}

As can be seen, even for a single tone signal, the phase, amplitude and frequency estimated from the WFT (by any method) are not strictly equal to the real ones (as when taking only the positive frequency part), with the error being given by $\epsilon_{neg}(\omega)$ (\ref{nerrwft}). This error is higher for a larger spread of the window function FT $\hat{g}(\omega)$ (or lower $f_0$), and for smaller $\nu$, which is quite clear: the former determines the frequency resolution, whereas the latter is half of the distance between positive and negative frequency parts (which is obviously $2\nu$).

This leads to restrictions on the minimal frequency $\nu_{\min}$ which can be reliably represented in a WFT, i.e.\ for which the interference-related error $\epsilon_{neg}$ is small. It depends on the chosen window parameters, and $\nu_{\min}$ can be quite high e.g.\ for a Gaussian window with small $f_0$. This in turn restricts the choice of parameters, mainly by imposing an upper bound on the window time-resolution (the lower bound on $f_0$). The interference error obviously propagates through synchrosqueezing, so the same considerations apply to the SWFT. Taking only the positive parts in WFT calculation, as done in (\ref{wft}), avoids all these issues and allows one to study low frequencies that are as low as we wish, using any window parameters.

\subsection{WT and SWT}

We now consider the problems that arise if we use the full signal $s(t)$ in the WT (\ref{wft}) instead of $s^{+}(t)$. In this case, i.e.\ taking the full frequency axis in the WT calculation, we find similar issues to those discussed for the WFT, although with some differences. Performing for the WT the procedure that was used for the WFT, we obtain
\begin{equation}\label{negwt0}
\begin{aligned}
s(t)=&A\cos(\nu t+\varphi)=\frac{1}{2}(e^{i(\nu t+\varphi)}+e^{-i(\nu t+\varphi)})\Rightarrow\\
W_s(\omega,t)=& \frac{A}{2}[\hat{\psi}^*(\omega_\psi\nu/\omega)e^{i(\nu t+\varphi)}+\hat{\psi}^*(-\omega_\psi\nu/\omega)e^{-i(\nu t+\varphi)}]\\
\equiv& \frac{A}{2}\hat{\psi}^*(\omega_\psi\nu/\omega)e^{i(\nu t+\varphi)}[1+\epsilon_{neg}(\omega)e^{-i(2\nu t+2\varphi)}],\\
\end{aligned}
\end{equation}
where now the error term is
\begin{equation}\label{nerrwt}
\epsilon_{neg}(\omega)\equiv \frac{\hat{\psi}^*(-\omega_\psi\frac{\nu}{\omega})}{\hat{\psi}^*(\omega_\psi\frac{\nu}{\omega})}\mbox{   }\left(=e^{-2\pi f_0\omega_\psi\frac{\nu}{\omega}}\mbox{ for Morlet wavelet}\right).
\end{equation}

From (\ref{negwt0}) one can also calculate the WT amplitude, phase and frequency as:
\begin{equation}\label{negwt}
\begin{aligned}
|W_s(\omega,t)|=& \frac{A}{2}|\hat{\psi}(\omega_\psi\nu/\omega)|\sqrt{1+\epsilon_{neg}^2-2\epsilon_{neg}\cos(2\nu t+2\varphi)}\\
=& \frac{A}{2}|\hat{\psi}(\omega_\psi\nu/\omega)|[1-\epsilon_{neg}\cos(2\nu t+2\phi)+O(\epsilon_{neg}^2)],\\
\operatorname{arg}[G_s(\omega,t)]=& \arctan\left(\frac{1-\epsilon_{neg}}{1+\epsilon_{neg}}\tan(\nu t+\varphi)\right)\\
=& \nu t+\varphi-\epsilon_{neg}\sin(2\nu t+2\varphi)+O(\epsilon_{neg}^2),\\
\nu_W(\omega,t)=&
\frac{\nu}{\frac{1+\epsilon_{neg}}{1-\epsilon_{neg}}\cos^2(\nu t+\varphi)+\frac{1-\epsilon_{neg}}{1+\epsilon_{neg}}\sin^2(\nu t+\varphi)}\\
=& \nu[1-2\epsilon_{neg} \cos(2\nu t+2\varphi)+O(\epsilon_{neg}^2)].\\
\end{aligned}
\end{equation}

Expressions (\ref{negwt}) are almost completely the same as their analog for WFT (\ref{negwft}), so the same considerations apply: the error related to interference with negative frequencies corrupts both the resultant representation and the quality of estimates that one can obtain from it by any method. The only difference of WT is the form of the error term (\ref{nerrwt}), which now depends on wavelet parameters and the ratio $\nu/\omega$ instead of their difference. As a result, this error corrupts \emph{all} frequencies equally, in contrast to what was seen for WFT. Therefore, contrary to WFT, there is no minimal `'reliable`' frequency $\nu_{\min}$, but only strict restrictions on the wavelet parameters. Thus, for Morlet wavelet with $f_0=1/3$ one should expect $\approx 1\%$ error due to negative frequencies interference, and it rapidly grows with further decrease of $f_0$. The same, obviously, concerns SWT. Nevertheless, many wavelets (called analytic) by construction have $\hat{\psi}(\xi\leq0)=0$ (so $\epsilon_{neg}=0$), implicitly performing negative frequencies filtering and thus avoiding these issues. If this is not so, one should set this manually. Thus, in any case and for any wavelet, one should calculate WT from only positive frequency part of the signal, as is done in (\ref{wt}). This removes (S)WT corruption by negative frequencies and therefore allows to use any wavelet parameters.

\subsection{The negative side effect: representation of time events}

Although taking only the positive signal part improves the representation of the oscillatory components, it simultaneously worsens the representation of sharp time-events such as delta-peaks, which can be interpreted as being another manifestation of the time-frequency uncertainty principle. Consider $s(t)=\delta(t-t_0)$. Using the full signal in WFT/WT, we obtain
\begin{equation}\label{dpfull}
\begin{aligned}
&G_s(\omega,t)=\int \delta(u-t_0)g(u-t)e^{-i\omega u}du=g(t_0-t)e^{i\omega (t-t_0)},\\
&W_s(\omega,t)=\int \delta(u-t_0)\psi^*{\Big(}\frac{\omega(u-t)}{\omega_\psi}{\Big)}\frac{\omega du}{\omega_\psi}=\frac{\omega}{\omega_\psi}\psi^*{\Big(}\frac{\omega(t_0-t)}{\omega_\psi}{\Big)}.\\
\end{aligned}
\end{equation}
However, taking only the positive frequency part of a signal, as we do in (\ref{wft}) and (\ref{wt}), will lead to convolution with $\delta^{+}(t)\equiv\delta(t)-\frac{1}{2\pi}\int_{-\infty}^{0}e^{i\xi t}d\xi$ instead of the full $\delta$-function, so the WFT/WT become
\begin{equation}\label{dppos}
\begin{aligned}
G_s(\omega,t)=& \int \delta^{+}(u-t_0)g(u-t)e^{-i\omega u}du\\
=& g(t_0-t)e^{i\omega (t-t_0)}[1-c(\omega,t)e^{-i\omega (t-t_0)}],\\
&\mbox{where }c(\omega,t)\equiv\frac{\int_{-\infty}^{0}e^{i\xi(t-t_0)}\hat{g}(\omega-\xi)d\xi}{\int e^{i\xi(t-t_0)}\hat{g}(\xi)d\xi},\\
W_s(\omega,t)=& \int \delta^{+}(u-t_0)\psi^*{\Big(}\frac{\omega(u-t)}{\omega_\psi}{\Big)}\frac{\omega du}{\omega_\psi}\\
=& \frac{\omega}{\omega_\psi}\Big[\psi^{+}{\Big(}\frac{\omega(t_0-t)}{\omega_\psi}{\Big)}\Big]^*
=\frac{\omega}{\omega_\psi}\psi^*{\Big(}\frac{\omega(t_0-t)}{\omega_\psi}{\Big)}\big[1-c(\omega,t)\big],\\
&\mbox{where }c(\omega,t)\equiv\frac{\int_{-\infty}^{0}e^{i\xi(t-t_0)}\hat{\psi}^*{\Big(}\frac{\omega_\psi\xi}{\omega}{\Big)}d\xi}
{\int e^{i\xi(t-t_0)}\hat{\psi}^*{\Big(}\frac{\omega_\psi\xi}{\omega}{\Big)}d\xi},\\
\end{aligned}
\end{equation}
where $\psi^{+}(t)$ is the positive frequency part of the wavelet $\psi(t)$.

As can be seen, due to filtering out the signal's negative frequency part, (\ref{dppos}) contains an additional terms $c(\omega,t)$ compared to (\ref{dpfull}). These terms corrupt delta-peaks in the WFT/WT representation and thus can be viewed as errors. They propagate to the SWFT/SWT as well. One can characterize the relative error $\epsilon_{neg}(\omega)$ of the representation (\ref{dppos}) as $c(\omega,t)$ at the time of delta-peak $t=t_0$, i.e.\ \begin{equation}\label{nerrdelta}
\begin{aligned}
\mbox{\textbf{WFT:}}\;\;&\epsilon_{neg}(\omega)\equiv c(\omega,t_0)=
\frac{\int_{\omega}^{+\infty} g(\xi)d\xi}{\int g(\xi)d\xi},\\
\mbox{\textbf{WT:}}\;\;&\epsilon_{neg}\equiv c(\omega,t_0)=
\frac{\int_{-\infty}^{0} \psi(\xi)d\xi}{\int \psi(\xi)d\xi}.\\
\end{aligned}
\end{equation}
The relative errors $\epsilon_{neg}(\omega)$ (\ref{nerrdelta}), arising for sharp time events when computing WFT/WT by (\ref{wft}) and (\ref{wt}), are in full analogy with the relative errors (\ref{nerrwft}) and (\ref{nerrwt}), seen previously for the representation of tones (i.e.\ frequency events) when using the full signal for WFT/WT computation. The former and latter appear in different cases, being mutually exclusive, and both of them grow with the increasing spread of $\hat{g}(\xi)$ (for WFT) or with the wavelet part at negative frequencies $\hat{\psi}(\xi<0)$ (for WT).

In the case of the WFT, $\epsilon_{neg}(\omega)$ (\ref{nerrdelta}) become negligible for large frequencies, so that one has $|\epsilon_{neg}(\omega\geq \omega_e(\epsilon))|\leq\epsilon$, with a threshold value $\omega_e(\epsilon)$ that is proportional to the spread of $\hat{g}(\xi)$ (it can be expressed in terms of the $\epsilon$-supports defined in (\ref{eswft}) as $\omega_e(\epsilon)=\xi_2(2\epsilon)$). Hence one observes a kind of boundary effect, but in the frequency domain. For the WT, $|\epsilon_{neg}|$ is determined by the relative amount of wavelet FT located at negative frequencies $\hat{\psi}(\xi<0)$ and does not depend on the frequency $\omega$. Note that $c(\omega,t)$ (\ref{dppos}) in the case of WT measures simply the quality of the approximation $\psi(t)\approx\psi^{+}(t)$: one can view WT (\ref{wt}) as calculated using the full signal but only analytic wavelet part $\psi^{+}(t)$. Therefore, for WT all errors are in fact ``illusory'' and can be removed by using $\psi^{+}(t)$ instead of $\psi(t)$ everywhere: the only problem is that $\psi^{+}(t)$ is rarely known in explicit form.

In general, the issues mentioned will arise for any pulses, not only for delta-peaks. Thus, for a Gaussian pulse $s(t)=\frac{1}{\sqrt{2\pi}a}e^{-(t-t_0)^2/2a^2}$ the relative errors $c(\omega,t)$ are obtained simply by replacing $e^{i\xi (t-t_0)}\rightarrow e^{-a^2\xi^2/2}e^{i\xi (t-t_0)}$ under the integral in their expressions (\ref{dppos}). Nevertheless, it is for delta-peaks that the errors are most pronounced.

Summarizing, taking only the positive frequency part of the signal in WFT/WT computation (as we do in (\ref{wft}) and (\ref{wt})) removes the errors in the representation of AM/FM components, but introduces similar errors into representation of the sharp time events. However, TFRs are mainly devoted to representing the signal's oscillatory dynamics, so we favor the accurate representation of tones rather than delta-peaks, thus choosing to filter out the negative frequency part of the signal when computing its TFRs.

\section{Derivation of the direct reconstruction formulas}\label{app:derdirrec}

Using the frequency domain forms of the WFT (\ref{wft}) and WT (\ref{wt}) one can show that
\begin{equation}\label{Dfwft}
\begin{aligned}
&\int \omega^n G_s(\omega,t)d\omega
=\int \omega^n d\omega \frac{1}{2\pi}\int_0^{\infty}\hat{s}(\xi)\hat{g}(\omega-\xi)e^{i\xi t}d\xi\\
&=\big/\omega=\wt{\omega}+\xi,\;(\wt{\omega}+\xi)^n=\sum_{k=0}^nC_n^k\wt{\omega}^{n-k}\xi^k\big/\\
&=\sum_{k=0}^nC_k^n\Big(\int\wt{\omega}^{n-k}\hat{g}(\wt{\omega})d\wt{\omega}\Big)
\Big(\frac{1}{2\pi}\int_0^\infty\xi^k\hat{s}(\xi)e^{i\xi t}d\xi\Big)\\
&=\sum_{k=0}^nC_k^n\Big(\frac{1}{2}\int\wt{\omega}^{n-k}\hat{g}(\wt{\omega})d\wt{\omega}\Big)
(-i)^k\partial_t^ks^a(t)\\
\end{aligned}
\end{equation}
\begin{equation}\label{Dfwt}
\begin{aligned}
&\int_0^\infty \omega W_s(\omega,t)\frac{d\omega}{\omega}=\int_0^\infty \omega^{n-1}d\omega
\frac{1}{2\pi}\int_0^\infty\hat{s}(\xi)\hat{\psi}^*(\omega_\psi\xi/\omega)d\xi\\
&=\big/\omega=\omega_\psi\xi/\wt{\omega}\mbox{, after which the integrals decouple}\big/\\
&=\Big( \omega_\psi^n\int_0^\infty\hat{\psi}^*(\wt{\omega})\frac{d\wt{\omega}}{\wt{\omega}^{n+1}} \Big)
\Big( \frac{1}{2\pi}\int_0^\infty \xi^n\hat{s}(\xi)e^{i\xi t}d\xi \Big)\\
&=\Big( \frac{\omega_\psi^n}{2}\int_0^\infty\hat{\psi}^*(\wt{\omega})\frac{d\wt{\omega}}{\wt{\omega}^{n+1}} \Big)(-i)^n\partial_t^ns^{a}(t),
\end{aligned}
\end{equation}
where it was taken into account that $\int_0^\infty\xi^k\hat{s}(\xi)e^{i\xi t}d\xi$ $=$ \\ $(-i\partial_t)^k\int_0^\infty\hat{s}(\xi)e^{i\xi t}d\xi$ $=$ $(-i)^k\partial_t^ks^{a}(t)/2$, and $C_k^n\equiv \frac{n!}{k!(n-k)!}$ are the binomial coefficients.

In a similar manner, using the time domain forms of the WFT (\ref{wft}) and WT (\ref{wt}) it can be shown that
\begin{equation}\label{Dtwft}
\begin{aligned}
&\int t^n G_s(\omega,t)e^{-i\omega t}dt
=\int t^n dt \int s^+(\tau)g(\tau-t)e^{-i\omega\tau}d\tau\\
&=\big/t=\tau-\tilde{t},\;(\tau-\tilde{t})^n=\sum_{k=0}^nC_n^k(-\tilde{t})^{n-k}\tau^k\big/\\
&=\sum_{k=0}^nC_k^n\Big(\int(-\tilde{t})^{n-k}g(\tilde{t})d\tilde{t}\Big)
\Big(\int\tau^ks^+(\tau)e^{-i\omega \tau}d\tau\Big)\\
&=
\left\{\begin{array}{l}
\underset{k=0}{\overset{n}{\sum}}(-1)^nC_k^n\Big(\int\tilde{t}^{n-k}g(\tilde{t})d\tilde{t}\Big)
(-i)^k\partial_\omega^k\hat{s}(\omega)\;\mbox{ if }\omega>0\\
0\mbox{ if }\omega<0\\
\end{array}\right.
\end{aligned}
\end{equation}
\begin{equation}\label{Dtwt}
\begin{aligned}
&\int t^n W_s(\omega,t)e^{-i\omega t}dt=
\int t^n dt \int s^+(\tau)\psi^*(\omega(\tau-t)/\omega_\psi)\frac{\omega d\tau}{\omega_\psi}\\
&=\big/t=\tau-\omega_\psi\tilde{t}/\omega,\;(\tau-\omega_\psi\tilde{t}/\omega)^n
=\sum_{k=0}^nC_n^k(-\omega_\psi\tilde{t}/\omega)^{n-k}\tau^k\big/\\
&=\sum_{k=0}^nC_k^n\Big(\frac{\omega_\psi^{n-k}}{\omega^{n-k}}\int(-\tilde{t})^{n-k}
\psi^*(\tilde{t})e^{i\omega_\psi\tilde{t}}d\tilde{t}\Big)
\Big(\int\tau^ks^+(\tau)e^{-i\omega \tau}d\tau\Big)\\
&=
\left\{\begin{array}{l}
\underset{k=0}{\overset{n}{\sum}}(-1)^nC_k^n\Big(\int\tilde{t}^{n-k}\psi^*(\tilde{t})e^{i\omega_\psi\tilde{t}}d\tilde{t}\Big)
(-i)^k\partial_\omega^k\hat{s}(\omega)\;\mbox{ if }\omega>0\\
0\mbox{ if }\omega<0\\
\end{array}\right.
\end{aligned}
\end{equation}
where it was taken into account that $\int \tau^ks^+(\tau)e^{-i\omega\tau}d\tau$ $=$ \\
$(-i\partial_\omega)^k\int s^+(\tau)e^{-i\omega\tau}d\tau$ $=$ $(-i)^k\partial_\omega^k\hat{s}(\omega)$ if $\omega>0$ and $=0$ otherwise (as only the positive frequency part of the signal $s^+(t)$ is used).

The expressions (\ref{Dfwft}-\ref{Dtwt}) provide a way to reconstruct any order derivatives of the signal's representation in time and frequency from its WFT/WT. Thus, the basic reconstruction formulas -- (\ref{iwft}) for the WFT and (\ref{iwt}) for the WT -- follow directly from (\ref{Dfwft}),(\ref{Dtwft}) and (\ref{Dfwt}),(\ref{Dtwt}) with $n=0$, respectively.

The formulas for the direct frequency estimation (\ref{dirrec}) follow from (\ref{Dfwft}) and (\ref{Dfwt}) with $n=1$, that can be rewritten as
\begin{gather}
\int \omega G_s(\omega,t)d\omega=C_gs^a(t)\Big[\overline{\omega}_g-i\frac{\partial_ts^a(t)}{s^a(t)}\Big],\;
\overline{\omega}_g\equiv \frac{\int \omega\hat{g}(\omega)d\omega}{2C_g}, \label{D1wft}\\
\int \omega W_s(\omega,t)\frac{d\omega}{\omega}=D_\psi(-i\partial_ts^a(t)),\;
D_\psi\equiv\frac{\omega_\psi}{2}\int_0^\infty \hat{\psi}^*(\omega)\frac{d\omega}{\omega^2}, \label{D1wt}
\end{gather}
where $C_{g,\psi}$ are as defined in (\ref{iwft}),(\ref{iwt}). For a single AM/FM component (\ref{amfm}) one has $s^a(t)=A^a(t)e^{i\phi^a(t)}$ and $\partial_ts^a(t)=s^a(t)\big[\frac{\partial_tA^a(t)}{A^a(t)}+i\partial_t\phi(t)\big]$, so that $\frac{\partial_tA^{a}(t)}{A^{a}(t)}={\rm Re}\frac{\partial_ts^{a}(t)}{s^{a}(t)}$, $\partial_t\phi^{a}(t)={\rm Im}\frac{\partial_t s^{a}(t)}{s^{a}(t)}$. Therefore, dividing both sides of (\ref{D1wft}) by $\int G_s(\omega,t)d\omega=C_gs^{a}(t)$ (\ref{iwft}) and separating real and imaginary parts of the result, one obtains
\begin{equation}\label{DDwft}
\begin{aligned}
\nu^{a}(t)\equiv\partial_t\phi^{a}(t)=&{\rm Re}\Bigg[\frac{\int \omega G_s(\omega,t)d\omega}{\int G_s(\omega,t)d\omega}-\overline{\omega}_g\Bigg],\\
\frac{\partial_t A^{a}(t)}{A^{a}(t)}=-&{\rm Im}\Bigg[\frac{\int \omega G_s(\omega,t)d\omega}{\int G_s(\omega,t)d\omega}-\overline{\omega}_g\Bigg].\\
\end{aligned}
\end{equation}
If the analytic amplitude and phase are the same as the actual ones (\ref{ap}), as is assumed in this work, the first of (\ref{DDwft}) coincides with the WFT-based direct frequency estimate (\ref{dirrec}), while the second one provides a way to estimate the time-derivative of the amplitude. The WT-based counterpart of (\ref{DDwft}) is derived in the same way, i.e.\ dividing both sides of (\ref{D1wt}) by $\int_0^\infty W_s(\omega,t)\frac{d\omega}{\omega}=C_\psi s^{a}(t)$ (\ref{iwt}) and separating the real and imaginary parts of the result, which gives
\begin{equation}\label{DDwt}
\begin{aligned}
\nu^{a}(t)\equiv\partial_t\phi^{a}(t)=&{\rm Re}\Bigg[\frac{D_\psi^{-1}\int \omega W_s(\omega,t)\frac{d\omega}{\omega}}{C_\psi^{-1}\int W_s(\omega,t)\frac{d\omega}{\omega}}\Bigg],\\
\frac{\partial_tA^{a}(t)}{A^{a}(t)}=-&{\rm Im}\Bigg[\frac{D_\psi^{-1}\int \omega W_s(\omega,t)\frac{d\omega}{\omega}}{C_\psi^{-1}\int W_s(\omega,t)\frac{d\omega}{\omega}}\Bigg].\\
\end{aligned}
\end{equation}

Obviously, using (\ref{Dfwft}) and (\ref{Dfwt}) with an appropriate $n$, one can derive the reconstruction formulas for higher-order derivatives $\partial_t^nA(t)$ and $\partial_t^n\phi(t)$ as well. The maximum order of $n$ that one can go to is determined by the condition $\int \omega^n\hat{g}(\omega)d\omega<\infty$ for the WFT, and by $\int_0^\infty\hat{\psi}^*(\omega)\frac{d\omega}{\omega^{n+1}}<\infty$ for the WT. For the WFT with a Gaussian window (\ref{gw}) and for the WT with a lognormal wavelet (\ref{lw}) this condition is satisfied for all $n\geq0$, while for the WT with a Morlet wavelet (\ref{mw}) one has $n<1$ (so that even the instantaneous frequency cannot be estimated by (\ref{dirrec}), and one needs to use hybrid reconstruction (\ref{hybridrec}) instead).

\section{Forecasting model for predictive padding}\label{app:forecasting}

The predictive padding strategy, introduced in Sec.\ \ref{sec:pract3}, aims to eliminate boundary effects in the signal's TFR by complementing the signal with its inferred/forecast past/future values. This appendix discusses the scheme used to predict the signal's behavior beyond its time limits.

Given a signal $s(t_n=(n-1)\Delta t),\;n=1,..,N$, consider first its extension for $t>T$, with $T=(N-1)\Delta t$ denoting its time duration. In the context of time-frequency analysis, it seems most appropriate to forecast the signal based on its spectral content. The simplest way of doing this is to represent signal as a sum of tones and continue this behavior to $t>T$. Thus, the signal is modelled as
\begin{equation}\label{sinmodel}
\begin{aligned}
s(t_n)=& x(t_n)+\sigma\zeta(t_n)=c_0+\sum_{m=1}^M c_m\cos(\omega_m t_n+\varphi_m)+\sigma\zeta(t_n)\\
=&a_0+\sum_{m=1}^M[a_m\cos\omega_m t_n+b_m\sin\omega_m t_n]+\sigma\zeta(t_n),
\end{aligned}
\end{equation}
where $M$ denotes the chosen model order, $\zeta(t_n)$ is Gaussian white noise of unit variance, and $x(t_n)$ stands for the noise-free signal, given by the sum of sinusoids. Having found appropriate $c_m,\varphi_m,\omega_m$, the signal is then padded for $t>T$ with values of $x(t)$.

However, it is not easy to find the parameters in (\ref{sinmodel}). It might seem at first glance that an approximation to $x(t)$ can readily be obtained from the signal's discrete FT $\hat{s}(\xi_n)$ as $x(t)=\sum_{m=1}^{N}\hat{s}(\xi_m)e^{i\xi_m t}$, but this is not so: the discrete FT represents a periodic spectrum estimate, being an exact FT only for signals that repeat themselves with period $T$. As a result, predicting the signal based on its discrete FT is equivalent to the usual periodic continuation.


There are many methods devoted to fitting the signal with the sinusoidal model \cite{Stoica:05}, using which the parameters $c_m,\varphi_m,\omega_m$ in (\ref{sinmodel}) can be estimated reliably (i.e.\ without the periodicity constraint). However, these methods, although very accurate, are usually quite expensive computationally, as well as giving rise to a variety of issues when applied in practice. Therefore, it seems reasonable to use the more convenient and computationally cheaper procedure of estimating tone frequencies by simple iterative fitting of the signal with the corresponding sinusoids, as described below.

Using least squares fitting, the residual error of the signal's fit with $q_0+q_1\cos\omega t_n+q_2\sin\omega t_n$ is first minimized over $\omega$; then the value of $\omega$ for which the local minimum occurs and the associated ``best-fit'' parameters $q_{1,2}$ are taken as one of the $\omega_m,a_m,b_m$ in (\ref{sinmodel}); these two steps are repeated $M-1$ times using in place of signal $s(t_n)\equiv s^{(1)}(t_n)$ the residual signal $s^{(m)}(t_n)$, obtained by subtracting from $s(t_n)$ the tones found at $m-1$ previous steps. To provide more accurate and faster optimization, at each iteration it is recommended to start the search for the optimal $\omega_m$ from the frequency of the highest peak in the residual signal's discrete FT $\hat{s}^{(m)}(\xi_n)$.

An important problem of the outlined approach is that the model (\ref{sinmodel}) is stationary, and so it cannot adequately describe nonstationary signals, where spectral content changes in time (e.g.\ when the tones persist only during some time intervals, as in Fig. \ref{fig:padding}(a,b,c,d)). Therefore, for adequate prediction the model (\ref{sinmodel}) should reflect mainly a ``local'' spectrum near the signal's end ($t$ around $T$). This can be achieved by using a weighted least squares procedure, with more weight being concentrated near the corresponding time boundary. Then the full procedure of estimating the ``local'' parameters of (\ref{sinmodel}) can be summarized as
\begin{equation}\label{fcast}
\begin{aligned}
&\begin{aligned}
\left[\begin{array}{l}
q_0^{(m)}(\omega)\\
q_1^{(m)}(\omega)\\
q_2^{(m)}(\omega)\\
\end{array}\right]
=\underset{[\tilde{q}_0,\tilde{q}_1,\tilde{q}_2]}{\rm argmin}\Bigg[
&\frac{1}{N}\sum_{n=1}^N{\rm w}(t_n)\Big(s^{(m)}(t_n)-\tilde{q}_0\\
&-\tilde{q}_1\cos\omega t_n-\tilde{q}_2\sin\omega t_n\Big)^2\Bigg],
\end{aligned}\\
&\begin{aligned}
\rho^{(m)}(\omega)\equiv \frac{1}{N}\sum_{n=1}^N{\rm w}(t_n)\Big[
&s^{(m)}(t_n)-q_0^{(m)}(\omega)\\
&-q_1^{(m)}(\omega)\cos\omega t_n-q_2^{(m)}(\omega)\sin\omega t_n\Big]^2,
\end{aligned}\\
&\omega_m= \underset{\omega}{\rm argmin}\big[\rho^{(m)}(\omega)\big],\\
&\begin{aligned}
s^{(m+1)}(t_n)\equiv
&s^{(m)}(t_n)-q_0^{(m)}(\omega_{m})\\
&-q_1^{(m)}(\omega_{m})\cos\omega_{m}t_n-q_2^{(m)}(\omega_{m})\sin\omega_{m}t_n,
\end{aligned}
\end{aligned}
\end{equation}
where one begins with $s^{(1)}(t_n)\equiv s(t_n)$, while ${\rm w}(t_n)$ denotes the chosen weighting function, and the ``best-fit'' values of $q_{0,1,2}^{(m)}(\omega)$ are found by weighted least squares for each $\omega$. The search for a minimum $\rho^{(m)}(\omega)$ is started from the frequency of the maximum in the discrete FT of $\sqrt{w(t_n)}s^{(m)}(t_n)$ and, by default, the $\omega_m$ for which such a minimum occurs is determined with accuracy $0.01\times2\pi/T$. The other parameters of (\ref{sinmodel}) are then $a_m=q_1^{(m)}(\omega_m)$, $b_m=q_2^{(m)}(\omega_m)$ and $a_0=\sum_{m=1}^M q_0^{(m)}(\omega_m)$.

Since the maximum amount of nonstationarity (i.e.\ the ``quickness'' of the spectrum's changes in time) which can be represented in TFR reliably is determined by the window/wavelet parameters (see Part II), the choice of the weighting function should be based on the window $g(t)$ or wavelet $\psi(t)$ used. Therefore, ${\rm w(t_n)}$ is selected as an exponential whose rate of decay is determined from the window/wavelet $0.5$-support (\ref{eswft}), (\ref{eswt}) in time:
\begin{equation}\label{weights}
\begin{aligned}
\mbox{\textbf{WFT:}}\;&{\rm w}(t_n)=\exp\Big[-\frac{(T-t_n)\log 2}{\tau_2(0.5)-\tau_1(0.5)}\Big],\\
\mbox{\textbf{ WT:}}\;&{\rm w}(t_n)=\exp\Big[-\frac{\omega_{\min}}{\omega_\psi}\frac{(T-t_n)\log 2}{\tau_2(0.5)-\tau_1(0.5)}\Big].\\
\end{aligned}
\end{equation}
Note, that the time resolution of the WT varies with frequency, so that the optimal ${\rm w}(t_n)$ will generally depend on frequency as well. Because using different weighting functions for each $\omega$ would be computationally very expensive, the optimal weights for the WT in (\ref{weights}) are taken as those for the minimum frequency $\omega_{\min}$, where the wavelet is most spread in time. Such a choice, however, might not be fully appropriate for higher frequencies, so the predictive padding is better defined for the WFT, whose time resolution is fixed (see Sec.\ \ref{sec:dwftwt}).

It remains to choose an appropriate order $M$ of the sinusoidal model (\ref{sinmodel}). This can be done using the Bayesian (Schwarz) information criterion (BIC) \cite{Schwarz:78}, which was empirically found to be superior to the (corrected) Akaike's information criterion \cite{Akaike:74,Hurvich:89,Burnham:02} in the present context. Thus, the number $M$ of sinusoidal components in (\ref{sinmodel}) can be selected by minimizing the functional
\begin{equation}\label{bic}
BIC=N\log\big[2\pi\rho^{(M)}(\omega_M)\big]+N+(3M+1)\log N,
\end{equation}
where $\rho^{(M)}(\omega_M)$ is given in (\ref{fcast}), while $3M+1$ is the number of parameters in the model (one $a_0$ plus $M$ of each $a_m,b_m,\omega_m$ in (\ref{sinmodel})). The order $M$ can be estimated ``on the fly'' while performing iterations (\ref{fcast}): at each such iteration one calculates the BIC (\ref{bic}), and if two of its consecutive values are higher than the current minimum (which in this case usually coincides with the global minimum), then the procedure is stopped.

The above considerations refer to forecasting the signal for $t>T$. To predict it for $t<0$, one can use the same procedure on the reversed version of the signal: $s(t)\rightarrow s(T-t)$. Finally, it should be noted that the approach outlined does not aim to provide the best or most rigorous predictive scheme, and other forecasting models can be employed instead. Nevertheless, the procedure discussed works well for the majority of signals, being at the same time very straightforward and computationally cheap. It is therefore a convenient choice given that the predictive padding represents a useful addition rather than an essential part of time-frequency analysis.

\begin{remark}
For sophisticated signals the order $M$ selected by BIC (\ref{bic}) can be quite high, so one should restrict its maximum value to ensure that the cost of the estimation (\ref{sinmodel}), $O(MN)$, is not much higher than the computational cost of the TFR calculation itself, which is $O(N_f\tilde{N}\log\tilde{N})$ (with $N_f$ being the number of frequencies at which TFR is calculated, and $\tilde{N}$ denoting the length of the padded signal). To ensure that this is the case, it is reasonable to limit $M$ to the maximum possible number of peaks in the TFR amplitude at each time, which is $N_f/2$. Additionally, since the $N$ data points composing the signal are modeled using $3M+1$ parameters in (\ref{sinmodel}), the value of $M$ cannot exceed $(N-1)/3$. Therefore, the maximum order can be chosen as $M_{\max}=\min[N_f/2,(N-1)/3]$.
\end{remark}

\section{Window and wavelet functions and their properties}\label{app:winwav}

This Appendix discusses the properties of the most common window and wavelet functions. All of them are implemented in the codes \cite{freecodes} under the specific names, which are provided in the tables below, though one can use any window/wavelet function by specifying either its frequency domain or time domain form (or both, if known).

\subsection{Window functions}

\begin{center}
\begin{figure*}[t!]
\centering\includegraphics[width=0.9\linewidth]{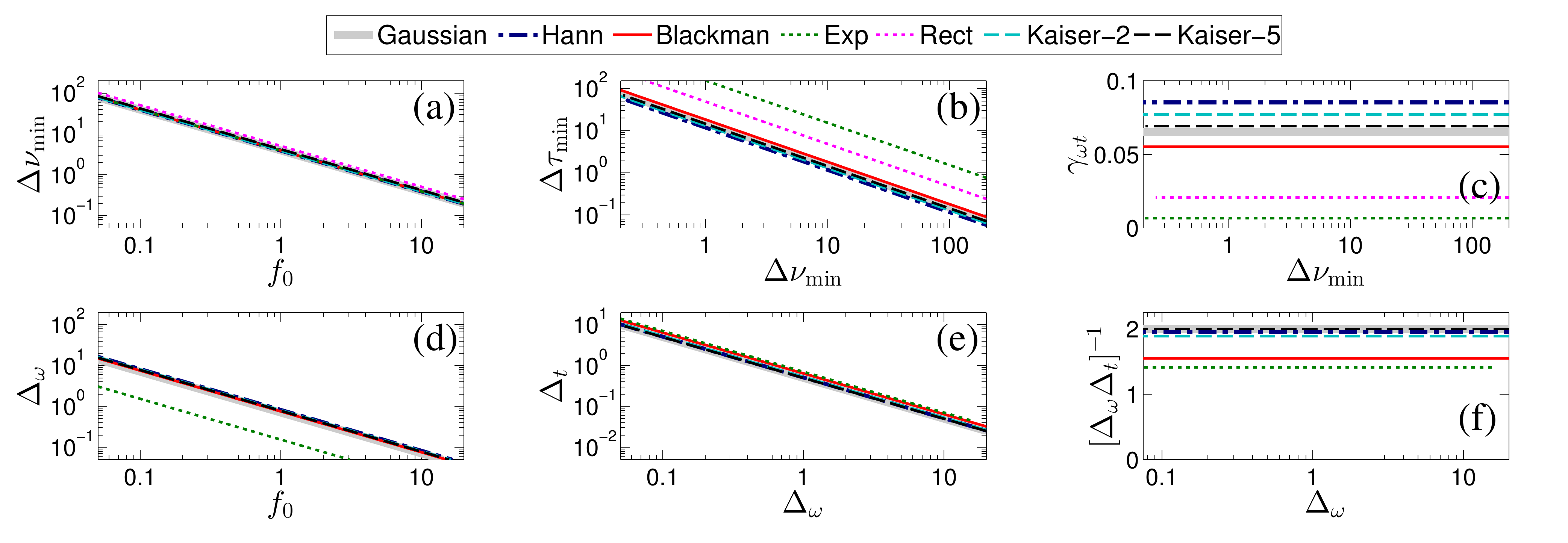}\\
\caption{(a-c): Resolution characteristics of the windows listed in Table \ref{tab:windows}, according to (\ref{tfreswft}) with $\epsilon_r=0.05$. The minimum resolvable frequency $\Delta\nu_{\min}$ in dependence on the window resolution parameter $f_0$ is shown in (a), while the dependences of the minimum resolvable time lag $\Delta\tau_{\min}$ and the joint time-frequency resolution $\gamma_{\omega t}$ on $\Delta\nu_{\min}$ are shown in (b) and (c), respectively. (d-f): The same as (a-c), but for the ``classic'' resolution characteristics (\ref{clreswft}), though they are not fully appropriate (see Sec.\ \ref{sec:tfresB}); there are no lines for the rectangular window here, as it is characterized by $\Delta_\omega=\infty$.}
\label{fig:wftres}
\end{figure*}
\end{center}

Table \ref{tab:windows} lists the commonest window forms and related quantities. The resolution characteristics (\ref{tfreswft}) for each window are presented in Fig.\ \ref{fig:wftres}(a-c). For completeness, the ``classic'' resolution characteristics (\ref{clreswft}) are also shown in Fig.\ \ref{fig:wftres}(d-f), although they do not have a straightforward relationship to the resolution of two components in the WFT, as discussed in Sec.\ \ref{sec:tfresB}.

\begin{center}
\begin{table*}[t!]
{\fontsize{9pt}{1.2\baselineskip}\selectfont
  \begin{tabular}{| >{\centering\arraybackslash}p{1.75cm} | >{\centering\arraybackslash}p{2.25cm} | p{13cm} |}
  \hline
Name & Name in codes & Description and characteristics \\ \hline
Gaussian & ``Gaussian'' (default) &
$g(t)=\frac{1}{\sqrt{2\pi}f_0}e^{-(t/f_0)^2/2},\;t\in(-\infty,\infty)$,\newline
$\hat{g}(\xi)=e^{-(f_0\xi)^2/2},\;\xi\in(-\infty,\infty)$,\newline
$R_g(\omega)=\frac{1}{2}\big[{\rm erf}\big(f_0^{-1}\omega/\sqrt{2}\big)+1\big],\;\xi_{1,2}(\epsilon)=\mp f_0^{-1}n_G(\epsilon)$,\newline
$P_g(\tau)=\frac{1}{2}\big[{\rm erf}\big(f_0\tau/\sqrt{2}\big)+1\big],\;\tau_{1,2}(\epsilon)=\mp f_0n_G(\epsilon)$,\newline
$C_g=\sqrt{\frac{\pi}{2}}f_0^{-1},\;\overline{\omega}_g=0$.
\\ \hline
Hann & ``Hann'' &
$q=4.4f_0$,\newline
$g(t)=(1+\cos(2\pi t/q))/2=\sin^2(\pi t/q+\pi/2),\;t\in(-q/2,q/2)$,\newline
$\hat{g}(\xi)=\frac{-(4\pi^2/q^2)\sin(\xi q/2)}{\xi(\xi^2-4\pi^2/q^2)},\;\xi\in(-\infty,\infty)$,\newline
$P_g(\tau)=2\tau/q-(1/2\pi)\sin(2\pi \tau/q)$,\newline
$C_g=\pi,\;\overline{\omega}_g=0$.
\\ \hline
Blackman & ``Blackman'' &
$q=5.6f_0,\;\alpha=0.16$,\newline
$g(t)=(1+\cos(2\pi t/q))/2-\alpha(1+\cos(4\pi t/q))/2,\;t\in(-q/2,q/2)$,\newline
$\hat{g}(\xi)=\frac{-(4\pi^2/q^2)\sin(\xi q/2)}{\xi}\left[\frac{1}{\xi^2-4\pi^2/q^2}-\frac{4\alpha}{\xi^2-16\pi^2/q^2}\right],\; \xi\in(-\infty,\infty)$,\newline
$P_g(\tau)=2\tau/q-[(1/2\pi)\sin(2\pi \tau/q)-(\alpha/4\pi)\sin(4\pi \tau/q)]/(1-\alpha)$,\newline
$C_g=\pi(1-\alpha),\;\overline{\omega}_g=0$.
\\ \hline
Exponential & ``Exp'' &
$q=6.5f_0$,\newline
$g(t)=e^{-|t|/q},\;t\in(-\infty,\infty)$,\newline
$\hat{g}(\xi)=2q^{-1}[\xi^2+q^{-2}]^{-1},\; \xi\in(-\infty,\infty)$,\newline
$R_g(\omega)=1/2+\pi^{-1}\arctan (q\omega),\;\xi_{1,2}(\epsilon)=\mp q^{-1}{\rm ctg}(\pi \epsilon/2)$,\newline
$P_g(\tau)=1/2+({\rm sign(\tau)}/2)[1-e^{-q^{-1}|\tau|}],\;\tau_{1,2}(\epsilon)=\pm q\log\epsilon$,\newline
$C_g=\pi,\;\overline{\omega}_g=0$.
\\ \hline
Rectangular & ``Rect'' &
$q=10f_0$,\newline
$g(t)=1,\;t\in[-q/2,q/2]$,\newline
$\hat{g}(\xi)=2\frac{\sin(q\xi/2)}{\xi},\;\xi\in(-\infty,\infty)$,\newline
$P_g(\tau)=\tau/q+0.5,\;\tau_{1,2}(\epsilon)=\pm q(1-\epsilon)/2$,\newline
$C_g=\pi,\;\overline{\omega}_g=0$.
\\ \hline
Kaiser & ``Kaiser-a'' (e.g.\ ``Kaiser-2.5'') &
$q=3\sqrt{1+|a-1/a|}f_0$,\newline
$g(t)=I_0(\pi a\sqrt{1-(2t/q)^2})/I_0(\pi a),\;t\in(-q/2,q/2)$,\newline
$C_g=\pi,\;\overline{\omega}_g=0$.
\\ \hline
  \end{tabular}}
\caption{Different window types and their characteristics (if known in analytic form). The names under which these windows are implemented in the codes \cite{freecodes} are given in the second column, but one can specify any window function there. The resolution parameter $f_0$ for each window is adjusted in such a way that, for the same $f_0$, all of them have similar frequency resolutions (as defined in (\ref{tfreswft}) with $\epsilon_r=0.05$), see Fig.\ \ref{fig:wftres}(a) below. Note, however, that it is hard to ensure this for all $a$ in the Kaiser window, so some deviations are possible in that case. An explanation of the mathematical notation is given in \ref{app:nomenclature}.
}
\label{tab:windows}
\end{table*}
\end{center}

Defined in any way, the time and frequency resolutions of all windows have simple linear proportionality to $f_0^{-1}$ and $f_0$, respectively (see Fig.\ \ref{fig:wftres}); this is an obvious result given that in all cases $g[f_0](t)=g[1](t/f_0)$. The joint time-frequency resolutions, which remain fixed, provide more useful information. Thus, it is clear that, for all conventions, exponential and rectangular windows have much worse resolution properties compared to other windows considered. Surprisingly, in terms of time-frequency resolution the Kaiser and Hann windows outperform slightly the Gaussian window (Fig.\ \ref{fig:wftres}(c)), which is widely believed to have the best resolution properties; these two windows are also very close to Gaussian even in terms of the (not fully appropriate) ``classic'' $\gamma_{\omega t}$ (Fig.\ \ref{fig:wftres}(f)), that is maximized for the latter \cite{Mallat:08,Kaiser:94,Addison:10}.

However, the Kaiser and Hann windows, apart from worse analytical tractability, have one very significant drawback as compared to the Gaussian window: they are not unimodal in the frequency domain, which makes the time-frequency support (see Sec.\ \ref{sec:tfsrec}) ill-defined and the reconstruction of components from the WFT more problematic. Because the corresponding difference in the time-frequency resolution is not huge (Fig.\ \ref{fig:wftres}(c)), the Gaussian window therefore remains a preferred choice and is implemented as default in the codes.

\subsection{Wavelet functions}

\begin{figure*}[t!]
\centering\includegraphics[width=0.9\linewidth]{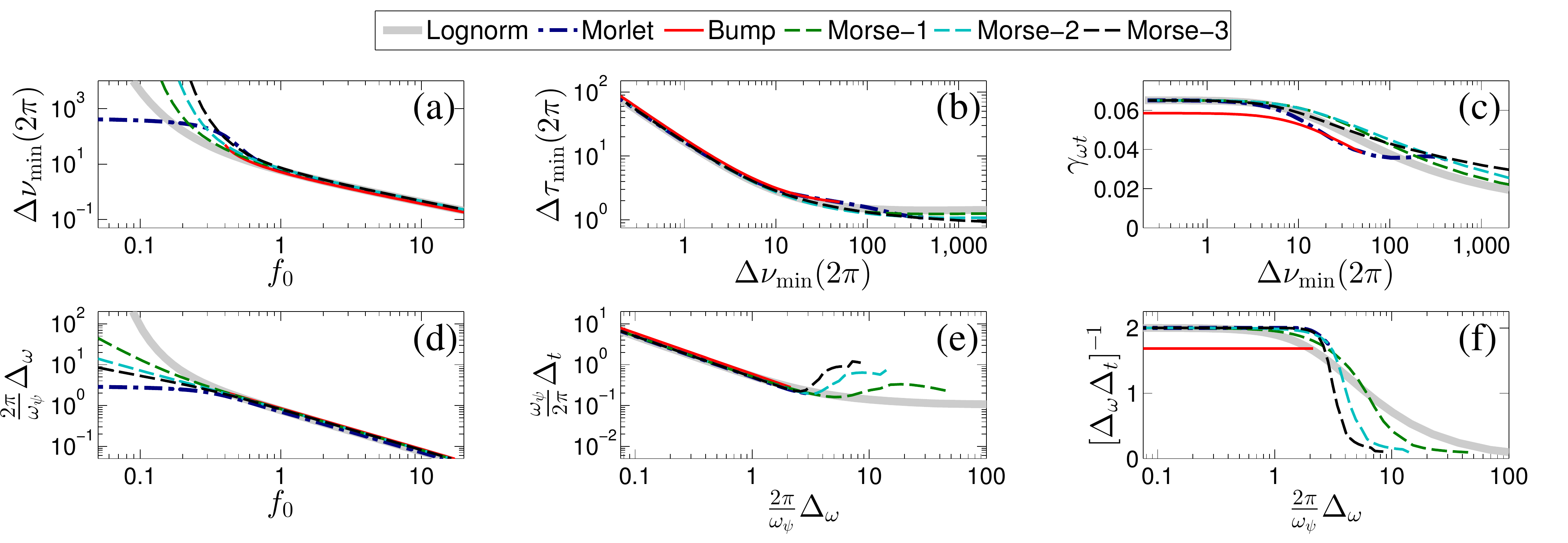}\\
\caption{(a-c): Resolution characteristics of the wavelets listed in Table \ref{tab:wavelets}, according to (\ref{tfreswt}) with $\epsilon_r=0.05$. The minimal resolvable frequency $\Delta\nu_{\min}(\omega)$ at $\omega=2\pi$ in dependence on the wavelet resolution parameter $f_0$ is shown in (a), while the dependences of the minimum resolvable time lag $\Delta\tau_{\min}(2\pi)$ and the joint time-frequency resolution $\gamma_{\omega t}$ on $\Delta\nu_{\min}(2\pi)$ are shown in (b) and (c), respectively. (d-f): The same as (a-c), but for the ``classic'' resolution characteristics (\ref{clreswt}). The latter, however, are completely inappropriate, see Sec.\ \ref{sec:tfresB}.}
\label{fig:wtres}
\end{figure*}

Table \ref{tab:wavelets} lists the commonest wavelet forms and related quantities. Note that many wavelets are included within the generalized Morse family, which was introduced in \cite{Olhede:02} and studied in detail in \cite{Lilly:09,Lilly:12}. Thus, this family includes Cauchy ($a=1$), the derivative of Gaussian ($a=2$), Airy ($a=3$) and other wavelets \cite{Lilly:12}. Different resolution characteristics for each wavelet are presented in Fig.\ \ref{fig:wtres}(a-c). For completeness, the ``classic'' resolution characteristics (\ref{clreswt}) are also shown in Fig.\ \ref{fig:wtres}(d-f) but, as discussed in Sec.\ \ref{sec:tfresB}, they are actually irrelevant and therefore will not be considered in what follows.

\begin{center}
\begin{table*}[t!]
{\fontsize{9pt}{1.2\baselineskip}\selectfont
  \begin{tabular}{| >{\centering\arraybackslash}p{1.75cm} | >{\centering\arraybackslash}p{2.25cm} | p{13cm} |}
  \hline
Name & Name in codes & Description and characteristics \\ \hline
Lognormal & ``Lognorm'' (default) &
$\hat{\psi}(\xi)=e^{-(2\pi f_0\log\xi)^2/2},\;\xi\in(0,\infty)$,\newline
$R_\psi(\omega)=\frac{1}{2}\left[{\rm erf}\Big((2\pi f_0)^{-1}\log\omega/\sqrt{2}\Big)+1\right],\;\xi_{1,2}(\epsilon)=\exp\left[\mp \frac{n_G(\epsilon)}{2\pi f_0}\right]$,\newline
$\omega_\psi=1$, $C_\psi=\sqrt{\frac{\pi}{2}}f_0^{-1}/2\pi,\;D_\psi=C_\psi e^{\frac{1}{2}\left(4\pi^2f_0^2\right)^{-1}}$.
\\ \hline
Morlet & ``Morlet'' &
$\hat{\psi}(\xi)=e^{-(\xi-2\pi f_0)^2}(1-e^{-2\pi f_0\xi}),\;\xi\in(0,\infty)$,\newline
$\psi(t)=\frac{1}{\sqrt{2\pi}}e^{-t^2/2}e^{i2\pi f_0 t}+O(e^{-(2\pi f_0)^2/2}),\;t\in(-\infty,\infty)$,\newline
$\omega_\psi=2\pi f_0+O(e^{-(2\pi f_0)^2/2}),\;D_\psi=\infty$.
\\ \hline
Bump & ``Bump'' &
$\Delta=0.4f_0^{-1}\leq 1$ (so that $f_0\geq0.4$),\newline
$\hat{\psi}(\xi)=\exp\left(1-\frac{1}{1-\Delta^{-2}(1-\xi)^2}\right),\;\xi\in(1-\Delta,1+\Delta)$,\newline
$\omega_\psi=1,\;D_\psi<\infty$.
\\ \hline
Generalized Morse family & ``Morse-a'' (e.g.\ ``Morse-2.5'') &
$q=30f_0/a$,\newline
$\hat{\psi}(\xi)=B\xi^{q}e^{-\xi^a}=e^{-\xi^a+q\log\xi+\log B},\;\xi\in(0,\infty),\;B\equiv(ea/q)^{q/a}$,\newline
$D_\psi=\frac{\omega_\psi B}{2a}\Gamma((q-1)/a)\;\;(=\infty\mbox{ for }q\leq1)$,\newline
$\omega_\psi=(q/a)^{1/a},\;C_\psi=\frac{B}{2a}\Gamma(q/a)$.
\\ \hline
  \end{tabular}}
\caption{Different wavelet types and their characteristics (if known in analytic form). The names under which these wavelets are implemented in the codes \cite{freecodes} are given in the second column, but one can specify any wavelet function there. The resolution parameter $f_0$ for each wavelet is adjusted in such a way that at $f_0=1$ all of them have similar frequency resolutions (as defined in (\ref{tfreswt}) with $\epsilon_r=0.05$), see Fig.\ \ref{fig:wtres}(a) below. An explanation of the mathematical notation is given in \ref{app:nomenclature}.
}
\label{tab:wavelets}
\end{table*}
\end{center}

From Fig.\ \ref{fig:wtres}(a,b) one can see that the widely used Morlet wavelet has the upper (lower) limit on its time (frequency) resolution. Thus, for this wavelet type, $\Delta\nu_{\min}(2\pi)$ saturates at some level when $f_0$ decreases, and the same happens with $\Delta\tau_{\min}(2\pi)$. As a result, lowering $f_0$ below $\approx0.05$ does not effectively change anything: neither the time, nor the frequency, nor the joint time-frequency resolution. This happens because the Morlet wavelet for $f_0\ll 1$ becomes $\hat{\psi}(\xi)\approx2\pi f_0\xi e^{-\xi^2/2}$ with $\omega_\psi\approx 1$, so that $f_0$ determines only the (unimportant) constant multiplier.

However, the impossibility of achieving high time resolution, as well as the worsening of joint time-frequency resolution for low $f_0$, appears to be a general property of wavelets. It is probably related to the fact that the WT effectively considers the wavelet on a linear scale in time (although rescaled at each frequency), and on a logarithmic scale in frequency (totally independent of time), as can be seen e.g.\ from (\ref{sfreq}) and (\ref{stime}). Thus, the properties of $\hat{\psi}(\xi)$ on both linear and logarithmic scales are important: the former determines the time resolution (since $\psi(t)$ is the inverse FT of $\hat{\psi}(\xi)$, considered on a linear scale), while the latter is responsible for the frequency resolution (see Sec.\ \ref{sec:tfres}).

Consider a few examples. Assuming for simplicity that the wavelet peak frequency $\omega_\psi$ is fixed, it is clear that to decrease the effective spread of the wavelet in time (and therefore increase its time resolution), one needs to increase its spread in frequency. If $\hat{\psi}(\xi)$ has a finite support in terms of $\log\xi$, then one can increase its spread around fixed $\omega_\psi$ only up to some limiting point determined by preserving the admissibility $\hat{\psi}(0)=0$; this is the case of the bump wavelet. On the other hand, if $\hat{\psi}(\xi)$ has an infinite support in terms of $\log\xi$, then increasing its spread will typically lead to increase of the asymmetry of $\hat{\psi}(\xi)$ on a linear scale and its more rapid drop to zero as $\xi\rightarrow 0$. Both of these have negative effects on the decay of $\psi(t)e^{-i\omega_\psi t}$, thus counteracting the desired increase in time resolution; this is the case of e.g.\ the lognormal wavelet. Hence, whether it is possible to achieve high time resolution with the WT, and whether there exists a wavelet form for which the joint time-frequency resolution is relatively unaffected by changing the tradeoff between its time and frequency resolutions, are open questions.

From Fig.\ \ref{fig:wtres} it can be seen that the Morse wavelets have slightly better resolution properties than the other forms (Fig.\ \ref{fig:wtres}(c)). However, with increasing $\epsilon_r$ in (\ref{tfreswt}) (Fig.\ \ref{fig:wtres} corresponds to $\epsilon_r=0.05$), the lognormal wavelet becomes progressively better in terms of the time-frequency resolution $\gamma_{\omega t}$ (not shown), outperforming Morse wavelets at $\epsilon_r\gtrsim 0.1$. It also possesses the many other advantages discussed in detail in Sec.\ \ref{sec:tfrWT}, thus being a preferred choice among those listed in Table \ref{tab:wavelets} (and so the lognormal wavelet is implemented as default in the codes). Nevertheless, it remains possible that some better wavelet form could be constructed that would greatly outperform those considered.

\section*{References}

\bibliographystyle{elsarticle-num}    %
\bibliography{resolutionbib}

\end{document}